\documentclass{article}

\usepackage{arxiv}

\usepackage[utf8]{inputenc} 
\usepackage[T1]{fontenc}    
\usepackage{hyperref}       
\usepackage{url}            
\usepackage{booktabs}       
\usepackage{amsfonts}       
\usepackage{nicefrac}       
\usepackage{microtype}      
\usepackage{graphicx}
\usepackage{doi}
\usepackage{graphicx}
\usepackage{bm}
\usepackage{subcaption}
\usepackage{siunitx}
\usepackage{float}
\usepackage{color}
\usepackage[font=normalsize]{caption}
\usepackage{amsmath}
\usepackage{amssymb}
\usepackage{mathtools}
\DeclarePairedDelimiter\norm{\lVert}{\rVert}

\title{Discovering Efficient Periodic Behaviours in Mechanical Systems via Neural Approximators}

\author{Yannik P. Wotte \thanks{Equal contribution author (order determined by coin flip)}\\
	Robotics and Mechatronics group \\
	University of Twente\\
	Enschede, Netherlands \\
	\texttt{y.p.wotte@utwente.nl} \\
    \And
    Sven C. Dummer$^*$ \\
	Mathematics of Imaging and AI group \\
	University of Twente\\
	Enschede, Netherlands \\
	\texttt{s.c.dummer@utwente.nl} \\
    \And
    Nicolò Botteghi$*$ \\
	Mathematics of Imaging and AI group \\
	University of Twente\\
	Enschede, Netherlands \\
	\texttt{n.botteghi@utwente.nl} \\
	\And
	Christoph Brune \\
	Mathematics of Imaging and AI group \\
	University of Twente\\
	Enschede, Netherlands \\
	\texttt{c.brune@utwente.nl} \\
    \And
	Stefano Stramigioli \\
	Robotics and Mechatronics group\\
	University of Twente \\
	Enschede, Netherlands \\
	\texttt{s.stramigioli@utwente.nl} \\
    \And
	Federico Califano \\
	Robotics and Mechatronics group\\
	University of Twente \\
	Enschede, Netherlands \\
	\texttt{f.califano@utwente.nl} \\
}

\renewcommand{\shorttitle}{\textit{arXiv} Template}

\newtheorem{theorem}{Theorem}[section]
\newtheorem{definition}{Definition}[theorem]
\newtheorem{lemma}[theorem]{Lemma}

\begin{document}
\maketitle

\begin{abstract}
It is well known that conservative mechanical systems exhibit local oscillatory behaviours due to their elastic and gravitational potentials, which completely characterise these periodic motions together with the inertial properties of the system. The classification of these periodic behaviours and their geometric characterisation are in an on-going secular debate, which recently led to the so-called \textit{eigenmanifold} theory. The eigenmanifold characterises nonlinear oscillations as a generalisation of linear eigenspaces. With the motivation of performing periodic tasks efficiently, we use tools coming from this theory to construct an optimization problem aimed at inducing desired closed-loop oscillations through a state feedback law. We solve the constructed optimization problem via gradient-descent methods involving neural networks. Extensive simulations show the validity of the approach.
\end{abstract}


\maketitle

\section{Introduction}

Mechanical systems, such as industrial robots or bio-inspired ones, often need to perform tasks exhibiting a periodic nature, e.g., pick and place or periodic locomotion. The ubiquity of these tasks, as well as the theoretical appeal of understanding and characterising periodic solutions of dynamical systems, made the study of repetitive motions and their control an immensely important branch in the system theoretic community.  

Abstracting from the specific class of mechanical systems and assuming a more general control theoretic perspective, the problem of tracking periodic signals, sometimes referred to as \textit{periodic regulation}, has been intensively tackled with different tools. Without the claim to be exhaustive, we refer to the surveys \cite{doi:10.1080/002071700405905,Wang2009SurveyControl} for an overview, to
\cite{califano2018stability,Astolfi2021NonlinearSystems} (and references therein) for more recent contributions, and to 
\cite{Kasac2008PassiveManipulators} for an application in robotics. 

Contrarily to what is pursued in this work, the mentioned approaches are mostly focused on the design of controllers which implement some steady-state cancellation of the plant dynamics to achieve tracking of specific periodic reference signals. As mentioned in \cite{Bjelonic2021ExperimentalRobot}, these approaches lack a biomimetic perspective in the sense that the design of the periodic regulator is focused on versatility rather than efficiency. In other words, the focus of these approaches is designing a controller that works for a large class of reference signals rather than designing efficient controllers for a smaller class of efficiently stabilisable periodic trajectories.
In \cite{Bjelonic2021ExperimentalRobot}, this efficiency objective is pursued by steering a mechanical system onto natural oscillations of the system itself, which are matched to the mechanical system's physics.

The existence of such periodic oscillations for nonlinear mechanical systems with conservative potentials (usually considered of elastic and gravitational type) is a well-known fact \cite{Rosenberg1966OnFreedom,Shaw1993,Avramov2013}, and the recent theory of \textit{eigenmanifolds} \cite{Albu-Schaffer2020ASystems} attempts at giving a geometric characterisation of these families of oscillations. These oscillations constitute an invariant of the system, i.e., when no dissipative effects or other disturbances are present, a system initialised on such a nonlinear mode would stay there autonomously, with no need of additional inputs. The control theoretic appeal for such structure is immediate once a nonlinear oscillation is understood as a desired periodic behaviour for the closed-loop system, which can vary from achieving energy efficient forms of locomotion, to industrial-like tasks like e.g., pick and place. 

In \cite{Bjelonic2021ExperimentalRobot} the authors successfully stabilised these periodic oscillations defined by eigenmanifold theory, claiming an efficient control design. In fact, a controller able to stabilise a specific invariant oscillation of the system only requires a minimal power consumption, as in principle only the energy to compensate for dissipative effects would be injected by the controller. In conclusion, the underlying biomimetic rational drives the designer in exploiting the natural physics (elastic joints, gravity, inertial parameters) to understand and stabilise an efficiently stabilisable behaviour with minimal energy consumption. We refer to the recent paper \cite{Albu-Schaffer2022WhatRobots} for further elaborations about the connection between efficiency in robotics and the exploitation of natural physics (referred to as "intrinsic dynamics" in that work) present in mechanical systems.

In \cite{Bjelonic2021ExperimentalRobot} the approach was limited to stabilise the open-loop nonlinear modes 
produced by the conservative elastic and gravitational potentials of the underlying mechanical system. Motivated by the fact that natural modes of the open-loop system might not correspond to desired task-specific oscillations, and that mechanical design of a system achieving specific desired oscillations might be very difficult, we introduce a new scheme, which can be seem as an extension of the one in \cite{Bjelonic2021ExperimentalRobot} to account for a broader class of periodic oscillations. In particular, we aim at learning and stabilising a desired oscillation which achieves the fulfillment of some periodic task, which is \textit{close} to the natural mode of the underlying system, but not necessarily coincident. The main contribution of this paper is to present a procedure aimed at finding a potential based state-feedback law which generates desired efficient oscillations in the closed-loop system. In order to do so we cast the control problem into an optimisation framework in which the decision variable is a control potential, approximated by a neural network and updated through gradient descent to minimise a task-dependent performance metric together with a metabolic cost. 
The learned potential uniquely defines a feedback law which generates a closed-loop system exhibiting the desired oscillations. These are then stabilised using an approach similar to \cite{Bjelonic2021ExperimentalRobot,Albu-Schaffer2020ASystems},
where non trivial adaptations have been made to improve the energetic behavior of the control.

Extensive simulations performed on a double pendulum show the validity of the approach.




\subsection{Structure of the paper}
\begin{figure}
    \centering
    \includegraphics[width=\textwidth]{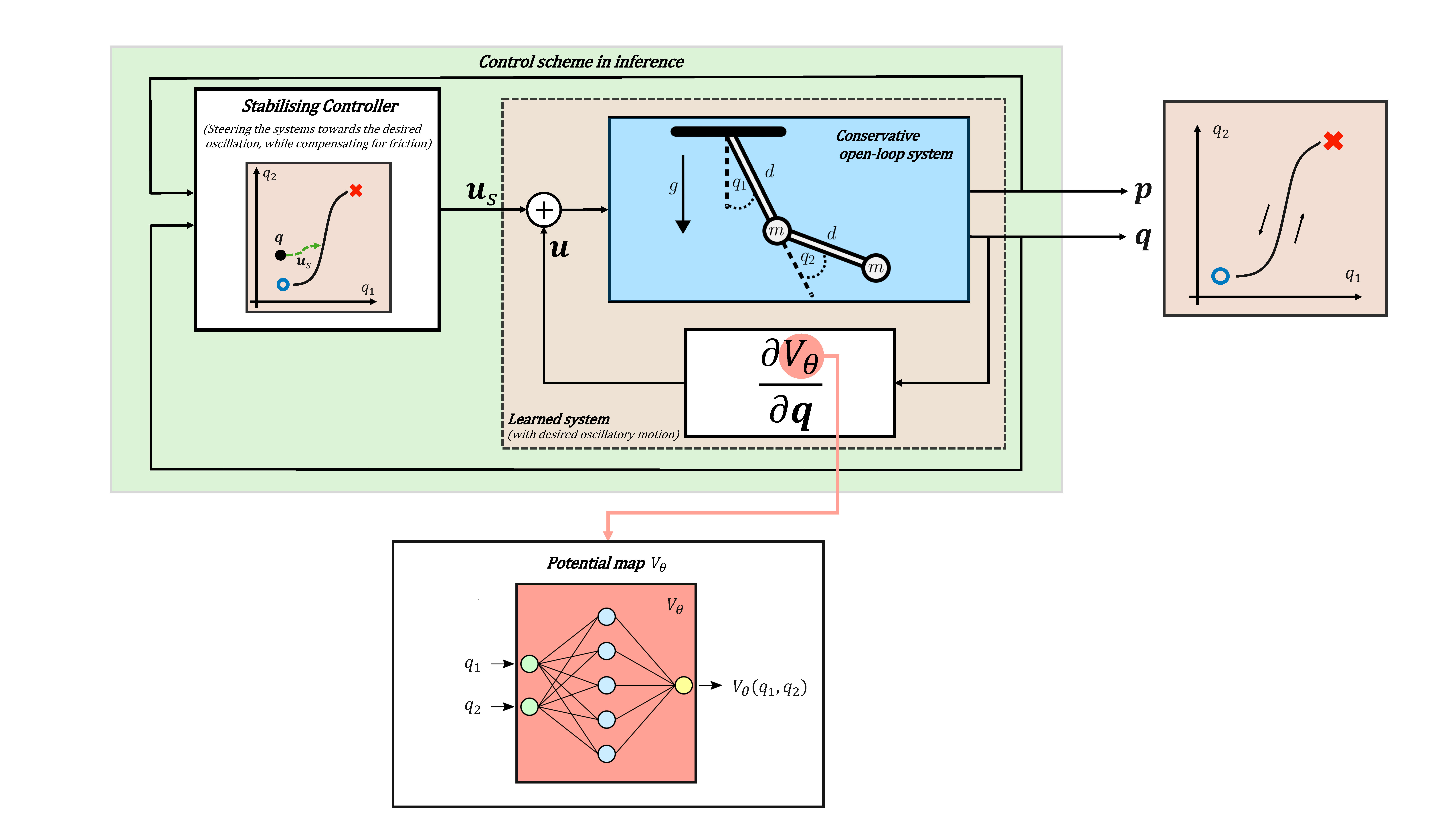}
    \caption{General architecture of the control scheme and synopsis of the paper}\label{fig:generalscheme}
\end{figure}

The structure of the paper is sketched in Fig. \ref{fig:generalscheme}. In Sec. \ref{sec:background}, we give some background material on the Hamiltonian formulation of controlled mechanical systems and on eigenmanifolds. In Sec. \ref{sec:learning}, the main contribution of this work, the optimisation of the control potential is presented and addressed through gradient descent methods involving neural networks as functional approximators. The section is concluded by defining the controller aimed at stabilising the mechanical system on the learned periodic mode and addressing the energetic behavior (in particular passivity) of the resulting closed-loop system. Sec. \ref{sec:sim} contains the simulations and discussions, while Sec. \ref{sec:conc} concludes the paper. The extensive appendices \ref{AppB:additional_results_config2} to \ref{AppE:Implementation} show further results of the proposed optimization.


\section{Background}
\label{sec:background}
\subsection{Hamiltonian formalism for controlled conservative mechanical systems}

In this work we deal with conservative mechanical systems, and we use the Hamiltonian formalism to describe their dynamics. Even if not standard in the eigenmanifold literature, this choice will provide technical advantages in formally presenting some properties of interest. In order to keep the focus on the relevant contributions, in this work we will present all the equations in "standard" coordinates with $\mathbb{R}^{n}$ as the configuration space for an $n$-dimensional mechanical system \footnote{The proper configuration space for a e.g. a double pendulum is the torus $T^2$, rather than $\mathbb{R}^2$ - in this work, the distinction has a negligible impact.}. However, all the concepts can be generalised at a manifold level. The Hamiltonian dynamics (with control) of an $n$-DoF conservative mechanical system with position $\bm{q}(t) \in \mathbb{R}^n$ and momentum $\bm{p}(t) \in \mathbb{R}^n$ is described (hiding time dependencies for lightening notation) by

\begin{equation}
    \begin{aligned}
        \frac{\mathrm{d}}{\mathrm{d}t} \begin{bmatrix}
           \bm{q} \\
           \bm{p}
         \end{bmatrix} & = \begin{bmatrix}
            0 & \bm{I} \\ -\bm{I} & 0
         \end{bmatrix} \nabla H(\bm{p}, \bm{q}) + \begin{bmatrix}
            \bm{0} \\ \bm{I}
         \end{bmatrix} \bm{u}\\
        \begin{bmatrix}
           \bm{q}(0) \\
           \bm{p}(0)
         \end{bmatrix} & = \begin{bmatrix}
           \bm{q}_0 \\
           \bm{p}_0
         \end{bmatrix}
    \end{aligned}
    \label{eq:system_to_consider}
\end{equation}
where $H(\bm{q},\bm{p})=K(\bm{q},\bm{p})+V(\bm{q})$ is the Hamiltonian, i.e., the total mechanical energy of the system. The total mechanical energy $H$ is given by the sum of kinetic energy $K(\bm{q},\bm{p})=\frac{1}{2} \bm{p}^T M^{-1}(\bm{q}) \bm{p}$ (where $M(\bm{q})\in \mathbb{R}^{n \times n}$ is the inertia tensor) and the potential energy $V(\bm{q})$, storing the conservative gravitational and elastic effects. As standard in this formalism, the gradient operator applied to the Hamiltonian is given by $\nabla H(\bm{p}, \bm{q})=\begin{bmatrix}
   \frac{\partial}{\partial {\bm{q}}}H(\bm{p}, \bm{q}) \quad \frac{\partial}{\partial {\bm{p}}}H(\bm{p}, \bm{q})
\end{bmatrix}^T \in \mathbb{R}^{2n}$, and $\bm{I}$ and $\bm{0}$ are the $n$-dimensional identity and zero matrices respectively. We consider an explicit control input $\bm{u}$, representing the generalised forces on the mechanical system collocated to the degrees of freedom defining the position coordinates $\bm{q}$. The usual corollary that the Hamiltonian function is conserved along autonomous evolutions ($\dot{H}=0$ holds along solutions of (\ref{eq:system_to_consider}) with $\bm{u}=0$) will be used in the rest of this work.

\subsection{Eigenmanifolds}\label{sec:eigen_man}
Eigenmanifold theory\cite{Albu-Schaffer2020ASystems} generalises the theory of oscillations present in linear mechanical systems to conservative, intrinsically nonlinear mechanical systems. Here, the essentials of this formalism are presented in its Hamiltonian form.

\begin{definition}\label{Definition:Eigenmode}
An isolated eigenmode $\bm{x}:\mathbb{R}\rightarrow \mathbb{R}^{2n}$ of an autonomous conservative mechanical system (i.e., a system in the form \eqref{eq:system_to_consider} with $\bm{u}=0$) is a trajectory $\bm{x}(t) = (\bm{q}(t),\bm{p}(t))$ with the properties:
\begin{itemize}
    \item $\bm{x}$ is periodic, i.e. $\exists T>0: (\bm{q}(t),\bm{p}(t)) = (\bm{q}(t+T),\bm{p}(t+T))$.
    \item within one period there must be two distinct points with zero momentum, i.e.  $\exists t_1 \neq t_2, t_2 - t_1 < T:\bm{p}(t_1) = 0  \textrm{ and } \bm{p}(t_2) = 0$.
    \item the set $\{\bm{q}(t) | t \in \mathbb{R}\}$ is "line-shaped", i.e., it is homeomorphic to the closed interval $[0,1]\subset\mathbb{R}$.
\end{itemize}
\end{definition}

\begin{figure}[h!]
    \centering
    \includegraphics[width=0.35\textwidth]{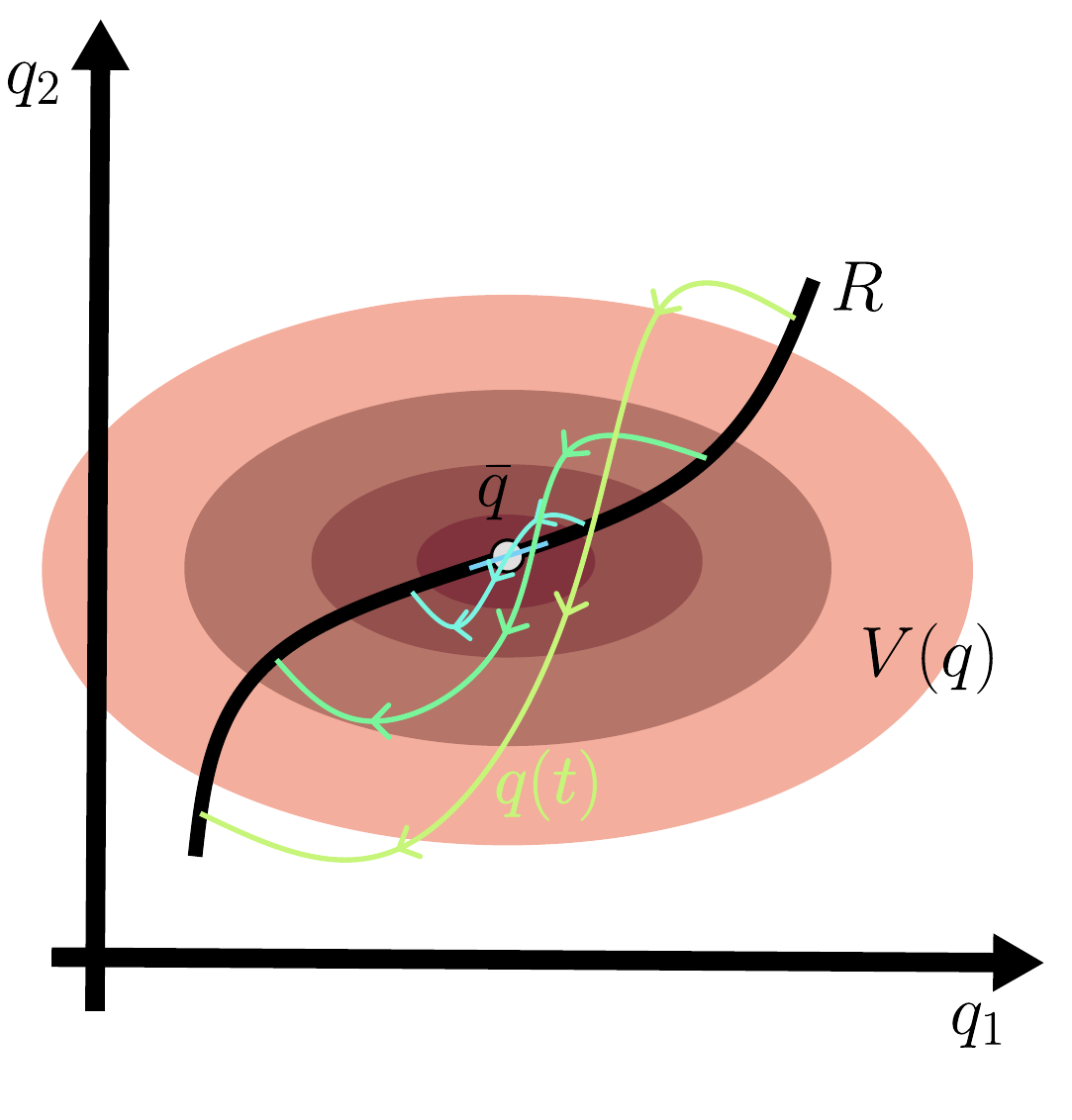}
    \caption{In a mechanical system with potential $V(\bm{q})$ and equilibrium $\bar{\bm{q}}$, an eigenmanifold is a collection of eigenmodes $x(t) = (\bm{q}(t),\bm{p}(t))$, such that the 1-dimensional generator $R$ (which collects particular initial conditions of different eigenmodes) contains the equilibrium $x_{eq} = (\bar{\bm{q}},\bm{0})$ as a limit point. Here, these concepts are depicted after their projection to a 2-dimensional configuration space. }\label{fig:eigenmanifold}
\end{figure}

An eigenmanifold $E\subseteq \mathbb{R}^{2n}$ is then a collection of such modes $\bm{x}$, defined with respect to an isolated, stable equilibrium $\bm{x}_{\textrm{eq}} =  (\bar{\bm{q}},0)$ of the system \eqref{eq:system_to_consider}. Such an equilibrium, which represents the "trivial mode" in the eigenmanifold, exists at a minimum of the potential $\bm{V}(\bm{q})$. 
The additional demand is that the collection $R = \bm{x}_{\textrm{eq}} \cup \{\bm{x}(t) | \bm{p}(t) = 0\}$, called the generator of the eigenmanifold, is a connected, 1-dimensional submanifold \footnote{Usually, generators are defined such that only one point of each mode appears on the generator, i.e., each eigenmanifold has two generators. The distinction is not important for the present work. }, see also Fig. \ref{fig:eigenmanifold}. The generator represents the collection of points which are the extrema of the oscillations of every mode in the eigenmanifold. These modes, for systems in the form (\ref{eq:system_to_consider}), partially characterise the periodic oscillations that a frictionless mechanical system can have. 

It is instructive to think about these modes as the collection of trajectories of (\ref{eq:system_to_consider}) factoring out i) the bounded and non periodic evolutions, whose behavior is commonly referred to as \textit{chaotic} and ii) the periodic evolutions for which no point with $\bm{p}(t)=0$ exists, i.e., those which do not qualify as oscillations. Eigenmanifolds are then of particular interest to factor out such trajectories in nonlinear mechanical systems with $n\geq2$ DoFs (e.g., double pendulum), where chaotic behavior is often present.

\textit{Remark:} Similar to linear oscillations the eigenmodes can often \footnote{The energy-level alone does not always induce a unique ordering of eigenmodes on an eigenmanifold.} be ordered in the eigenmanifold for increasing levels of energy along a mode (starting with the zero energy level corresponding to the trivial mode which is the equilibrium $\bm{x}_{\textrm{eq}}$). Contrarily to the linear case, the eigenmanifold can be bounded (it is not a linear space and it can also not be extended indefinitely for high energy levels), and every mode has in general a different period $T$ (while in the linear case it is constant).

The problem of \textit{existence} and the complete characterisation of eigenmanifolds for conservative mechanical systems is an open problem and is out of the scope of this paper. We refer to \cite{Albu-Schaffer2020ASystems} for the latest developments in this direction. Nevertheless, as main motivation of this work, both experimental and numerical evidence \cite{Wotte2022SufficientType} are confirming that such nonlinear oscillations are structurally present in mechanical systems of any dimension, and can be detected and stabilised, as shown in \cite{DellaSantina2021ExcitingStabilization,Bjelonic2021ExperimentalRobot}. 

To formulate the proposed eigenmanifold optimization method in Section \ref{sec:learning}, we need three lemmas involving conservative mechanical systems which can be verified using the Hamiltonian formulation (\ref{eq:system_to_consider}). In fact, the latter system (with $\bm{u}=0$) is subject to the discrete symmetry $(\bm{q},\bm{p},t) \rightarrow (\bm{q},-\bm{p},-t)$, i.e., if $(\bm{q}(t), \bm{p}(t))$ is a forward in time solution for (\ref{eq:system_to_consider}), then $(\bm{q}(-t), -\bm{p}(-t))$ is likewise a forward in time solution for (\ref{eq:system_to_consider}). The following lemmas, which are proven in \cite{Wotte2022SufficientType}, act as corollaries.
\begin{lemma}
Any trajectory with $\bm{q}(0) = \bm{q}_0$ arbitrary and initial momentum $\bm{p}(0) = 0$ will have the property that $\bm{q}(t) = \bm{q}(-t)$ and $\bm{p}(t) = -\bm{p}(-t)$. 
\label{lemma:symmetry}
\end{lemma}

\begin{lemma}
Any \textit{periodic} trajectory with $\bm{p}(0) = 0$ and period $T$ will have the property that $\bm{p}(T/2) = 0$.
\label{lemma:zero_momentum_half_way}
\end{lemma}

\begin{lemma}
Any trajectory with two distinct points $\bm{q}(t_1)  \neq \bm{q}(t_2)$ such that $\bm{p}(t_1) = 0$ and $\bm{p}(t_2) = 0$ will be periodic, with period $T = 2|(t_1-t_2)|$.
\label{lemma:periodicity_identification}
\end{lemma}

Combining the definition of an eigenmode and the previous lemmas, valid for any conservative mechanical system, the following can be concluded. An \textit{eigenmode} with initial conditions\footnote{Note that due to the definition of eigenmode this choice of initial condition does not induce a loss of generality.} $(\bm{q}(0)=\bm{q}_0,\bm{p}(0)=0)$ has necessarily a period $T=2\bar{t}$, where $\bar{t}$ is the time instant of the other extremum of the oscillation, i.e., $\bm{p}(\bar{t})=0$. Conversely, if a periodic trajectory with properties defined in Lemma \ref{lemma:periodicity_identification} presents a line shaped set $\{\bm{q}(t) | t \in \mathbb{R}\}$, it is necessarily a (possibly isolated) eigenmode. It is worth to note that this condition of line-shapedness was rarely violated in practice.

\section{Related Work: Neural Networks in Dynamical Systems}
The recent developments of artificial intelligence and machine learning has opened the door to new approaches for understanding and controlling dynamical systems by relying on data. In particular, data-driven methods, e.g. neural networks, have often been used as function approximators for learning the dynamics of systems \cite{Brunton2016DiscoveringSystems, StevenLBrunton2022Data-drivenControl, Champion2019Data-drivenEquations, Kalia2021LearningEquations} or for representing control strategies \cite{RichardS.Sutton2018ReinforcementIntroduction, Lillicrap2015ContinuousLearning, Arulkumaran2017DeepSurvey} even in high-dimensional optimal control problems \cite{Ruthotto2020AProblems}. Furthermore, neural networks can be used to approximate the Lyapunov functions in the case of autonomous \cite{Serpen2005EmpiricalNets, PetridisVasilios2006ConstructionFunctions, GabyLyapunov-Net:Approximation, Grune2020ComputingNetworks, RichardsTheSystems} and non-autonomous dynamical systems \cite{Chang2020NeuralControl, MittalNeuralExamples, MohammadKhansari-Zadeh2014LearningMotions, Long1993FeedbackNetworks} for stability and control purposes.  

However, purely data-driven methods often learn physically-inconsistent models that do not respect physical conservation laws. Therefore, the most recent research trends have shifted towards encoding physical principles into neural networks. Examples are hamiltonian, symplectic, and lagrangian neural networks \cite{GreydanusGoogleBrain2019HamiltonianNetworks, CranmerLagrangianNetworks, Zhong2019SymplecticControl} and the physics-informed neural networks \cite{Raissi2019Physics-informedEquations}, aiming at exploiting the best of both worlds, namely the expressive power of nonlinear function approximators with grounded physical knowledge. 

Another important step in this direction has been the introduction of Neural Ordinary Differential Equations (Neural ODEs) \cite{Chen2018NeuralEquations}. The neural ODE framework allows the study of a neural network and its training phase as ODE, opening many possibilities for analysis and understanding of black-box methods.

A closely-related approach to our methods is the work of \cite{massaroli2021optimal}, where a neural ODE is used for learning an optimal passive controller in the port-Hamiltonian framework. The learned controller is composed of a learned potential energy term and a learned damping injection term. However, differently from \cite{massaroli2021optimal} which solve the problem of the stabilization of an inverted pendulum, we focus on a more complex problem, namely the learning of energy-efficient eigenmodes for optimally solving pick and place tasks with a double pendulum. Additionally, instead of learning a damping injection term, we introduce a passive controller injecting only the energy lost by the system due to dissipative elements.

\section{Discovering and stabilizing optimal eigenmodes}
In this work, we want to control trajectories efficiently towards periodic signals that perform some task. To do so, we first need to find a periodic signal that 1) represents the execution of a task and 2) allows for efficient control towards it. We discuss this in Section \ref{sec:learning}. Once such an oscillation is found, we discuss a controller that steers trajectories onto this orbit in Section \ref{ssec:stablizing_control_theory}

\subsection{Discovering optimal eigenmodes via Neural Approximators}
\label{sec:learning}



To find an oscillatory motion that allows for efficient control towards it, we consider eigenmodes (see definition \ref{Definition:Eigenmode}) of the system in \eqref{eq:system_to_consider} where we restrict the input to the gradient of a control potential $\bm{u}=\nabla_q{V_{\bm{\theta}}(\bm{q})}$, so that the closed-loop system will have the form of an autonomous mechanical system \eqref{eq:system_to_consider} ($\bm{u}=0$) with Hamiltonian $H+V_{\bm{\theta}}$. The rationale behind this choice is to modify the system dynamics from \eqref{eq:system_to_consider} as little as possible, avoiding potential cancellation approaches and exploiting the natural physics in the most efficient way. 

The aim is to find the map $V_{\bm{\theta}}: \mathbb{R}^n \rightarrow \mathbb{R}$, such that for a fixed initial condition the resulting motion is an eigenmode of the closed-loop system and minimizes a task-dependent cost term $L_{\text{task}}$. This yields a constrained optimization problem whose decision variable is the map $V_{\bm{\theta}}$. To solve this problem with gradient descent methods, a finite-dimensional parametrisation of $V_{\bm{\theta}}$ is necessary. We denote $\bm{\theta}$ the vector collecting the (finitely many) parameters characterising the map $V_{\bm{\theta}}$, which motivates the notation for the latter. In this work $\bm{\theta}$ will collect the parameters of a neural network, which will be used as functional approximator for $V_{\bm{\theta}}(\bm{q})$. In this perspective, the closed-loop system in optimisation phase becomes a so-called Neural ODE\cite{Chen2018NeuralEquations}.

Summarising the above considerations, the optimisation that we aim to solve is then represented as:
\begin{equation}
\begin{aligned}
\min_{\bm{\theta}} \quad & L_{\text{task}}(\bm{x}) \\
\textrm{s.t.} \quad &
        \frac{\mathrm{d}}{\mathrm{d}t} \begin{bmatrix}
           \bm{q} \\
           \bm{p}
         \end{bmatrix}  = \begin{bmatrix}
            0 & \bm{I} \\ -\bm{I} & 0
         \end{bmatrix} \nabla (H+V_{\bm{\theta}})(\bm{p}, \bm{q}) \quad
        , \begin{bmatrix}
           \bm{p}(0) \\
           \bm{q}(0)
         \end{bmatrix} & = \begin{bmatrix}
           \bm{0} \\
           \bm{q}_0
         \end{bmatrix} \\
  & L_{\textrm{eigen}}(\bm{x}) = 0
\label{eq:optimization_problem_eigmode_formal_fixed_q0_and_fixed_T}
\end{aligned}
\end{equation}

where $L_{\text{task}}(\bm{x})$ is the loss function of the problem, $L_{\text{eigen}}(\bm{x})=0$ represents the constraint forcing the closed-loop trajectory to be an eigenmode, and $\bm{x} = (\bm{q},\bm{p})$. Notice that the choice $\bm{p}_0 = \bm{0}$ happens without loss of generality since we are dealing with periodic orbits, and by Def. \ref{Definition:Eigenmode} an eigenmode is always characterised by $\bm{p}(t) = \bm{0}$ for some $t$. Unless specified otherwise, in this work we assume both the initial position $\bm{q}_0$ of the eigenmode and its period $T$ to be fixed.

In Section \ref{sec:sim}, we solve the optimization problem \eqref{eq:optimization_problem_eigmode_formal_fixed_q0_and_fixed_T} for a pick and place experiment where we move from initial task space position $h(\bm{q}_0)$ (being $h(\bm{q})$ the forward kinematic map) to a desired position $h^*$. In this case we design $L_{\text{task}}(\bm{x})$ as:
\begin{equation}
    L_{\text{task}}(\bm{x}) = \frac{1}{2}\alpha_{\text{task}}\norm{h(\bm{q}\left(\frac{T}{2}\right))  - h^*}_2^2 + \alpha_{\text{eff}} \int_0^T \norm{\bm{u}}_2^2 \mathrm{d}t,
    \label{eq:L_task}
\end{equation}
where $\norm{\cdot}_2$ is the 2-norm, such that the first term promotes the minimisation of the distance between the end-effector position at time $t=\frac{T}{2}$ and the target position $h^*$, and the second term is of metabolic nature and penalises high control efforts $\bm{u}(t)= \nabla V_\theta(\bm{q}(t))$ \footnote{While this cost is well-defined in the given coordinates, we point out the implicit choice of a distance $\norm{h(\bm{q}\left(\frac{T}{2}\right))  - h^*}_2$ on the task space $\mathbb{R}^2$ and a norm $\norm{\bm{u}}_2$ in the vector space $T_q^*\mathbb{R}^2 \cong \mathbb{R}^2$.}. Here, $\alpha_{\text{eff}}$ is a positive scalar balancing the contribution of the two terms, whose effect is analysed in Appendix \ref{AppB1:ablation_study_alphaeff}.

The constraint $L_{\text{eigen}}(\bm{x})=0$ in \eqref{eq:optimization_problem_eigmode_formal_fixed_q0_and_fixed_T} is designed in a way to force the evolution of the closed-loop system to be an oscillation: 
the construction of the function $L_{\textrm{eigen}}(\bm{x})$ is inspired by Lemma's \ref{lemma:symmetry}, \ref{lemma:zero_momentum_half_way}, and \ref{lemma:periodicity_identification}. In particular, given the initial conditions in (\ref{eq:optimization_problem_eigmode_formal_fixed_q0_and_fixed_T}), by Lemma's \ref{lemma:zero_momentum_half_way} and \ref{lemma:periodicity_identification}, it suffices to enforce $\bm{p}\left(\frac{T}{2}\right) = \bm{0}$ to get a periodic trajectory of period $T$. Moreover, given a periodic trajectory with period $T$, Lemma \ref{lemma:symmetry} shows that $\bm{q}(t) = \bm{q}(T - t)$ and $\bm{p}(t) = -\bm{p}(T-t)$. Finally, by periodicity, the trajectory satisfies $\bm{q}(T) = \bm{q}_0$ and $\bm{p}(T) = \bm{0}$. Combining all these observations, we chose the following form for $L_{\textrm{eigen}}$:
\begin{equation}
    L_{\text{eigen}}(\bm{x}(t)) = \lambda_1 \left(\norm{\bm{q}(t) - \bm{q}(T-t)}_{\infty, T}  + \alpha_1 \norm{\bm{p}(t) + \bm{p}(T - t)}_{\infty, T} \right) + \frac{\lambda_2}{2} \norm{\bm{p}\left(\frac{T}{2}\right)}^2_2 
    \label{eq:L_eigen}
\end{equation}
where $\lambda_i\in \mathbb{R}_+$ ($i = 1,$), $\alpha_1 \in \mathbb{R}_+$, and where $\norm{\cdot}_{\infty,T}$ is defined by:
$
    \norm{\bm{y}(\cdot)}_{\infty,T} := \max_{t\in [0, \frac{T}{2}]}(\norm{\bm{y}(t)}_1)
$
with $\bm{y}:[0, \infty) \rightarrow \mathbb{R}^n$.


As an alternative to \eqref{eq:optimization_problem_eigmode_formal_fixed_q0_and_fixed_T}, the eigenmode constraint can be relaxed into a soft one by solving the optimisation:
\begin{equation}
\begin{aligned}
\min_{\bm{\theta}} \quad & L_{\text{task}}(\bm{x}) + \beta L_{\textrm{eigen}}(\bm{x}) \\
\textrm{s.t.} \quad &
        \frac{\mathrm{d}}{\mathrm{d}t} \begin{bmatrix}
           \bm{q}(t) \\
           \bm{p}(t)
         \end{bmatrix}  = \begin{bmatrix}
            0 & \bm{I} \\ -\bm{I} & 0
         \end{bmatrix} \nabla (H+V_{\bm{\theta}})(\bm{p}, \bm{q}) \quad
        , \begin{bmatrix}
           \bm{p}(0) \\
           \bm{q}(0)
         \end{bmatrix} & = \begin{bmatrix}
           \bm{0} \\
           \bm{q}_0
         \end{bmatrix}
\end{aligned}
\label{eq:optimization_problem_eigmode_to_solve}
\end{equation}
with $\beta \in \mathbb{R}^+$ a positive constant.

We stress that even though this version of the optimisation does not present $L_{\textrm{eigen}} = 0$ as a hard constraint, in the moment in which a line shaped periodic trajectory results as a solution of the optimization problem, we are able to assess the learning of an eigenmode with the same confidence as for \eqref{eq:optimization_problem_eigmode_formal_fixed_q0_and_fixed_T} by considering Definition \ref{Definition:Eigenmode}, although we can in principle not ensure that the optimization will result in a periodic orbit.

\textit{Remark:} In the pick and place experiment, we want the end-effector to stop at a specific location $h^* \neq h(\bm{q_0})$ at some arbitrary time $t$. By choosing $t = \frac{T}{2}$ in \eqref{eq:L_task}, the constraint $L_{\textrm{eigen}}(\bm{x}) = 0$ guarantees that the end-effector will actually stop at $h^*$. 

\textit{Remark:} With the lemmas \ref{lemma:symmetry}, \ref{lemma:zero_momentum_half_way}, \ref{lemma:periodicity_identification}, and the definition of an eigenmode in mind, it is easy to check that the trajectory $\bm{x}(t)$ will correspond to an eigenmode in the sense of eigenmanifold theory if and only if $\lambda_2 > 0$, $L_{\textrm{eigen}}(\bm{x}(t)) = 0$ and $\{\bm{q}(t) | t \in \mathbb{R}\}$ is line shaped. As a consequence, the term $\lambda_1$ is strictly speaking redundant, but was found to improve the convergence (together with the specific choice of norms in \eqref{eq:L_eigen}) in the optimisation. In conclusion, if the solution of the optimization problem above yields a line-shaped periodic trajectory, we can conclude that the orbit indeed corresponds to an eigenmode. We furthermore stress the practical scarcity of non line-shaped periodic trajectories, which, to the knowledge and the experience of the authors, have been rarely found in the previously studied cases. 

\subsubsection{Solving the optimisation}

Given the finite-dimensional parametrisation of the map $V_{\bm{\theta}}(\bm{q})$, in this work the optimisation is solved through gradient descent methods, i.e., the optimal parameters $\bm{\theta}$ are found by iterating:
\begin{equation}
    \bm{\theta}_{k+1}=\bm{\theta}_{k}-\eta_k \frac{\partial}{\partial \bm{\theta}}L(\bm{x}(\bm{\theta}))
\end{equation}
where $L(\bm{x})=L_{\text{task}}(\bm{x}) + \beta L_{\textrm{eigen}}(\bm{x})$ is the cost in \eqref{eq:optimization_problem_eigmode_to_solve}.
If $\eta_k$, a positive scalar referred to as \textit{learning rate}, is suitably chosen, and $L(\bm{x})$ is convex, $\bm{\theta}$ converges
to the minimiser of $L(\bm{x})$ as $k \rightarrow \infty$. Although global convergence is no longer guaranteed in the
nonconvex case (which is the case of this work), gradient descent techniques are widely used in practical applications, especially
among the machine learning community, due to their scalability and computational efficiency.

In order to implement the gradient descent procedure, the sensitivity $\frac{\partial}{\partial \bm{\theta}}L(\bm{x}(\bm{\theta}))$ needs to be computed. This is where the so called neural ODE framework, an extension of the continuous depth framework for recurrent neural networks, is used. In particular, the dynamic constraint in \eqref{eq:optimization_problem_eigmode_formal_fixed_q0_and_fixed_T} has the structure of a neural ODE, i.e, an ordinary differential equations parametrised by a neural network $V_{\bm{\theta}}(\bm{q})$ with parameters $\bm{\theta}$. The \textit{training} of this continuous network corresponds to solving the optimisation problem \eqref{eq:optimization_problem_eigmode_formal_fixed_q0_and_fixed_T}. The sensitivities $\frac{\partial}{\partial \bm{\theta}}L(\bm{x}(\bm{\theta}))$ are calculated via the backpropagation method, in particular via automatic differentiation\cite{PaszkeAutomaticPyTorch}, that is commonly used for training neural networks. Utilising the adjoint method \cite{Chen2018NeuralEquations} for computing the exact sensitivies, rather than the approximate ones computed by backpropagation, is an option for future investigation.

\subsection{Stabilising Controller and Analysis of the Closed-Loop System}\label{ssec:stablizing_control_theory}
We formally introduced the optimisation that aims at learning a closed-loop conservative mechanical system exhibiting desired oscillations. In real applications, where dissipative effects and parametric disturbances are present, it is important to design a controller able to robustly stabilise the closed-loop system onto the learned eigenmode. With the motivation of interpreting the learned oscillations as "efficient" (minimizing a certain cost-function), it would furthermore be desirable that the stabilising controller acts in a energetically convenient way (i.e., the control effort is equal to zero on the desired trajectory, and the controller is passive, if no dissipation is present). In other words, the controller should inject the mechanical energy needed to stay on the eigenmode into the system and it should compensate for unavoidable dissipative effects only, resembling a clear biomimetic approach. In \cite{Bjelonic2021ExperimentalRobot} such a controller was successfully implemented to stabilise the (open-loop) eigenmodes of a $7$-DoF KUKA iiwa robot. Here we propose an alternative stabilising controller that is likewise split into an energy-injecting and an eigenmode stabilizing part. Contrary to \cite{Bjelonic2021ExperimentalRobot}, the latter is not allowed to inject energy in this work. The effect of this splitting will simplify the analysis of the controller.

The system with stabilizing feedback $\bm{u}_s:\mathbb{R}^{2n}\rightarrow\mathbb{R}^n$ is of the form 
\begin{equation}\label{eq:sys_plus_ctrl}
        \frac{\mathrm{d}}{\mathrm{d}t} \begin{bmatrix}
           \bm{q}(t) \\
           \bm{p}(t)
         \end{bmatrix}  = \begin{bmatrix}
            0 & \bm{I} \\ -\bm{I} & 0
         \end{bmatrix} \nabla (H+V_{\bm{\theta}})(\bm{p}, \bm{q}) 
         + 
         \begin{bmatrix}
           \bm{0} \\
           \bm{I}
         \end{bmatrix} \bm{u}_s(\bm{q},\bm{p})
\end{equation}
The purpose of this feedback is to stabilize an eigenmode $\bar{\bm{x}}:\mathbb{R}\rightarrow\mathbb{R}^{2n}$ ($\bar{\bm{x}}(t)=(\bar{\bm{q}}(t),\bar{\bm{p}}(t))$), the latter being itself a solution of the learned autonomous system \eqref{eq:optimization_problem_eigmode_to_solve}. To this end, the desired requirements are  
\begin{align} \label{eq:ctrl_properties}
    \lim_{t\rightarrow\infty}\text{dist}(\bm{q}(t),\bar{\bm{q}}(\bar{t})) = 0 \,, \\
    \lim_{t\rightarrow\infty}(\|\bm{p}(t)-\sigma\bar{\bm{p}}(\bar{t})\|) = 0 \,, \\
    \bar{t} = \arg\min_{s}\text{dist}(\bm{q}(t),\bar{\bm{q}}(s)) \,, \label{eq:ctrl_bart}\\
    \sigma = \text{sign}(\bm{p}^T(t) M^{-1}(\bm{q}(t))\bm{\bar{p}}(\bar{t})) \label{eq:ctrl_sigma}\,.
\end{align}
Here $\text{dist}(\bm{a},\bm{b})$ returns the Euclidean distance\footnote{In a differential geometric context, $\text{dist}(\bm{x},\bm{y})$ would implement the geodesic distance depending on a choice of metric tensor and connection, while the second requirement would read $\lim_{t\rightarrow\infty}(\|\bm{p}(t)-\rho^*\bar{\bm{p}}(\bar{t})\|) = 0$, with $\rho^*:T^*_{\bar{\bm{q}}(\bar{t})}\mathcal{M}\rightarrow T^*_{\bm{q}(t)}\mathcal{M}$ implementing the parallel transport of the momentum $\bm{p}$ along the geodesic $\rho$ from $\bm{q}(t)$ to $\bar{\bm{q}}(\bar{t})$. Here, instead, $\rho^*$ is chosen to be the identity map in the given coordinate system.} 
of points $\bm{a},\bm{b}\in \mathbb{R}^{n}$. 

Intuitively speaking, $\bar{t}$ in equation \eqref{eq:ctrl_bart} is the parameter at which the desired trajectory $\bar{\bm{q}}$ is closest to the current position $\bm{q}(t)$. In practice, $\bar{t}$ is implemented as a function $\bar{t}:\mathbb{R}^{n}\rightarrow \mathbb{R}$ that takes as input $\bm{q}\in \mathbb{R}^n$. Although $\bar{\bm{q}}(\bar{t})$ is uniquely determined, $\bar{\bm{p}}(\bar{t})$ is only determined up to a sign for an eigenmode, which is chosen according to the sign function $\sigma:\mathbb{R}^{2n}\rightarrow \{-1,0,1\} $ in equation \eqref{eq:ctrl_sigma} to be aligned with the current system momentum $\bm{p}(t)$.

The choice is made to split the control
\begin{equation}\label{eq:ctrl_full}
    \bm{u}_s = \bm{u}_E + \bm{u}_M
\end{equation}
into an energy-controlling feedback $\bm{u}_E$ (cf. \cite{Folkertsma2014Power-continuousSystems}) and an eigenmode stabilizing feedback $\bm{u}_M$. For an analogous control splitting see \cite{Bjelonic2021ExperimentalRobot,DellaSantina2021ExcitingStabilization}. 

\subsubsection{Energy-controlling feedback}
The energy-controlling feedback $\bm{u}_E$ steers the system's energy $E=H+V_{\bm{\theta}}$ towards a desired energy $\bar{E} = E(\bar{\bm{q}}(0),\bar{\bm{p}}(0))$. 
The form chosen is 
\begin{equation}\label{eq:u_E}
     \bm{u}_E = \alpha_E(\bar{E}-E) \hat{\bm{p}} \,, 
\end{equation}
with $\alpha_E \in \mathbb{R}^+$ a positive control gain and the normalized momentum \footnote{To avoid numerical issues in practice, $\hat{\bm{p}}$ is chosen as 0 when $\bm{p}^T M(\bm{q}) \bm{p} = 0$. }
\begin{equation}
    \hat{\bm{p}} = \frac{1}{\sqrt{\bm{p}^T M^{-1}(\bm{q}) \bm{p}}} {\bm{p}}
\end{equation}
Since $\dot{\bm{q}}=M^{-1}(\bm{q}) \bm{p}$, it holds that the mechanical power $\bm{u}_E^T \dot{\bm{q}}$ injected by the energy controller is given by
\begin{equation}
    \bm{u}_E^T \dot{\bm{q}} = \alpha_E (\bar{E} - E)\sqrt{\bm{p}^T M^{-1}(\bm{q}) \bm{p}} \,.
\end{equation} 


\subsubsection{Eigenmode stabilizing feedback}
The eigenmode stabilizing feedback $\bm{u}_M$ is defined as

\begin{equation}\label{eq:u_M}
    \bm{u}_M = \alpha_M \pi_{\bm{p}}(\sigma \bar{\bm{p}}(\bar{t}))\,,   
\end{equation}
where $\alpha_M\in\mathbb{R}^+$ is the positive control gain. Furthermore, $\sigma(\bm{q},\bm{p})\in \{-1,0,1\}$ and $\bar{t}(\bm{q}) \in \mathbb{R}$ are as defined in \eqref{eq:ctrl_sigma} and \eqref{eq:ctrl_bart}, respectively. $\bar{\bm{p}}:\mathbb{R}\rightarrow \mathbb{R}^n$ is the momentum component of the desired eigenmode. Last, $\pi_{\bm{p}}$ is the projection defined by
\begin{equation}
    \pi_{\bm{p}}(\bm{X}) := \bm{X} - \frac{\bm{p}^T M^{-1}(\bm{q})\bm{X}}{ \bm{p}^T M^{-1}(\bm{q}) \bm{p}} \bm{p}\,.
\end{equation}
This projection is such that 

\begin{equation}\label{eq:noPower}
    \bm{u}_M^T \dot{\bm{q}} = 0\,,
\end{equation}
which means that $\bm{u}_M$ cannot change the energy content of the system, and thus cannot interfere with the control-task of $\bm{u}_E$.

\textit{Remark:} This is a D-type controller analogous to \cite{Bjelonic2021ExperimentalRobot,DellaSantina2021ExcitingStabilization}, with the only adaptation being that the energy injection is restricted (compare e.g. \cite{Duindam2004PassiveDynamics}). The controller of the form \eqref{eq:u_M} follows from 

\begin{equation}
    \bm{u}_M = \alpha_M \pi_{\bm{p}}(\sigma\bar{\bm{p}}(\bar{t})-\bm{p})\,,
\end{equation}
by using the property of the projection that $\pi_{\bm{p}}(\bm{p}) = 0$.

\subsubsection{Stability}
We first investigate the energetic behavior of the combined controller $\bm{u}_s=\bm{u}_E+\bm{u}_M$, and investigate the stability of the trajectory afterwards.
The energy injected by the controller is equal to the mechanical power $\bm{u}_s^T\dot{q}$:  
\begin{equation}
    \dot{E} = \frac{\partial E}{\partial \bm{q}} \dot{\bm{q}} + \frac{\partial E}{\partial \bm{p}} \dot{\bm{p}}  = \frac{\partial E}{\partial \bm{p}} \bm{u} = \dot{\bm{q}}^T \bm{u}_s = \bm{u}_s^T \dot{\bm{q}} \,.
\end{equation}
Here, the second equality holds because the system without feedback $\bm{u}_s$ conserves $E$, while the third equality follows from the definition of momentum $M(\bm{q})^{-1}p = \dot{\bm{q}}$.

Combining the expressions shows that 
\begin{equation}
    \dot{E} = \bm{u}_s^T \dot{\bm{q}} = \alpha_{E} (\bar{E} - E)\sqrt{\bm{p}^T M^{-1}(\bm{q}) \bm{p}}\,.
\end{equation}
Hence, the energy converges to the desired energy level $\bar{E}$ 
almost always, i.e. as long as $\bm{u}_E \neq 0$, and otherwise $E$ is constant. 
Moreover, as the combined actions $\bm{u}_E$ and $\bm{u}_M$ vanish only on the desired mode, we get highly efficient control behavior as highlighted in \cite{DellaSantina2021ExcitingStabilization, Bjelonic2021ExperimentalRobot}. 

However, the above does not prove either global or local stability. This work restricts itself to a guarantee of local stability, which can be obtained by evaluating the cycle multipliers of the stabilized periodic orbit. Let $\Psi_t(\bm{x}(0)) := \bm{x}(t)$ define the flow of the dynamic system \eqref{eq:ctrl_full}, then cycle multipliers can be defined as the ratio of partial derivatives \footnote{Typically, cycle multipliers are defined as the eigenvalues of $\frac{\partial}{\partial \bm{x}} (\Psi_T(\bm{x}) - \bm{x}))_{|\bm{x}=\bm{x}_0}$. The authors found the alternative definition to be more robust, numerically.} 
\begin{equation}
    \frac{\frac{\partial}{\partial x^i} \text{dist}(\Psi_T(\bm{x}),\bar{\bm{x}}(\bar{t}))_{|\bm{x}=\bm{x}_0}}{\frac{\partial}{\partial x^i} \text{dist}(\bm{x},\bar{\bm{x}}(\bar{t}))_{|\bm{x}=\bm{x}_0}} \,. 
\end{equation}
Here, $x^i$ denotes the $i$-th component of $\bm{x}$ and $\bm{x}_0$ is a starting point of the stabilized periodic orbit, while $\bar{t}$ and $\bar{\bm{x}}$ are as defined in and above equation \eqref{eq:ctrl_bart}. 
As will be shown along the result section, if these cycle multipliers have absolute values smaller than $1$, the periodic orbit is stable.

\section{Simulations}
\label{sec:sim}
In this section, we perform numerical experiments for the case of a double pendulum. More precisely, we consider a pick and place experiment where we want the end-effector of the double pendulum to move between two points in an oscillatory fashion. To achieve this, an optimal eigenmode is learned via the optimization strategy in Section \ref{sec:learning}. For our numerical experiments, we solve the optimization problem in \eqref{eq:optimization_problem_eigmode_to_solve} with loss functions given in Equations \eqref{eq:L_task} and \eqref{eq:L_eigen}. Subsequently, we stabilize the eigenmode using the control strategy in Section \ref{ssec:stablizing_control_theory}.

\subsection{Double Pendulum Model}

\begin{figure}[h!]
    \centering
    \includegraphics[width=0.3\textwidth]{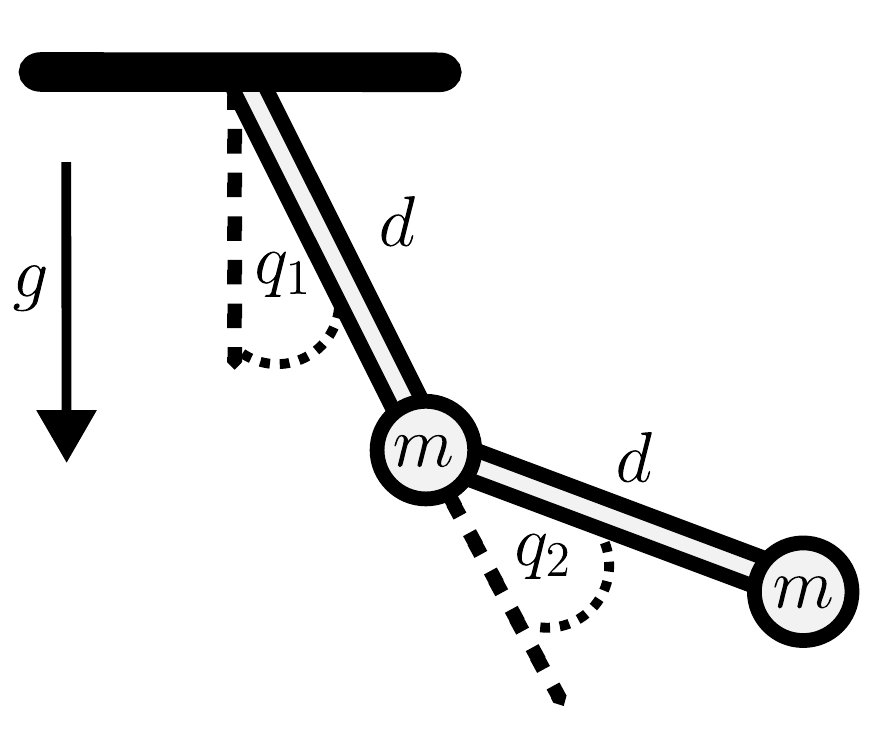}
    \caption{Double Pendulum corresponding to Equations \eqref{eq:DP_Hamiltonian}-\eqref{eq:DP_Potential}.}\label{fig:DP_Conventions}
\end{figure}

The double pendulum is one of the simplest mechanical systems with non-trivial eigenmanifolds (see also \cite{Wotte2022SufficientType}).
The presented double pendulum is under the influence of gravity and has a linear spring at the second joint. The equations of motion correspond to the conventions shown in Figure \ref{fig:DP_Conventions}. They are fully determined by \eqref{eq:system_to_consider} and the Hamiltonian $H:\mathbb{R}^2\times \mathbb{R^2}\rightarrow\mathbb{R}$ given as in Equations \eqref{eq:DP_Hamiltonian} to \eqref{eq:DP_Potential}.

\begin{align} 
    & H(\bm{q},\bm{p}) = \bm{p}^T M^{-1}(\bm{q}) \bm{p} + V(\bm{q}) \,, \label{eq:DP_Hamiltonian} \\
    & M((q_1,q_2)) = m d^2 \begin{bmatrix} 
             (3+2\cos(q_2)) & \cos(q_2)+1 \\ \cos(q_2)+1 & 1  \end{bmatrix} \,, \label{eq:DP_Mass} \\
    & V((q_1,q_2)) = V_{\bm{\theta}}((q_1,q_2)) - mdg (2 \cos(q_1) + \cos(q_1+q_2)) + k (q_2 -\pi/2)^2 \label{eq:DP_Potential} \,.
\end{align}
Here, $V_{\bm{\theta}}(\bm{q}):\mathbb{R}^2\rightarrow\mathbb{R}$ is a potential function, that will be constructed as a neural net with parameters $\bm{\theta} \in \mathbb{R}^m$. The actual equations of motion are reported for completeness in Appendix \ref{AppA:YipYip}. 

\subsection{Results}

\subsubsection{Learning Eigenmodes}
\begin{figure}[h!]
\centering
\hfill
\begin{subfigure}{.33\textwidth}
  \centering
  \includegraphics[width=1.0\linewidth]{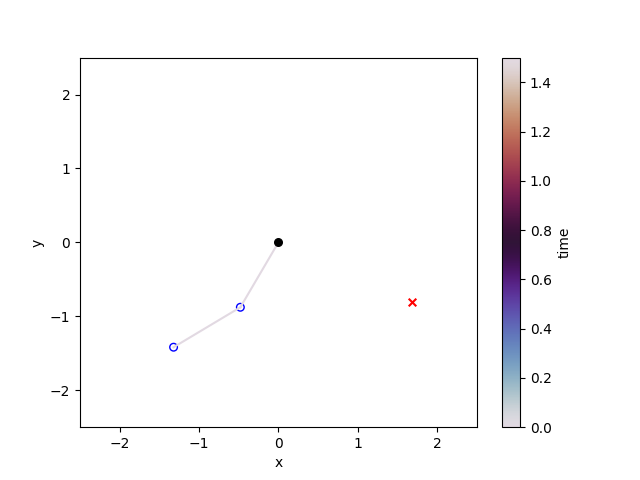}
  \caption{t=\SI{0.0}{\second}}
  \label{fig:sub1}
\end{subfigure}%
\hfill
\begin{subfigure}{.33\textwidth}
  \centering
  \includegraphics[width=1.0\linewidth]{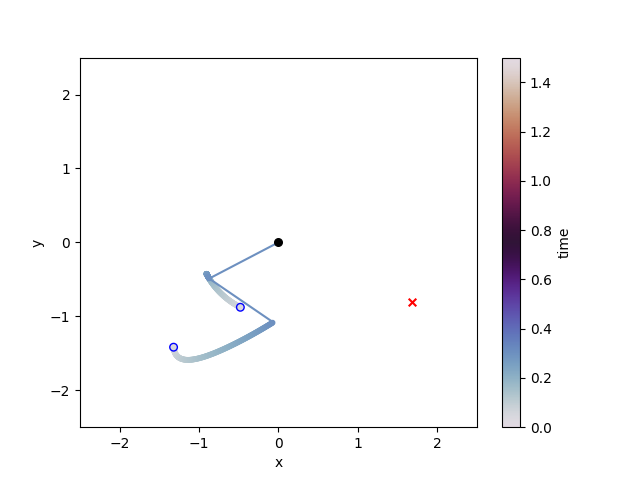}
  \caption{t=\SI{0.332}{\second}}
  \label{fig:sub3}
\end{subfigure}
\hfill
\begin{subfigure}{.33\textwidth}
  \centering
  \includegraphics[width=1.0\linewidth]{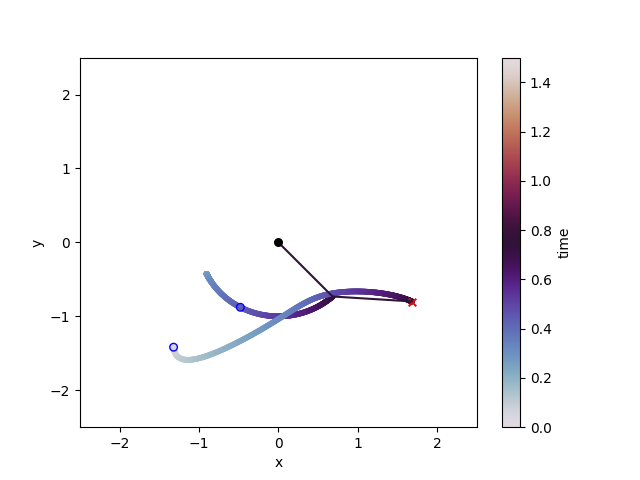}
  \caption{t=\SI{0.83}{\second}}
  \label{fig:sub6}
\end{subfigure}%
\hfill
\begin{subfigure}{.33\textwidth}
  \centering
  \includegraphics[width=1.0\linewidth]{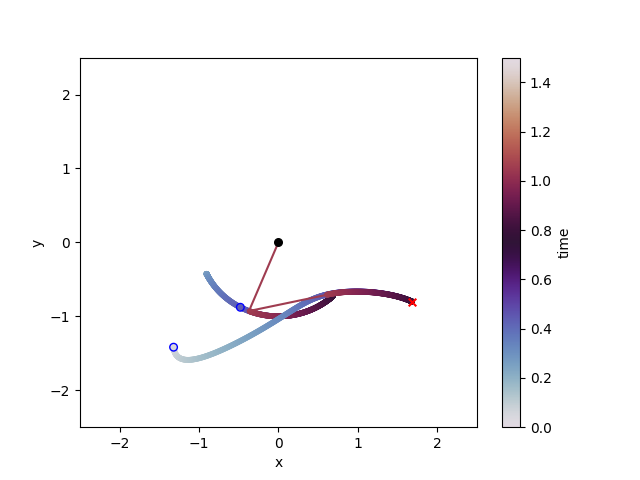}
  \caption{t=\SI{1.162}{\second}}
  \label{fig:sub8}
\end{subfigure}
\begin{subfigure}{.33\textwidth}
  \centering
  \includegraphics[width=1.0\linewidth]{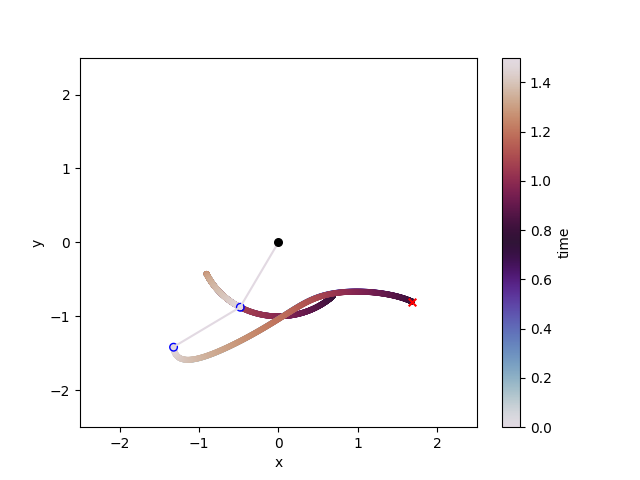}
  \caption{T=\SI{1.5}{\second}}
  \label{fig:sub10}
\end{subfigure}
\caption{Learned eigenmode at different time steps. The blue circles represent the initial position of the joints of the pendulum, while the red cross represents the end-effector target used for computing the first term of $L_{\text{task}}(\bm{x})$ in \eqref{eq:L_task}.}
\label{fig:snapshots_trajectory}
\end{figure}

\begin{figure}[h!]
     \centering
          \begin{subfigure}[b]{0.49\textwidth}
         \centering
         \includegraphics[width=0.9\textwidth]{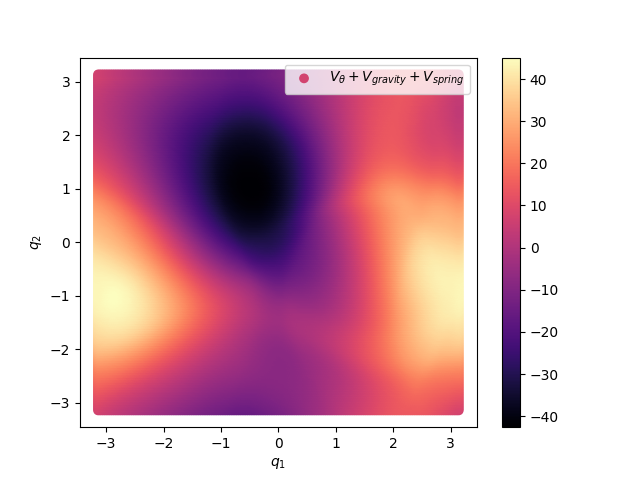}
         \caption{Overall potential $V_{\bm{\theta}}+V_{\text{spring}}+V_{\text{gravity}}$.}
         \label{fig:overall_potential_T=1.5}
     \end{subfigure}
     \hfill
     \begin{subfigure}[b]{0.49\textwidth}
         \centering
         \includegraphics[width=0.9\textwidth]{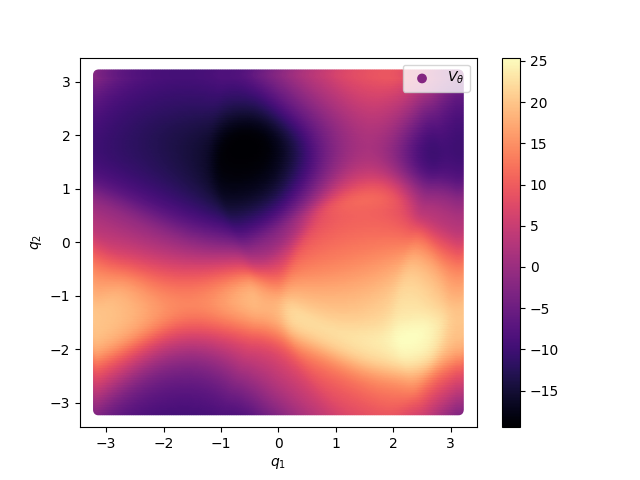}
         \caption{Learned potential $V_{\bm{\theta}}$.}
         \label{fig:learned_potential_T=1.5}
     \end{subfigure}
        \caption{Potentials over $\bm{q}\in[-\pi, \pi]$.}
        \label{fig:three graphs}
\end{figure}

\begin{figure}[h!]
     \centering
     \begin{subfigure}[b]{0.33\textwidth}
         \centering
         \includegraphics[width=1.0\textwidth]{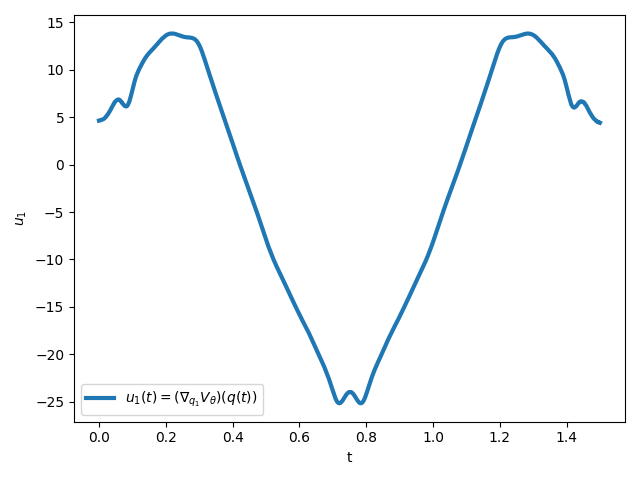}
         \caption{First component of control $u(t)$}
         \label{fig:gradV_q1_T=1.5}
     \end{subfigure}
     \hfill
     \begin{subfigure}[b]{0.33\textwidth}
         \centering
         \includegraphics[width=1.0\textwidth]{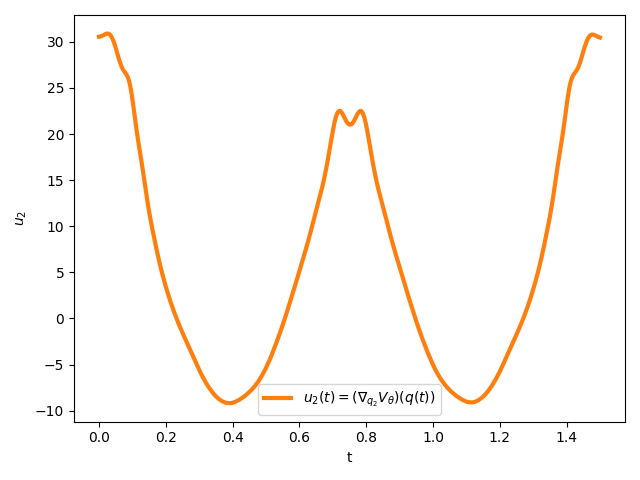}
         \caption{Second component of $u(t)$.}
         \label{fig:gradV_q2_T=1.5}
     \end{subfigure}
          \hfill
     \begin{subfigure}[b]{0.33\textwidth}
         \centering
         \includegraphics[width=1.0\textwidth]{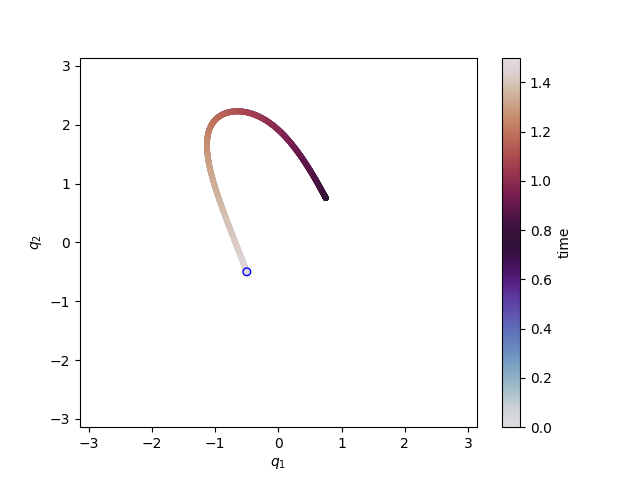}
         \caption{Trajectory in the configuration space.}
         \label{fig:traj_q_lineshaped}
     \end{subfigure}
        \caption{Control inputs over time (Figure \ref{fig:gradV_q1_T=1.5} and \ref{fig:gradV_q2_T=1.5}), and trajectory in the configuration space (Figure \ref{fig:traj_q_lineshaped}).}
        \label{fig:grad_learned_potential}
\end{figure}

In Figure \ref{fig:snapshots_trajectory}, we visualise the trajectory of the inner closed-loop conservative system (see Figure \ref{fig:generalscheme}) at different time instants after the training of the learned potential, via the optimisation procedure described in Section \ref{sec:learning}, for 500 epochs and a given period $T=\SI{1.5}{\second}$. The results are obtained with the set of loss function hyperparameters reported in Table \ref{tab:loss_hyper}\footnote{The complete list of hyperparameters is shown in Table \ref{tab_app:hyper_list}.}. 
\begin{table}[h!]
\centering
\begin{tabular}{||c | c ||} 
\hline
 Hyperparameter & Value \\ [0.5ex] 
 \hline\hline
 $\alpha_{\text{task}}$ & 10\\
 \hline
 $\alpha_{\text{eff}}$ & 0.0001 \\
 \hline
 $\alpha_{\text{task}}$ & 10\\
 \hline
$\lambda_1$ & 0.05\\
 \hline
  $\alpha_1$ & 0.0005\\
 \hline
 $\lambda_2$ & 0.95\\
  \hline 
 $\beta$ & 1\\
 \hline
\end{tabular}
\caption{Loss function hyperparameters used in the experiments.}
\label{tab:loss_hyper}
\end{table}
The potential $V_{\bm{\theta}}$ (see Figure \ref{fig:learned_potential_T=1.5}) is capable of shaping the systems potential (\ref{fig:overall_potential_T=1.5}), such that the trajectory of the system is an energy-efficient eigenmode of the desired period $T$. Additionally in Figure \ref{fig:grad_learned_potential}, we depict the control inputs $\bm{u}=\nabla_{\bm{q}}V_{\bm{\theta}}(\bm{q}(t))$ inducing the desired periodic behaviour, and the trajectory in the configuration space (Figure \ref{fig:traj_q_lineshaped}) from which it is possible to notice the line-shaped property of the eigenmode described in Def. \ref{Definition:Eigenmode}. 


\subsubsection{Stabilization of the Learned Eigenmode} \label{sec:stabilization_results}
\begin{figure}[h!]
     \centering
     \begin{subfigure}[b]{0.33\textwidth}
         \centering
         \includegraphics[width=0.9\textwidth]{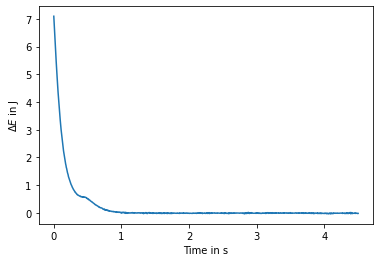}
         \caption{Closed loop energy error.}
         \label{fig:Energy_b=0}
     \end{subfigure}
     \hfill
     \begin{subfigure}[b]{0.33\textwidth}
         \centering
         \includegraphics[width=0.9\textwidth]{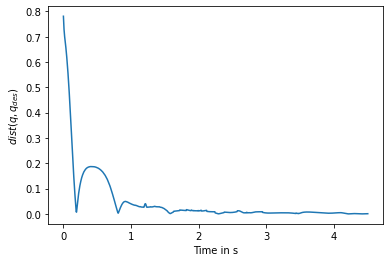}
         \caption{$\textrm{dist}(q(t),\bar{q}(\bar t))$.}
         \label{fig:Distance_b=0}
     \end{subfigure}
     \hfill
     \begin{subfigure}[b]{0.33\textwidth}
         \centering
         \includegraphics[width=0.9\textwidth]{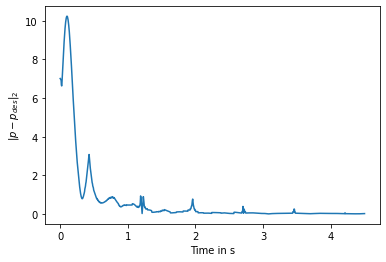}
         \caption{$\|p(t) - \bar{p}(\bar{t})\|_2 $.}
         \label{fig:Momentum_Discrepancy_b=0}
     \end{subfigure}
     \begin{subfigure}[b]{0.33\textwidth}
         \centering
         \includegraphics[width=0.9\textwidth]{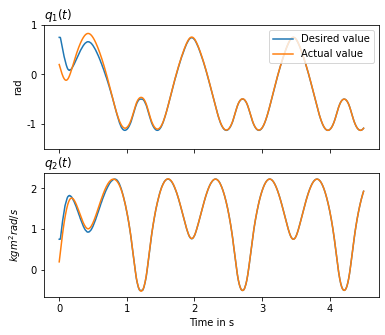}
         \caption{$(q_1(t),q_2(t))$}
         \label{fig:Pos_b=0}
     \end{subfigure}
     \begin{subfigure}[b]{0.33\textwidth}
         \centering
         \includegraphics[width=0.9\textwidth]{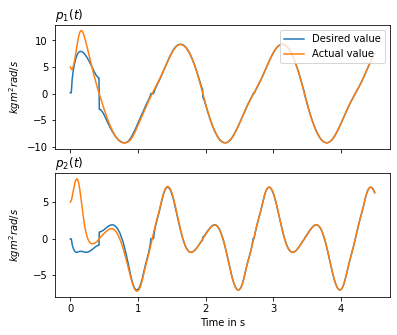}
         \caption{$(p_1(t),p_2(t))$}
         \label{fig:Mom_b=0}
     \end{subfigure}
     \begin{subfigure}[b]{0.33\textwidth}
         \centering
         \includegraphics[width=0.9\textwidth]{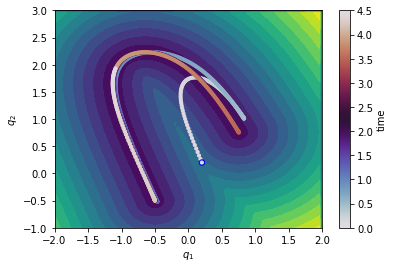}
         \caption{Trajectory in configuration space and level sets of $\textrm{dist}(q,\bar{q})$.}
         \label{fig:Conf_Space_b=0}
     \end{subfigure}
        \caption{Various features of the stabilised system with learned potential as in Figure \ref{fig:three graphs}, stabilizing the mode shown in Figure \ref{fig:snapshots_trajectory} with gains $\alpha_M =10$, $\alpha_E = 1$. The starting condition is $\bm{q}(0) = (0.2,0.2)$, $\bm{p}(0) = (5,5)$, shown here over three periods of oscillation.}
        \label{fig:ctrl_b=0}
\end{figure}

Figure \ref{fig:ctrl_b=0} shows the results of applying the control structure introduced in Section \ref{ssec:stablizing_control_theory} to the learned trajectory shown in Figure \ref{fig:snapshots_trajectory}, for coefficients $\alpha_M = 10$ and $\alpha_E = 1$. 

In the example, we use the starting condition $\bm{q} = (0.2,0.2)$, $\bm{p} = (5,5)$. In particular, Figure \ref{fig:Pos_b=0} and Figure \ref{fig:Mom_b=0} show the development of $\bm{q}(t)$ and $\bm{p}(t)$ over time, which approach the desired $\bar{\bm{q}}(\bar{t})$ and $\bar{\bm{p}}(\bar{t})$ (see Section \ref{ssec:stablizing_control_theory} for their definition). Figure \ref{fig:Energy_b=0} shows the energy $H+V_\theta$ of the closed loop system, which approaches the constant energy level of the learned mode. Figures \ref{fig:Distance_b=0} and \ref{fig:Momentum_Discrepancy_b=0} show the distance of the trajectory from the desired trajectory in position and momentum space respectively (i.e. $\|\bm{q}(t) - \bar{\bm{q}}(\bar{t})\|_2$ and $\|\bm{p}(t) - \bar{\bm{p}}(\bar{t})\|_2$), in both cases approaching 0. The cycle multipliers of the closed loop system are less than 1: for this example, it was found that they are bounded by 0.5, which guarantees that the learned periodic orbit is locally stable. 

\begin{figure}[h!]
     \centering
     \begin{subfigure}[b]{0.33\textwidth}
         \centering
         \includegraphics[width=0.9\textwidth]{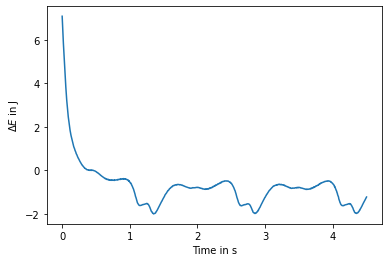}
         \caption{Closed loop energy error.}
         \label{fig:Energy_b=0.1}
     \end{subfigure}
     \hfill
     \begin{subfigure}[b]{0.33\textwidth}
         \centering
         \includegraphics[width=0.9\textwidth]{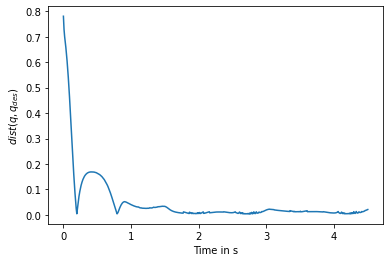}
         \caption{$\textrm{dist}(q(t),\bar{q}(\bar t))$.}
         \label{fig:Distance_b=0.1}
     \end{subfigure}
     \hfill
     \begin{subfigure}[b]{0.33\textwidth}
         \centering
         \includegraphics[width=0.9\textwidth]{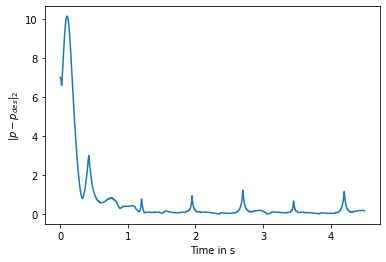}
         \caption{$\|p(t) - \bar{p}(\bar{t})\|_2 $.}
         \label{fig:Momentum_Discrepancy_b=0.1}
     \end{subfigure}
     \begin{subfigure}[b]{0.33\textwidth}
         \centering
         \includegraphics[width=0.9\textwidth]{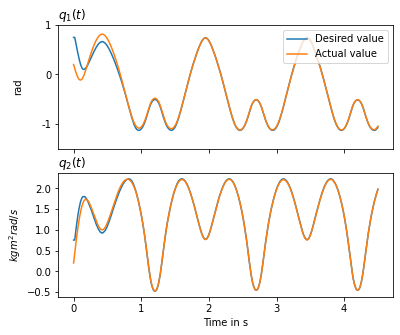}
         \caption{$(q_1(t),q_2(t))$}
         \label{fig:Pos_b=0.1}
     \end{subfigure}
     \begin{subfigure}[b]{0.33\textwidth}
         \centering
         \includegraphics[width=0.9\textwidth]{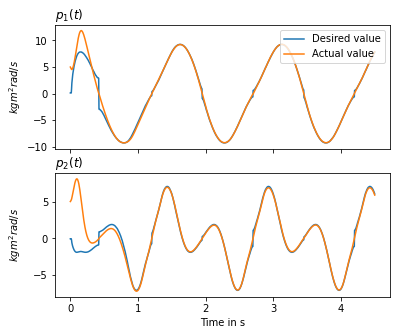}
         \caption{$(p_1(t),p_2(t))$}
         \label{fig:Mom_b=0.1}
     \end{subfigure}
     \begin{subfigure}[b]{0.33\textwidth}
         \centering
         \includegraphics[width=0.9\textwidth]{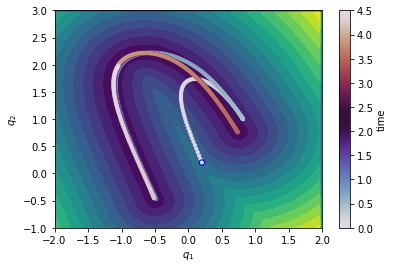}
         \caption{Trajectory in configuration space and level sets of $\textrm{dist}(q,\bar{q})$.}
         \label{fig:Conf_Space_b=0.1}
     \end{subfigure}
        \caption{Various features of the stabilized system shown in \ref{fig:ctrl_b=0}, but including damping linear in system velocity with damping coefficient $b=0.1$}
        \label{fig:ctrl_b=0.1}
\end{figure}

\begin{figure}[h!]
     \centering
     \begin{subfigure}[b]{0.33\textwidth}
         \centering
         \includegraphics[width=0.9\textwidth]{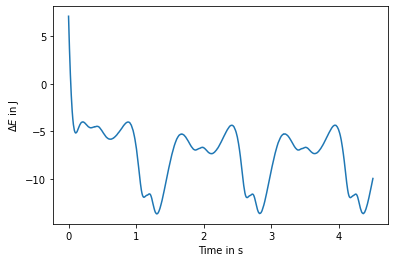}
         \caption{Closed loop energy error.}
         \label{fig:Energy_b=1}
     \end{subfigure}
     \hfill
     \begin{subfigure}[b]{0.33\textwidth}
         \centering
         \includegraphics[width=0.9\textwidth]{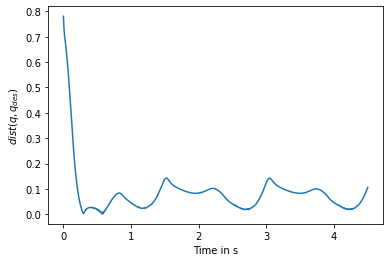}
         \caption{$\textrm{dist}(q(t),\bar{q}(\bar t))$.}
         \label{fig:Distance_b=1}
     \end{subfigure}
     \hfill
     \begin{subfigure}[b]{0.33\textwidth}
         \centering
         \includegraphics[width=0.9\textwidth]{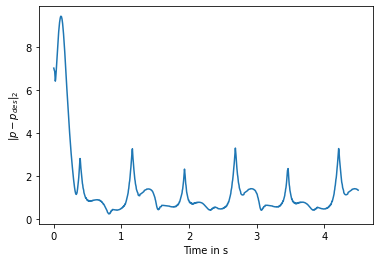}
         \caption{$\|p(t) - \bar{p}(\bar{t})\|_2 $.}
         \label{fig:Momentum_Discrepancy_b=1}
     \end{subfigure}
     \begin{subfigure}[b]{0.33\textwidth}
         \centering
         \includegraphics[width=0.9\textwidth]{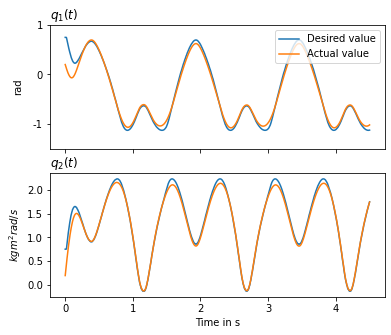}
         \caption{$(q_1(t),q_2(t))$}
         \label{fig:Pos_b=1}
     \end{subfigure}
     \begin{subfigure}[b]{0.33\textwidth}
         \centering
         \includegraphics[width=0.9\textwidth]{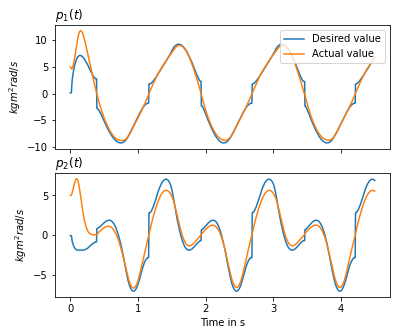}
         \caption{$(p_1(t),p_2(t))$}
         \label{fig:Mom_b=1}
     \end{subfigure}
     \begin{subfigure}[b]{0.33\textwidth}
         \centering
         \includegraphics[width=0.9\textwidth]{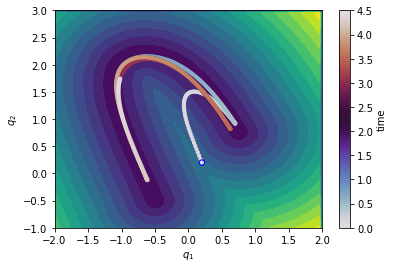}
         \caption{Trajectory in configuration space and level sets of $\textrm{dist}(q,\bar{q})$.}
         \label{fig:Conf_Space_b=1}
     \end{subfigure}
        \caption{Various features of the stabilized system shown in \ref{fig:ctrl_b=0}, but including damping linear in system velocity with damping coefficient $b=1$}
        \label{fig:ctrl_b=1}
\end{figure}

To observe the robustness of the controller in the presence of damping, viscous damping is introduced. With $b$ the damping coefficient, the input $\bm{u}_s$ in \eqref{eq:sys_plus_ctrl} is adapted to read 
\begin{equation}
    \bm{u}_s = \bm{u}_E + \bm{u}_M - b M^{-1}(q) p \,,
\end{equation}
which corresponds to velocity dependent damping. The cases $b = 0.1$ and $b = 1$ are shown in Figures \ref{fig:Conf_Space_b=0.1} and \ref{fig:Conf_Space_b=1}, respectively. Notably, the damping causes the energy shown in Figures \ref{fig:Energy_b=0.1} and \ref{fig:Energy_b=1} to continue to fluctuate in the eventual periodic evolution, about a value lower than the desired energy. It is worth noting that the systems remain close to the desired mode, even for such large cases of damping. However, it should be considered to adapt the energy controlling term $\bm{u}_E$ to compensate for  damping more accurately, as was done e.g. in \cite{Bjelonic2021ExperimentalRobot} for a particular case of damping that was, among others, linear in velocity.

\subsubsection{Additional Results}
To strengthen the numerical contribution, we include additional results and ablation studies in Appendices. In particular, in Appendix \ref{AppA:YipYip}, we show the equations of motion for the double pendulum used in our simulations, while in Appendix \ref{AppB:additional_results_config2}, we study the effect of varying the regularisation coefficient $\alpha_{\text{eff}}$ and the period $T$ on the resulting eigenmode and control inputs. In Appendix \ref{AppC:additional_results_config1}, we apply the method with different initial and target positions, periods, and regularisation coefficients, and in Appedix \ref{AppD:learnableT}, we show a more advanced version of the optimisation problem in \eqref{eq:optimization_problem_eigmode_formal_fixed_q0_and_fixed_T}, where we learn the potential $V_{\bm{\theta}}$ jointly with the period $T$. Eventually, Appendix \ref{AppE:Implementation} reports the implementation details for reproducing our experiments.

\section{Conclusions and Future Work}
\label{sec:conc}

In this paper we present a procedure aiming at shaping desired periodic oscillations for mechanical systems. In particular, using tools from eigenmanifold theory and neural networks as function approximators, a state feedback law is learned in such a way to produce a closed-loop system exhibiting a desired periodic motion. This is done by minimising the effort of the learned control law and exploiting at best the natural physical properties of the underlying open-loop systems, characterised by its inertia and its conservative potentials. A stabilising controller able to steer the system along the learned oscillation in presence of parametric disturbances is presented. Extensive simulations show the validity of the approach. 

Concerning future developments, besides an experimental validation of the scheme, the proposed approach opens the way to co-design of the mechanical system along the desired periodic task. In fact, by constraining the search space of the learned potential to a form which can be reproduced mechanically with e.g., nonlinear springs at the joints, it would be possible to exploit the described learning procedure as a preliminary phase for a mechanical design which would produce a mechanical system achieving the desired behaviour in an open-loop fashion.

\bibliographystyle{unsrt}
\bibliography{references}

\begin{thebibliography}{10}

\bibitem{doi:10.1080/002071700405905}
Richard~W Longman.
\newblock {Iterative learning control and repetitive control for engineering
  practice}.
\newblock {\em International Journal of Control}, 73(10):930--954, 2000.

\bibitem{Wang2009SurveyControl}
Youqing Wang, Furong Gao, and Francis~J. Doyle.
\newblock {Survey on iterative learning control, repetitive control, and
  run-to-run control}.
\newblock {\em Journal of Process Control}, 19(10):1589--1600, 2009.

\bibitem{califano2018stability}
Federico Califano, Michelangelo Bin, Alessandro Macchelli, and Claudio
  Melchiorri.
\newblock {Stability analysis of nonlinear repetitive control schemes}.
\newblock {\em IEEE control systems letters}, 2(4):773--778, 2018.

\bibitem{Astolfi2021NonlinearSystems}
Daniele Astolfi, Laurent Praly, Lorenzo Marconi, and Mines Paristech.
\newblock {Nonlinear Robust Periodic Output Regulation of Minimum Phase
  Systems}.
\newblock 2021.

\bibitem{Kasac2008PassiveManipulators}
Josip Kasac, Branko Novakovic, Dubravko Majetic, and Danko Brezak.
\newblock {Passive finite-dimensional repetitive control of robot
  manipulators}.
\newblock {\em IEEE Transactions on Control Systems Technology},
  16(3):570--576, 2008.

\bibitem{Bjelonic2021ExperimentalRobot}
Filip Bjelonic, Arne Sachtler, Alin Albu-sch, and Cosimo~Della Santina.
\newblock {Experimental Closed-Loop Excitation of Nonlinear Normal Modes on an
  Elastic Industrial Robot}.
\newblock pages 1--8, 2021.

\bibitem{Rosenberg1966OnFreedom}
Reinhard~M. Rosenberg.
\newblock {On nonlinear vibrations of systems with many degrees of freedom}.
\newblock 1966.

\bibitem{Shaw1993}
S.~W. Shaw and C.~Pierre.
\newblock {Normal Modes for Non-Linear Vibratory Systems}.
\newblock {\em Journal of Sound and Vibration}, 164(1):85--124, 6 1993.

\bibitem{Avramov2013}
Konstantin~V. Avramov and Yuri~V. Mikhlin.
\newblock {Review of applications of Nonlinear Normal Modes for Vibrating
  Mechanical Systems}, 3 2013.

\bibitem{Albu-Schaffer2020ASystems}
Alin Albu-Sch{\"{a}}ffer and Cosimo Della~Santina.
\newblock {A review on nonlinear modes in conservative mechanical systems}.
\newblock {\em Annual Reviews in Control}, 50:49--71, 2020.

\bibitem{Albu-Schaffer2022WhatRobots}
Alin Albu-Sch{\"{a}}ffer and Arne Sachtler.
\newblock {What Can Algebraic Topology and Differential Geometry Teach Us About
  Intrinsic Dynamics and Global Behavior of Robots?}
\newblock 2022.

\bibitem{Wotte2022SufficientType}
Yannik Wotte, Arne Sachtler, Alin Albu-Sch{\"{a}}ffer, and Cosimo
  Della~Santina.
\newblock {Sufficient conditions for an eigenmanifold to be of the extended
  Rosenberg type}.
\newblock 2022.

\bibitem{DellaSantina2021ExcitingStabilization}
Cosimo~Della Santina and Alin Albu-Schaeffer.
\newblock {Exciting Efficient Oscillations in Nonlinear Mechanical Systems
  Through Eigenmanifold Stabilization}.
\newblock {\em IEEE Control Systems Letters}, 5(6):1916--1921, 2021.

\bibitem{Brunton2016DiscoveringSystems}
Steven~L. Brunton, Joshua~L. Proctor, J.~Nathan Kutz, and William Bialek.
\newblock {Discovering governing equations from data by sparse identification
  of nonlinear dynamical systems}.
\newblock {\em Proceedings of the National Academy of Sciences of the United
  States of America}, 113(15):3932--3937, 4 2016.

\bibitem{StevenLBrunton2022Data-drivenControl}
{Steven L Brunton} and {J Nathan Kutz}.
\newblock {\em {Data-driven science and engineering: Machine learning,
  dynamical systems, and control}}.
\newblock 2022.

\bibitem{Champion2019Data-drivenEquations}
Kathleen Champion, Bethany Lusch, J.~Nathan~Kutz, and Steven~L. Brunton.
\newblock {Data-driven discovery of coordinates and governing equations}.
\newblock {\em Proceedings of the National Academy of Sciences of the United
  States of America}, 116(45):22445--22451, 11 2019.

\bibitem{Kalia2021LearningEquations}
Manu Kalia, Steven~L. Brunton, Hil G.~E. Meijer, Christoph Brune, and J.~Nathan
  Kutz.
\newblock {Learning normal form autoencoders for data-driven discovery of
  universal,parameter-dependent governing equations}.
\newblock 6 2021.

\bibitem{RichardS.Sutton2018ReinforcementIntroduction}
{Richard S. Sutton} and {Andrew G. Barto}.
\newblock {Reinforcement Learning: An Introduction}, 2018.

\bibitem{Lillicrap2015ContinuousLearning}
Timothy~P. Lillicrap, Jonathan~J. Hunt, Alexander Pritzel, Nicolas Heess, Tom
  Erez, Yuval Tassa, David Silver, and Daan Wierstra.
\newblock {Continuous control with deep reinforcement learning}.
\newblock 9 2015.

\bibitem{Arulkumaran2017DeepSurvey}
Kai Arulkumaran, Marc~Peter Deisenroth, Miles Brundage, and Anil~Anthony
  Bharath.
\newblock {Deep Reinforcement Learning: A Brief Survey}.
\newblock {\em IEEE Signal Processing Magazine}, 34(6):26--38, 11 2017.

\bibitem{Ruthotto2020AProblems}
Lars Ruthotto, Stanley~J. Osher, Wuchen Li, Levon Nurbekyan, and Samy~Wu Fung.
\newblock {A machine learning framework for solving high-dimensional mean field
  game and mean field control problems}.
\newblock {\em Proceedings of the National Academy of Sciences of the United
  States of America}, 117(17):9183--9193, 4 2020.

\bibitem{Serpen2005EmpiricalNets}
Gursel Serpen.
\newblock {Empirical Approximation for Lyapunov Functions with Artiicial Neural
  Nets}, 2005.

\bibitem{PetridisVasilios2006ConstructionFunctions}
{Petridis Vasilios} and {Petridis Stavros}.
\newblock {Construction of Neural Network Based Lyapunov Functions}, 2006.

\bibitem{GabyLyapunov-Net:Approximation}
Nathan Gaby, Fumin Zhang, and Xiaojing Ye.
\newblock {Lyapunov-Net: A Deep Neural Network Architecture for Lyapunov
  Function Approximation}.

\bibitem{Grune2020ComputingNetworks}
Lars Gr{\"{u}}ne.
\newblock {Computing Lyapunov functions using deep neural networks}.
\newblock 2020.

\bibitem{RichardsTheSystems}
Spencer~M Richards, Felix Berkenkamp, and Andreas Krause.
\newblock {The Lyapunov Neural Network: Adaptive Stability Certification for
  Safe Learning of Dynamical Systems}.

\bibitem{Chang2020NeuralControl}
Ya~Chien Chang, Nima Roohi, and Sicun Gao.
\newblock {Neural Lyapunov Control}.
\newblock {\em Advances in Neural Information Processing Systems}, 32, 5 2020.

\bibitem{MittalNeuralExamples}
Mayank Mittal, Marco Gallieri, Alessio Quaglino, Sina~Mirrazavi Salehian, and
  Jan Koutn{\'{i}}k.
\newblock {Neural Lyapunov Model Predictive Control: Learning Safe Global
  Controllers from Sub-optimal Examples}.

\bibitem{MohammadKhansari-Zadeh2014LearningMotions}
S.~Mohammad Khansari-Zadeh and Aude Billard.
\newblock {Learning control Lyapunov function to ensure stability of dynamical
  system-based robot reaching motions}.
\newblock {\em Robotics and Autonomous Systems}, 62(6):752--765, 2014.

\bibitem{Long1993FeedbackNetworks}
Y.~Long and {Bayoumi M.M.}
\newblock {Feedback Stabilization: Control Lyapunov Functions Modelled by
  Neural Networks}, 1993.

\bibitem{GreydanusGoogleBrain2019HamiltonianNetworks}
Sam Greydanus Google~Brain, Misko Dzamba~PetCube, and Jason Yosinski.
\newblock {Hamiltonian Neural Networks}.
\newblock {\em Advances in Neural Information Processing Systems}, 32, 2019.

\bibitem{CranmerLagrangianNetworks}
Miles Cranmer, , Sam Greydanus, , Stephan Hoyer, , Peter Battaglia, , David
  Spergel, , and Shirley Ho.
\newblock {Lagrangian Neural Networks}.

\bibitem{Zhong2019SymplecticControl}
Yaofeng~Desmond Zhong, Biswadip Dey, and Amit Chakraborty.
\newblock {Symplectic ODE-Net: Learning Hamiltonian Dynamics with Control}.
\newblock 9 2019.

\bibitem{Raissi2019Physics-informedEquations}
M.~Raissi, P.~Perdikaris, and G.~E. Karniadakis.
\newblock {Physics-informed neural networks: A deep learning framework for
  solving forward and inverse problems involving nonlinear partial differential
  equations}.
\newblock {\em Journal of Computational Physics}, 378:686--707, 2 2019.

\bibitem{Chen2018NeuralEquations}
Ricky T.~Q. Chen, Yulia Rubanova, Jesse Bettencourt, and David~K. Duvenaud.
\newblock {Neural Ordinary Differential Equations}.
\newblock {\em Advances in Neural Information Processing Systems}, 31, 2018.

\bibitem{massaroli2021optimal}
Stefano Massaroli, Michael Poli, Federico Califano, Jinkyoo Park, Atsushi
  Yamashita, and Hajime Asama.
\newblock {Optimal Energy Shaping via Neural Approximators}.
\newblock {\em SIAM Journal on Applied Dynamical Systems}, 21(3):2126--2147,
  2022.

\bibitem{PaszkeAutomaticPyTorch}
Adam Paszke, Sam Gross, Soumith Chintala, Gregory Chanan, Edward Yang,
  Zachary~Devito Facebook, A~I Research, Zeming Lin, Alban Desmaison, Luca
  Antiga, Orobix Srl, and Adam Lerer.
\newblock {Automatic differentiation in PyTorch}.

\bibitem{Folkertsma2014Power-continuousSystems}
Gerrit~A. Folkertsma, Arjan~J. Van Der~Schaft, and Stefano Stramigioli.
\newblock {Power-continuous synchronisation of oscillators: A novel,
  energy-free way to synchronise dynamical systems}.
\newblock In {\em Proceedings - IEEE International Conference on Robotics and
  Automation}, pages 1493--1498. Institute of Electrical and Electronics
  Engineers Inc., 9 2014.

\bibitem{Duindam2004PassiveDynamics}
Vincent Duindam, Stefano Stramigioli, and Jacquelien~M.A. Scherpen.
\newblock {Passive compensation of nonlinear robot dynamics}.
\newblock {\em IEEE Transactions on Robotics and Automation}, 20(3):480--487,
  2004.

\bibitem{PoliTorchDyn:PyTorch}
Michael Poli, Stefano Massaroli, Atsushi Yamashita, Hajime Asama, and Stefano
  Ermon.
\newblock {TorchDyn: Implicit Models and Neural Numerical Methods in PyTorch}.

\bibitem{Paszke2019PyTorch:Library}
Adam Paszke, Sam Gross, Francisco Massa, Adam Lerer, James Bradbury~Google,
  Gregory Chanan, Trevor Killeen, Zeming Lin, Natalia Gimelshein, Luca Antiga,
  Alban Desmaison, Andreas~Köpf Xamla, Edward Yang, Zach Devito, Martin
  Raison~Nabla, Alykhan Tejani, Sasank Chilamkurthy, Qure Ai, Benoit Steiner,
  Lu~Fang Facebook, Junjie~Bai Facebook, and Soumith Chintala.
\newblock {PyTorch: An Imperative Style, High-Performance Deep Learning
  Library}.
\newblock 2019.

\bibitem{KingmaADAM:OPTIMIZATION}
Diederik~P Kingma and Jimmy Lei~Ba.
\newblock {ADAM: A METHOD FOR STOCHASTIC OPTIMIZATION}.

\end{thebibliography}
\clearpage
\appendix 

\section{Equations of motion for double pendulum}\label{AppA:YipYip}
\begin{align}
    \dot{\bm{q}} = & M(\bm{q})^{-1} \bm{p} \\
    \dot{\bm{p}} = & C(\bm{q},\bm{p})\bm{p}
            - \frac{\partial }{ \partial \bm{q}} V_{\bm{\theta}}(\bm(q)) - \begin{bmatrix} mdg (2\sin(q_1) + \sin(q_1+q_2)) \\  
             mdg\sin(q_1+q_2) + k(\pi - 2q_2)\end{bmatrix} 
\end{align}

where the inverse mass matrix $M^{-1}(\bm{q})$ and Coriolis terms $C(\bm{q},\bm{p})\bm{p}$ are given as
\begin{equation}
    M^{-1}((q_1,q_2)) = \frac{1}{m d^2(3+2\cos(q_2) - (\cos(q_2)+1)^2)} \begin{bmatrix} 
             1 & -\cos(q_2)-1 \\ -\cos(q_2)-1 & (3+2\cos(q_2))  \end{bmatrix} \,,
\end{equation}

\begin{equation}
    C(\bm{q},\bm{p})\bm{p} = \frac{\sin(q_2)}{2d^2m(1+\sin(q_2)^2)^2}\begin{bmatrix} 0 \\ 
    2\cos(q_2)p_1^2 - (5+4\cos(q_2)+\cos(2q_2))p_1p_2 + (5+6\cos(q_2)+\cos(2q_2))p_2^2 \\
    
    \end{bmatrix} \,.
\end{equation} 

\clearpage
\section{Eigenmodes for different penalisation of the control effort and periods}\label{AppB:additional_results_config2}

\subsection{Effect of the Control Effort Penalty}\label{AppB1:ablation_study_alphaeff}
\begin{figure}[h!]
     \centering
     \begin{subfigure}[b]{0.33\textwidth}
         \centering
         \includegraphics[width=1.0\textwidth]{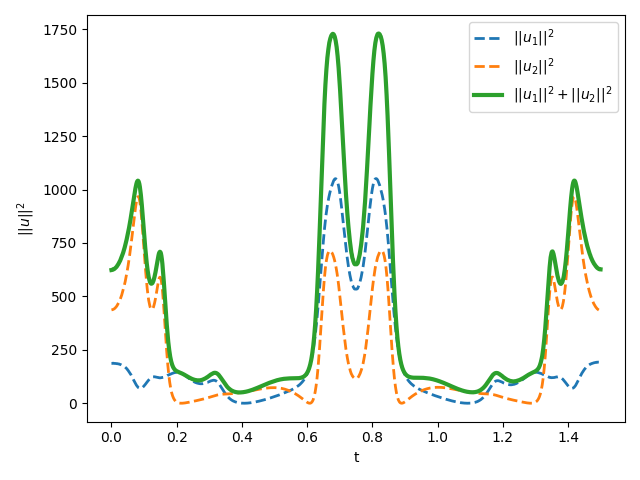}
         \caption{Control effort penalty $\alpha_{\text{eff}}=0.0$.}
         \label{fig:no_reg}
     \end{subfigure}
     \hfill
     \begin{subfigure}[b]{0.33\textwidth}
         \centering
         \includegraphics[width=1.0\textwidth]{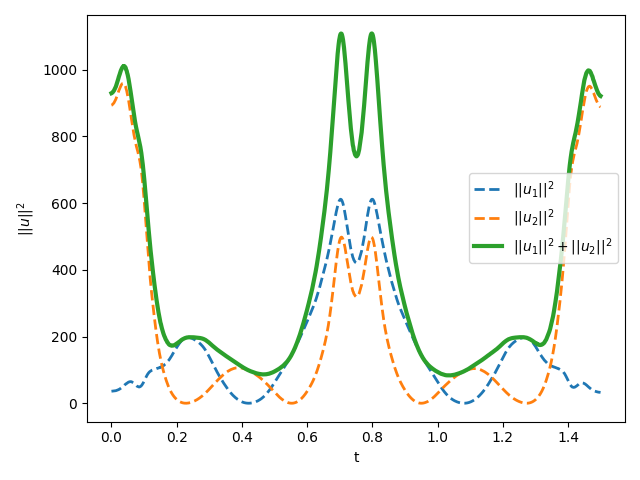}
         \caption{Control effort penalty $\alpha_{\text{eff}}=0.00001$.}
         \label{fig:reg=0.00001}
     \end{subfigure}
     \hfill
     \begin{subfigure}[b]{0.33\textwidth}
         \centering
         \includegraphics[width=1.0\textwidth]{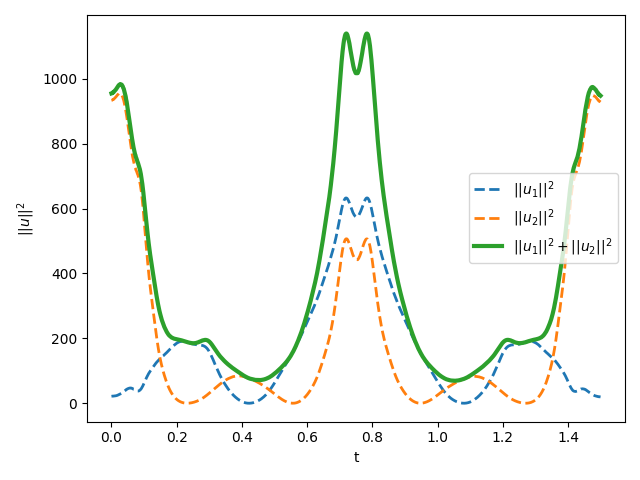}
         \caption{Control effort penalty $\alpha_{\text{eff}}=0.0001$.}
         \label{fig:reg=0.0001}
     \end{subfigure}
     \begin{subfigure}[b]{0.33\textwidth}
         \centering
         \includegraphics[width=1.0\textwidth]{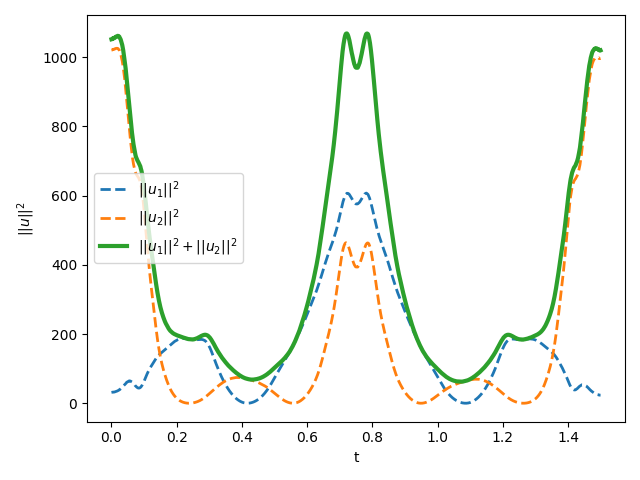}
         \caption{Control effort penalty $\alpha_{\text{eff}}=0.001$.}
         \label{fig:reg=0.001}
     \end{subfigure}
     \begin{subfigure}[b]{0.33\textwidth}
         \centering
         \includegraphics[width=1.0\textwidth]{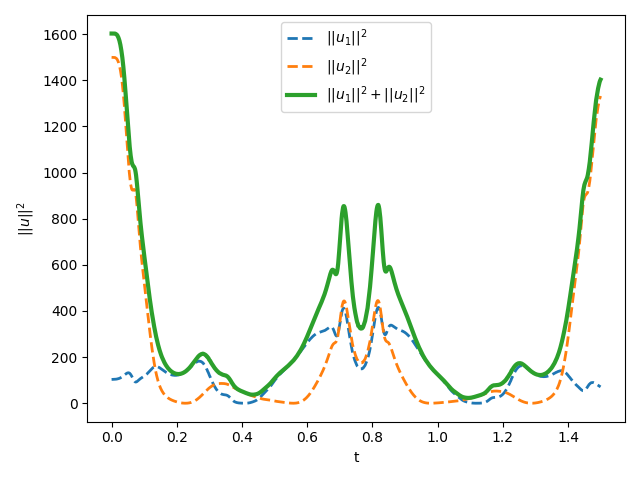}
         \caption{Control effort penalty $\alpha_{\text{eff}}=0.01$.}
         \label{fig:reg=0.01}
     \end{subfigure}
        \caption{Control effort squared for different control effort penalty coefficients.}
        \label{fig:control_effort_penalizaiton}
\end{figure}

In Figure \ref{fig:control_effort_penalizaiton}, we show the squared control effort $||u_1||^2 + ||u_2||^2$ derived from the gradient of the learned potential $V_{\bm{\theta}}$ for different values of the regularization coefficient $\alpha_{\text{eff}} \in \{0.0, 0.00001, 0.0001, 0.001, 0.01\}$. We used this grid-search experiment to find a suitable value for $\alpha_{\text{eff}}$.

\begin{figure}[h!]
     \centering
     \begin{subfigure}[b]{0.33\textwidth}
         \centering
         \includegraphics[width=1.0\textwidth]{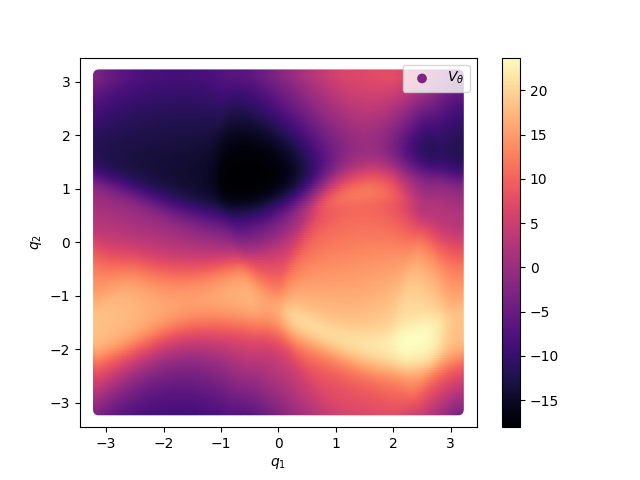}
         \caption{Control effort penalty $\alpha_{\text{eff}}=0.0$.}
         \label{fig:no_reg}
     \end{subfigure}
     \hfill
     \begin{subfigure}[b]{0.33\textwidth}
         \centering
         \includegraphics[width=1.0\textwidth]{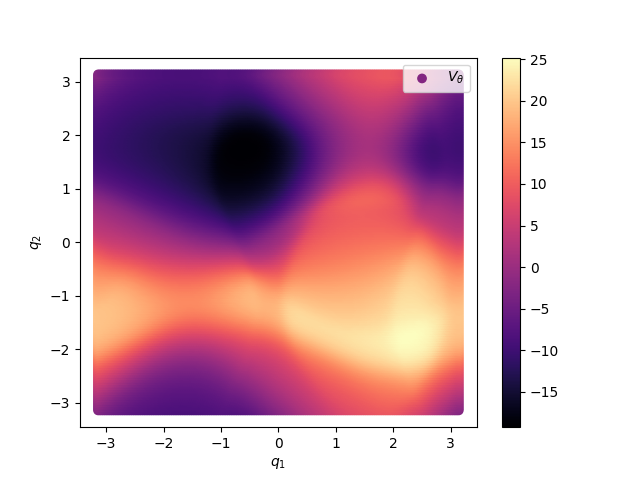}
         \caption{Control effort penalty $\alpha_{\text{eff}}=0.00001$.}
         \label{fig:reg=0.00001}
     \end{subfigure}
     \hfill
     \begin{subfigure}[b]{0.33\textwidth}
         \centering
         \includegraphics[width=1.0\textwidth]{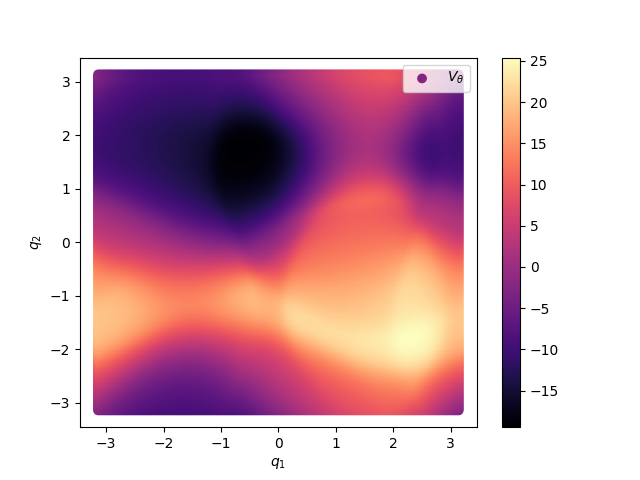}
         \caption{Control effort penalty $\alpha_{\text{eff}}=0.0001$.}
         \label{fig:reg=0.0001}
     \end{subfigure}
     \begin{subfigure}[b]{0.33\textwidth}
         \centering
         \includegraphics[width=1.0\textwidth]{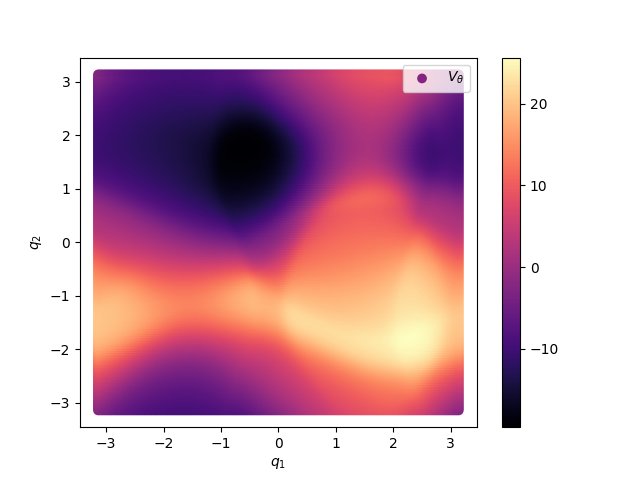}
         \caption{Control effort penalty $\alpha_{\text{eff}}=0.001$.}
         \label{fig:reg=0.0001}
     \end{subfigure}
     \begin{subfigure}[b]{0.33\textwidth}
         \centering
         \includegraphics[width=1.0\textwidth]{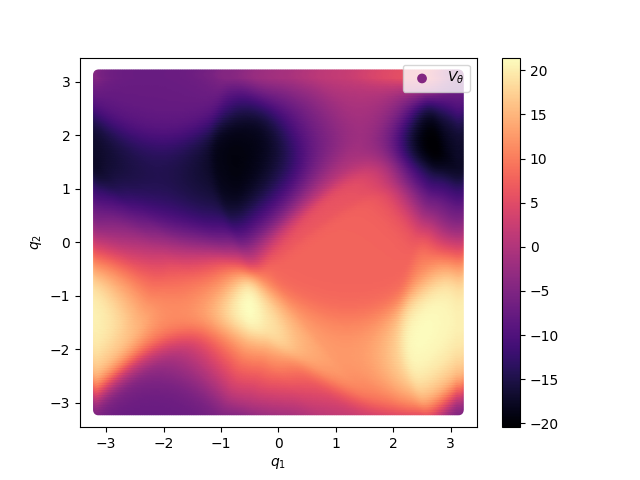}
         \caption{Control effort penalty $\alpha_{\text{eff}}=0.01$.}
         \label{fig:reg=0.0001}
     \end{subfigure}
        \caption{Learned potential for different effort penalty coefficients.}
        \label{fig:control_effort_penalizaiton_potential}
\end{figure}

Moreover, for a fixed period $T=\SI{1.5}{\second}$, it is possible to notice from Figure \ref{fig:control_effort_penalizaiton} and \ref{fig:control_effort_penalizaiton_potential} that the increase of the regularization penalty decreases the control effort (as expected) by improving the smoothness of the potential. 

\subsection{Learned Eigenmodes for Different Fixed Periods $T$}

\begin{figure}[h!]
     \centering
     \begin{subfigure}[b]{0.24\textwidth}
         \centering
         \includegraphics[width=1.0\textwidth]{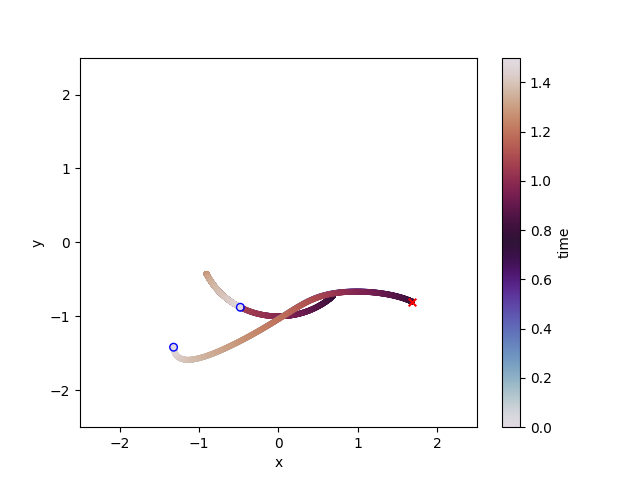}
         \caption{$T=\SI{1.5}{\second}$.}
         \label{}
     \end{subfigure}
     \hfill
     \begin{subfigure}[b]{0.24\textwidth}
         \centering
         \includegraphics[width=1.0\textwidth]{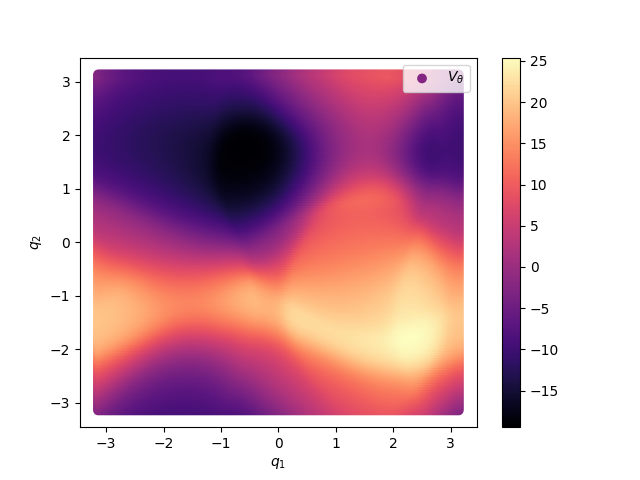}
         \caption{$T=\SI{1.5}{\second}$.}
         \label{}
     \end{subfigure}
     \hfill
     \begin{subfigure}[b]{0.24\textwidth}
         \centering
         \includegraphics[width=1.0\textwidth]{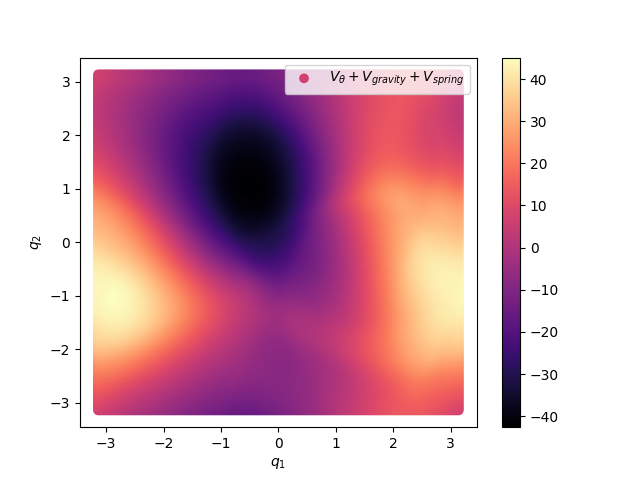}
         \caption{$T=\SI{1.5}{\second}$.}
         \label{}
     \end{subfigure}
          \hfill
     \begin{subfigure}[b]{0.24\textwidth}
         \centering
         \includegraphics[width=1.0\textwidth]{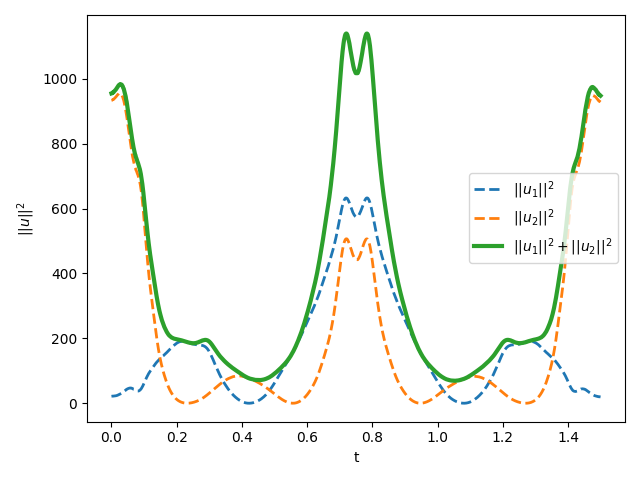}
         \caption{$T=\SI{1.55}{\second}$.}
         \label{}
     \end{subfigure}
     \hfill
     \begin{subfigure}[b]{0.24\textwidth}
         \centering
         \includegraphics[width=1.0\textwidth]{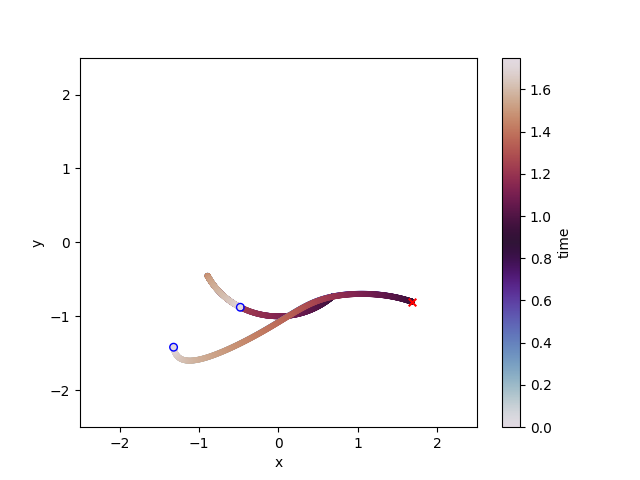}
         \caption{$T=\SI{1.75}{\second}$.}
         \label{}
     \end{subfigure}
     \hfill
     \begin{subfigure}[b]{0.24\textwidth}
         \centering
         \includegraphics[width=1.0\textwidth]{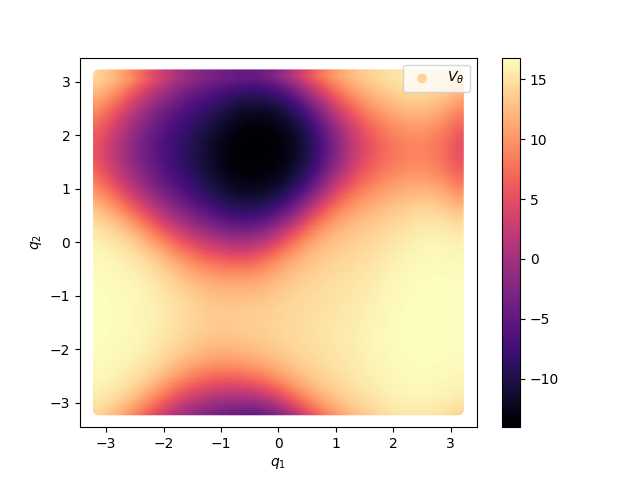}
         \caption{$T=\SI{1.75}{\second}$.}
         \label{}
     \end{subfigure}
     \hfill
     \begin{subfigure}[b]{0.24\textwidth}
         \centering
         \includegraphics[width=1.0\textwidth]{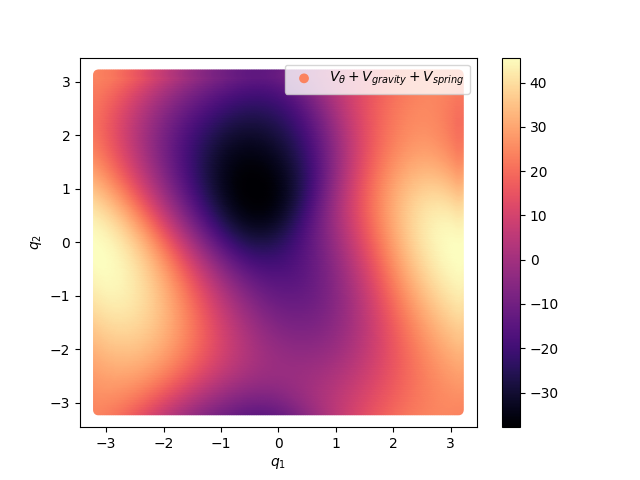}
         \caption{$T=\SI{1.75}{\second}$.}
         \label{}
     \end{subfigure}
          \hfill
     \begin{subfigure}[b]{0.24\textwidth}
         \centering
         \includegraphics[width=1.0\textwidth]{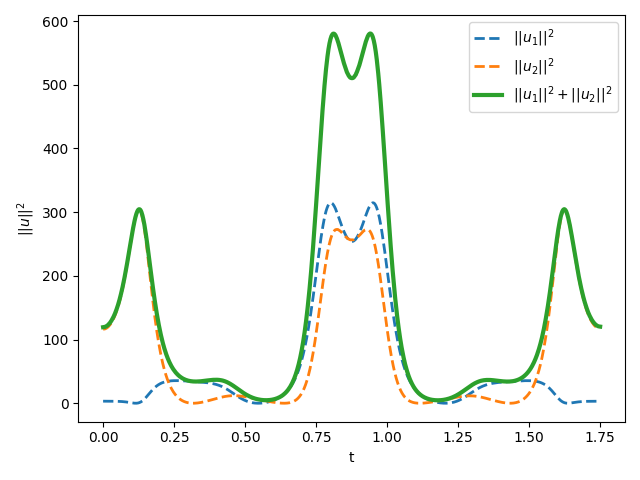}
         \caption{$T=\SI{1.75}{\second}$.}
         \label{}
     \end{subfigure}
     \hfill
          \begin{subfigure}[b]{0.24\textwidth}
         \centering
         \includegraphics[width=1.0\textwidth]{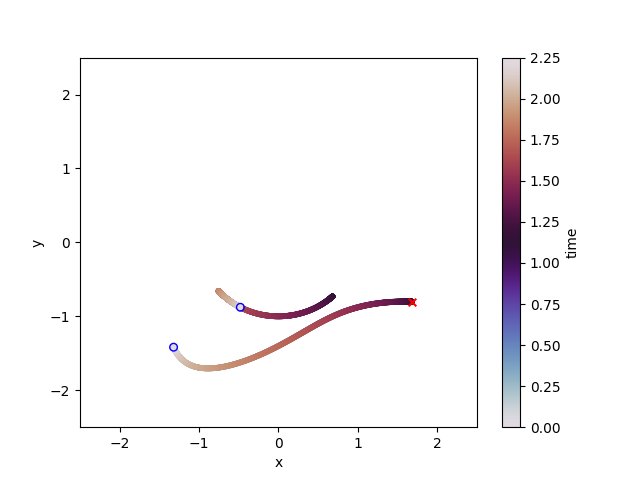}
         \caption{$T=\SI{2.25}{\second}$.}
         \label{}
     \end{subfigure}
     \hfill
     \begin{subfigure}[b]{0.24\textwidth}
         \centering
         \includegraphics[width=1.0\textwidth]{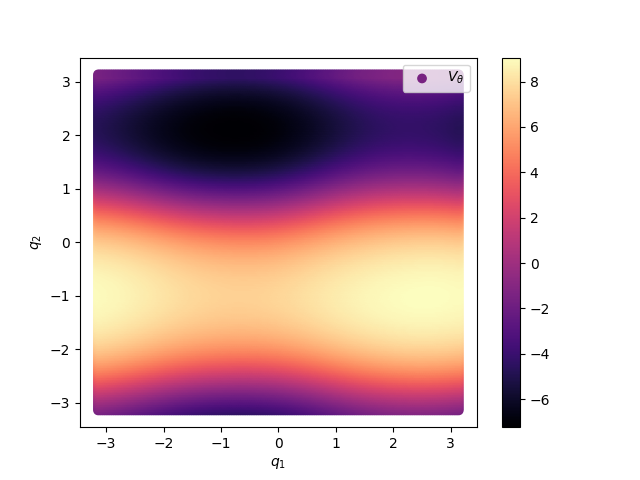}
         \caption{$T=\SI{2.25}{\second}$.}
         \label{}
     \end{subfigure}
     \hfill
     \begin{subfigure}[b]{0.24\textwidth}
         \centering
         \includegraphics[width=1.0\textwidth]{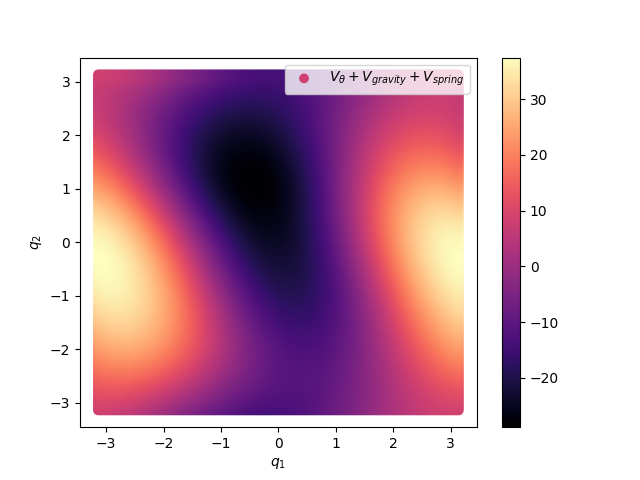}
         \caption{$T=\SI{2.25}{\second}$.}
         \label{}
     \end{subfigure}
          \hfill
     \begin{subfigure}[b]{0.24\textwidth}
         \centering
         \includegraphics[width=1.0\textwidth]{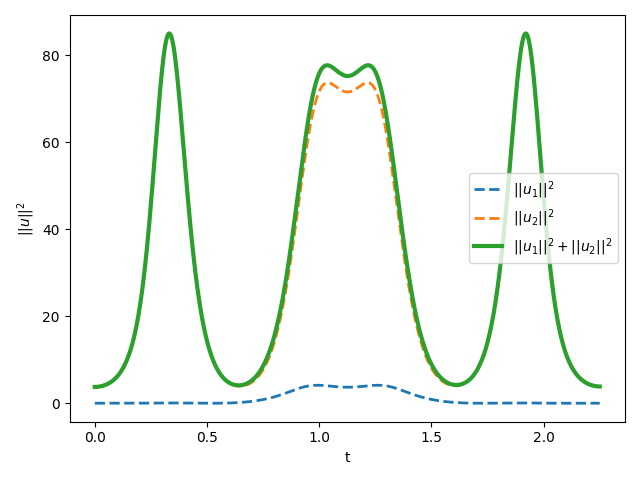}
         \caption{$T=\SI{2.25}{\second}$.}
         \label{}
     \end{subfigure}
     \hfill
          \begin{subfigure}[b]{0.24\textwidth}
         \centering
         \includegraphics[width=1.0\textwidth]{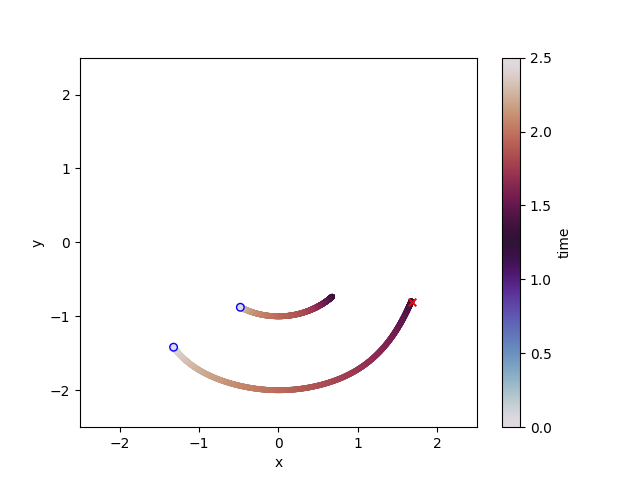}
         \caption{$T=\SI{2.5}{\second}$.}
         \label{}
     \end{subfigure}
     \hfill
     \begin{subfigure}[b]{0.24\textwidth}
         \centering
         \includegraphics[width=1.0\textwidth]{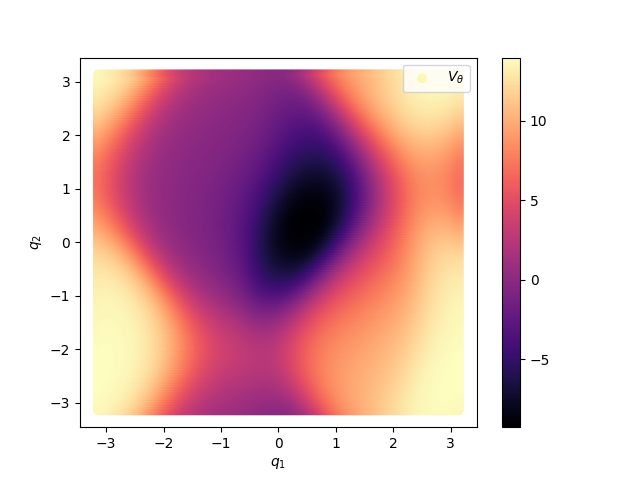}
         \caption{$T=\SI{2.5}{\second}$.}
         \label{}
     \end{subfigure}
     \hfill
     \begin{subfigure}[b]{0.24\textwidth}
         \centering
         \includegraphics[width=1.0\textwidth]{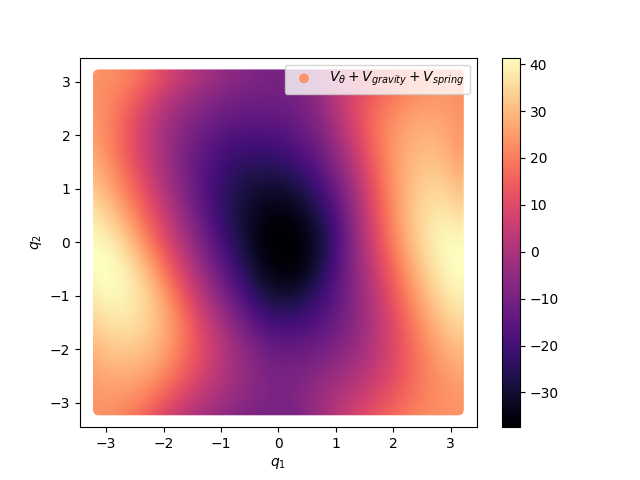}
         \caption{$T=\SI{2.5}{\second}$.}
         \label{}
     \end{subfigure}
          \hfill
     \begin{subfigure}[b]{0.24\textwidth}
         \centering
         \includegraphics[width=1.0\textwidth]{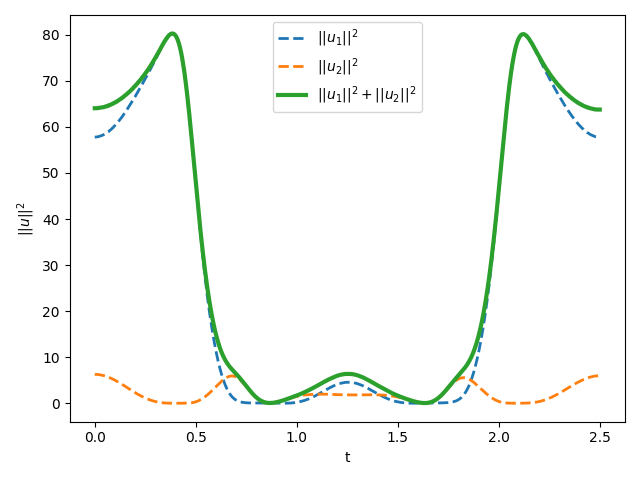}
         \caption{$T=\SI{2.5}{\second}$.}
         \label{}
     \end{subfigure}
          \hfill
          \begin{subfigure}[b]{0.24\textwidth}
         \centering
         \includegraphics[width=1.0\textwidth]{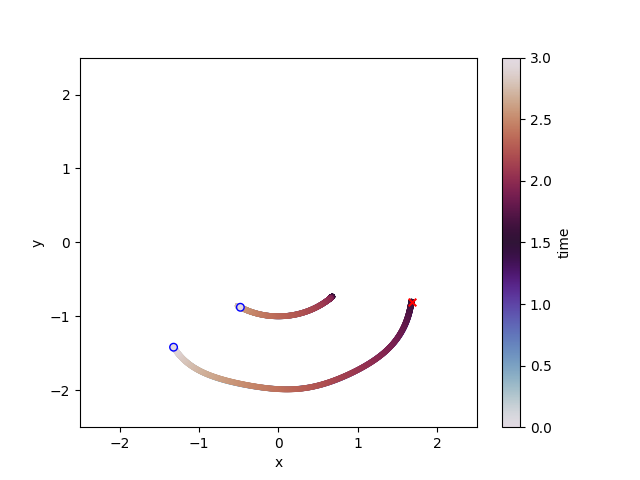}
         \caption{$T=\SI{3.0}{\second}$.}
         \label{fig:T=3_traj}
     \end{subfigure}
     \hfill
     \begin{subfigure}[b]{0.24\textwidth}
         \centering
         \includegraphics[width=1.0\textwidth]{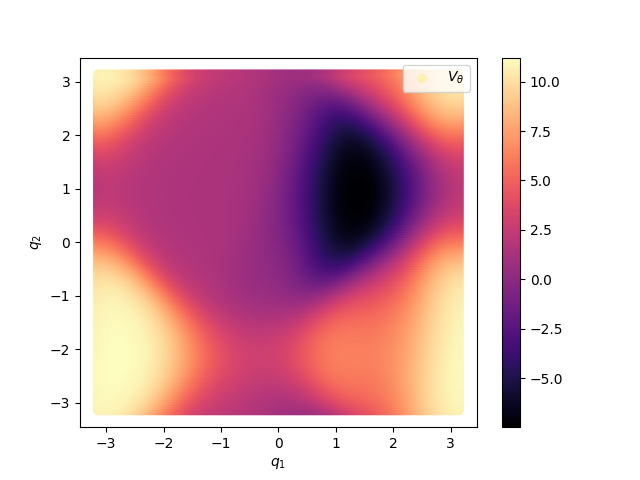}
         \caption{$T=\SI{3.0}{\second}$.}
         \label{}
     \end{subfigure}
          \hfill
     \begin{subfigure}[b]{0.24\textwidth}
         \centering
         \includegraphics[width=1.0\textwidth]{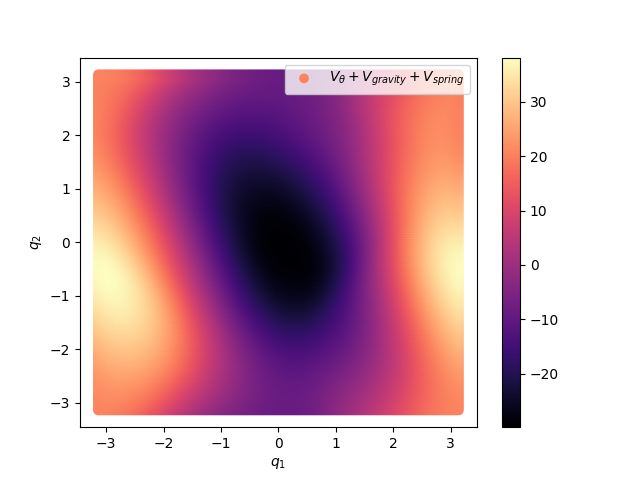}
         \caption{$T=\SI{3.0}{\second}$.}
         \label{fig:T=3_ceff}
     \end{subfigure}
     \hfill
     \begin{subfigure}[b]{0.24\textwidth}
         \centering
         \includegraphics[width=1.0\textwidth]{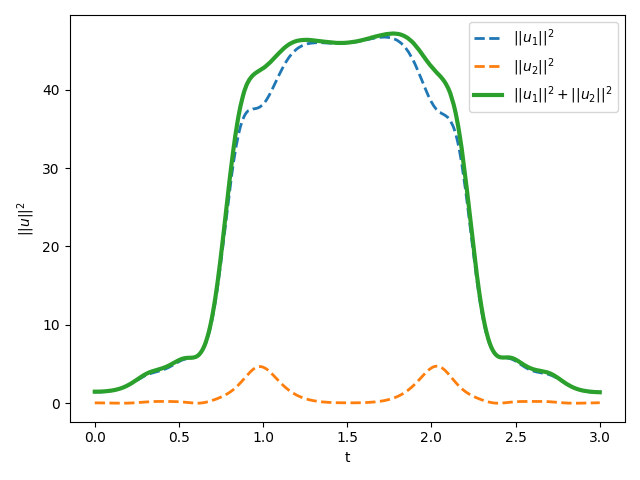}
         \caption{$T=\SI{3.0}{\second}$.}
         \label{fig:T=3_ceff}
     \end{subfigure}
        \caption{Trajectory, learned and total potential, and squared control effort for different period length $T \in \{\SI{1.5}{\second}, \SI{1.75}{\second}, \SI{2.25}{\second}, \SI{2.5}{\second}, \SI{3.0}{\second}\}$.}
        \label{fig:traj_for_varying_T}
\end{figure}

In Figure \ref{fig:traj_for_varying_T}, we show the resulting trajectories, learned potentials $V_{\bm{\theta}}$, and squared control effort $\bm{u}$ for different period length $T$. Our approach is capable of finding eigenmodes for different periods $T$. It is noticed that the learned potential combines with gravitational and elastic 
potentials in non trivial ways to steer the system on oscillatory modes with the desired period. In Figure \ref{fig:T=3_traj}-\ref{fig:T=3_ceff}, the period of oscillation is close to the natural evolution of the system, i.e. when only the gravitation potential is active and no learned potential is present, the learned potential is such that the resulting control effort is extremely small.

\clearpage
\section{Eigenmodes for a Different Target and Initial Position}\label{AppC:additional_results_config1}
In this appendix, we redo the eigenmode discovery experiment in Section \ref{sec:sim} with a different initial position and target. The numerical experiments are done for different values of the fixed period and different values of the control effort regularization  $\alpha_{\text{eff}}=0.0001$.

\begin{figure}[h!]
\centering
\begin{subfigure}{.33\textwidth}
  \centering
  \includegraphics[width=1.0\linewidth]{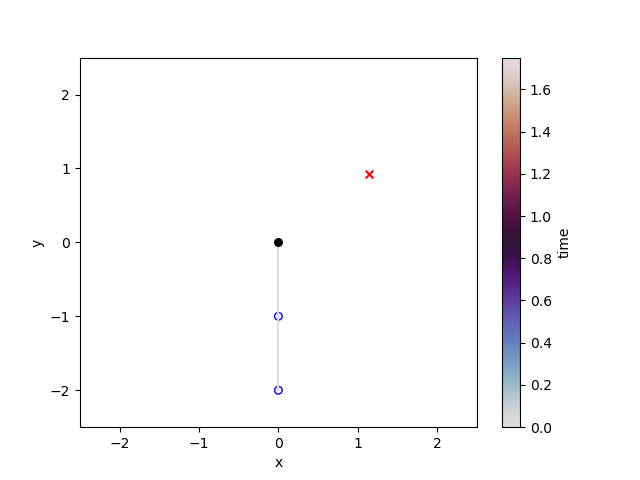}
  \caption{t=\SI{0.0}{\second}}
\end{subfigure}%
\begin{subfigure}{.33\textwidth}
  \centering
  \includegraphics[width=1.0\linewidth]{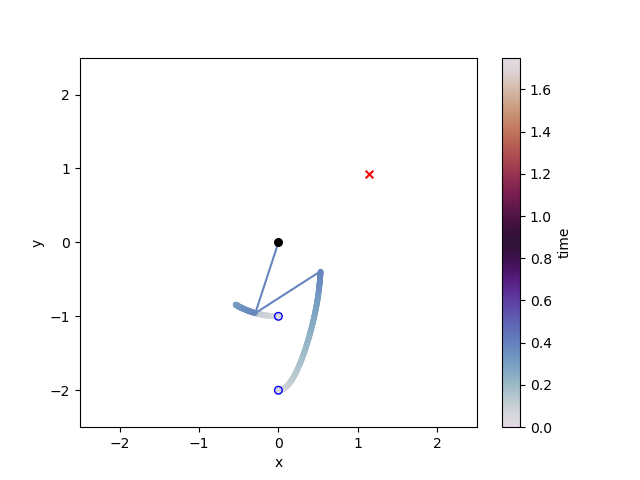}
  \caption{t=\SI{0.389}{\second}}
\end{subfigure}
 \begin{subfigure}{.33\textwidth}
   \centering
   \includegraphics[width=1.0\linewidth]{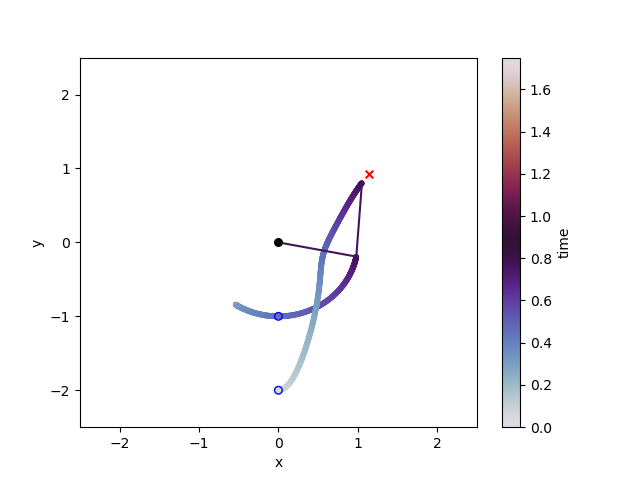}
   \caption{t=\SI{0.778}{\second}}
 \end{subfigure}
\begin{subfigure}{.33\textwidth}
  \centering
  \includegraphics[width=1.0\linewidth]{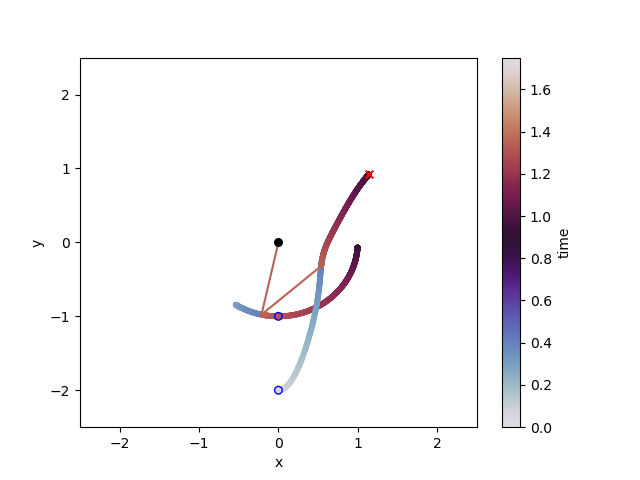}
  \caption{t=\SI{1.361}{\second}}
\end{subfigure}
\begin{subfigure}{.33\textwidth}
  \centering
  \includegraphics[width=1.0\linewidth]{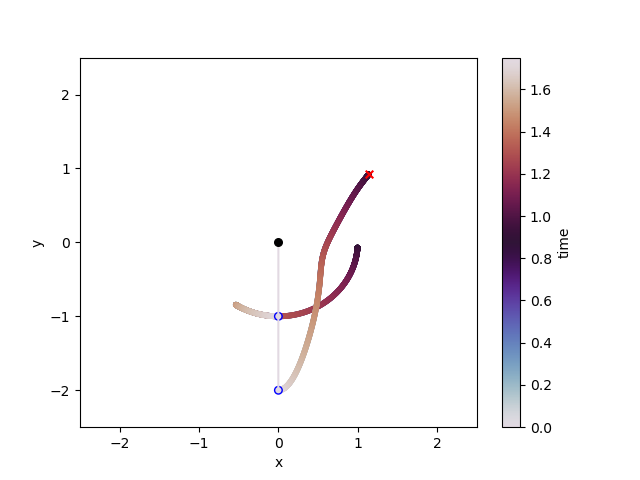}
  \caption{T=\SI{1.75}{\second}}
\end{subfigure}
\caption{Eigenmode at different time steps.}
\label{appC:T=1.75_traj}
\end{figure}

\begin{figure}[h!]
     \centering
     \begin{subfigure}[b]{0.33\textwidth}
         \centering
         \includegraphics[width=1.0\textwidth]{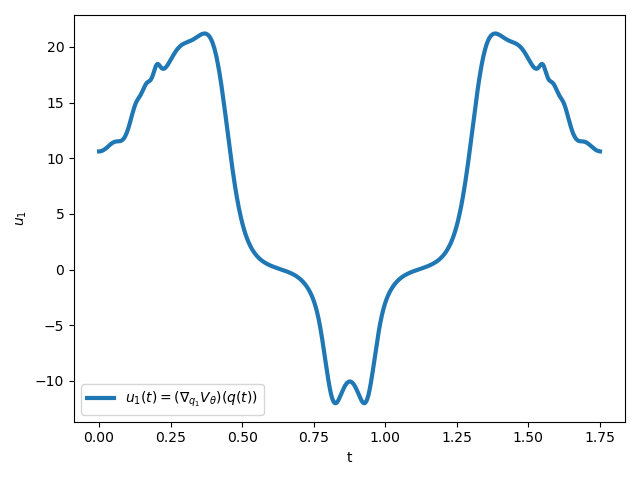}
         \caption{First component of $u(t)$.}
     \end{subfigure}
     \hfill
     \begin{subfigure}[b]{0.33\textwidth}
         \centering
         \includegraphics[width=1.0\textwidth]{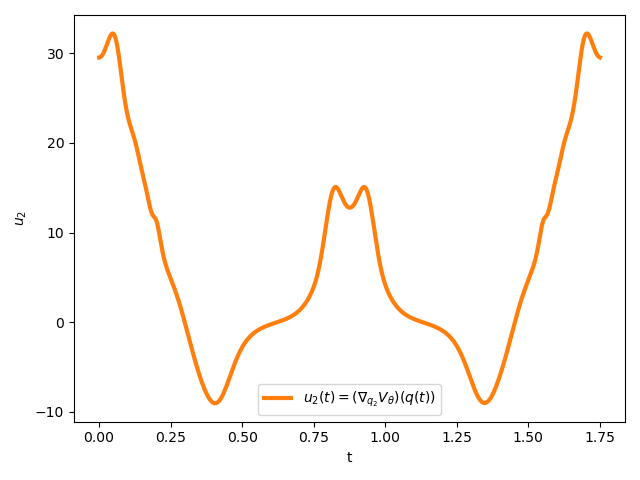}
         \caption{Second component of  $u(t)$.}
     \end{subfigure}
     \hfill
     \begin{subfigure}[b]{0.33\textwidth}
        \centering
        \includegraphics[width=1.0\textwidth]{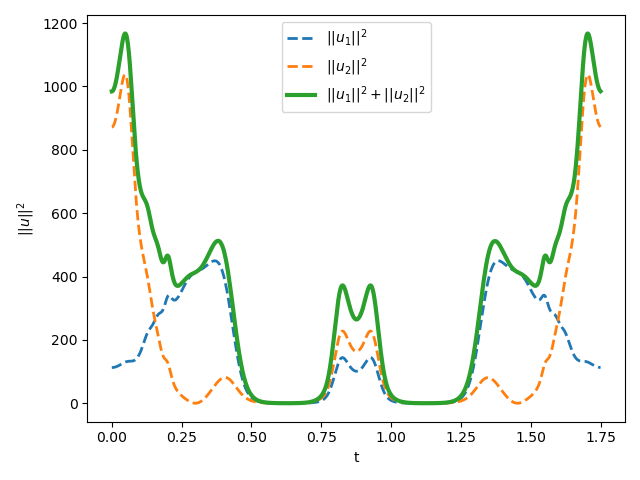}
         \caption{Squared control effort penalty.}
     \end{subfigure}
        \caption{Control inputs and control effort.}
        \label{appC:T=1.75_control_inp}
\end{figure}

\begin{figure}[h!]
     \centering
     \begin{subfigure}[b]{0.33\textwidth}
         \centering
         \includegraphics[width=1.0\textwidth]{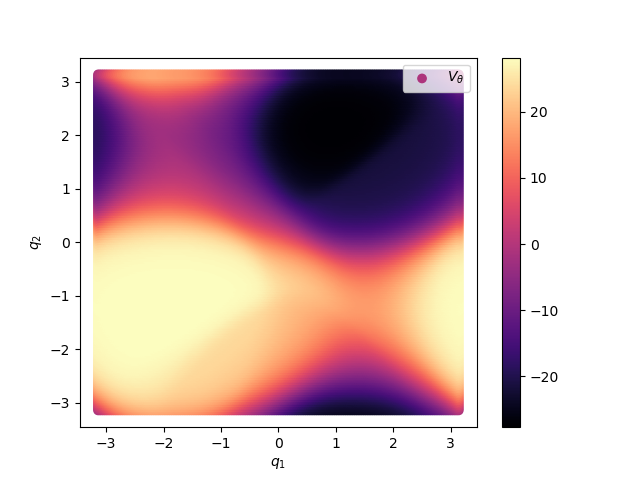}
         \caption{Learned potential $V_{\bm{\theta}}$.}
     \end{subfigure}
     \hfill
     \begin{subfigure}[b]{0.33\textwidth}
         \centering
          \includegraphics[width=1.0\textwidth]{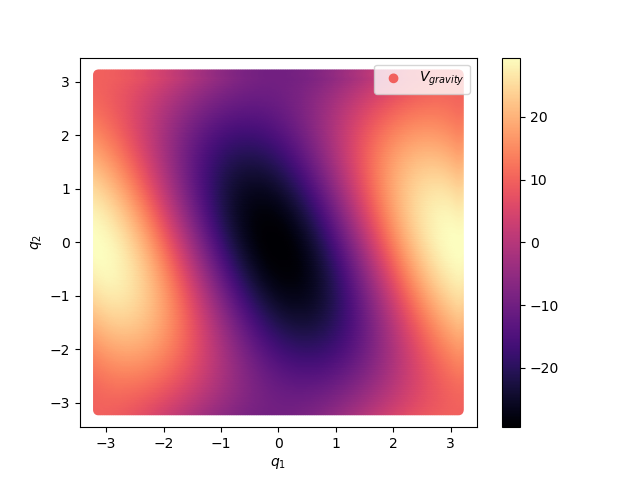}
         \caption{Gravitational potential $V_{\text{gravity}}$.}
     \end{subfigure}
     \hfill
     \begin{subfigure}[b]{0.33\textwidth}
         \centering
          \includegraphics[width=1.0\textwidth]{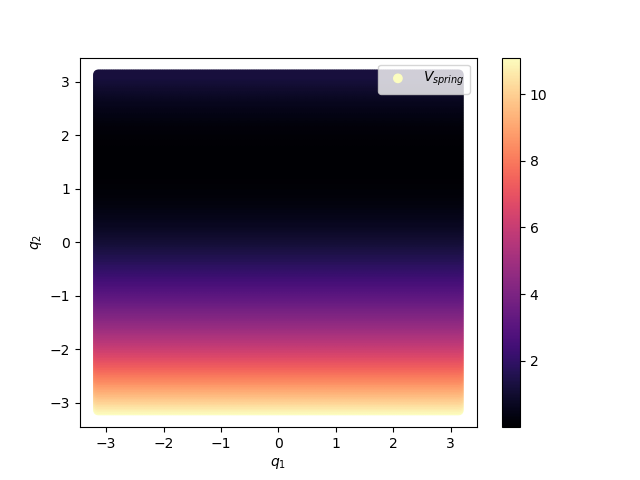}
         \caption{Spring potential $V_{spring}$}
     \end{subfigure}

     \begin{subfigure}[b]{0.33\textwidth}
         \centering
         \includegraphics[width=1.0\textwidth]{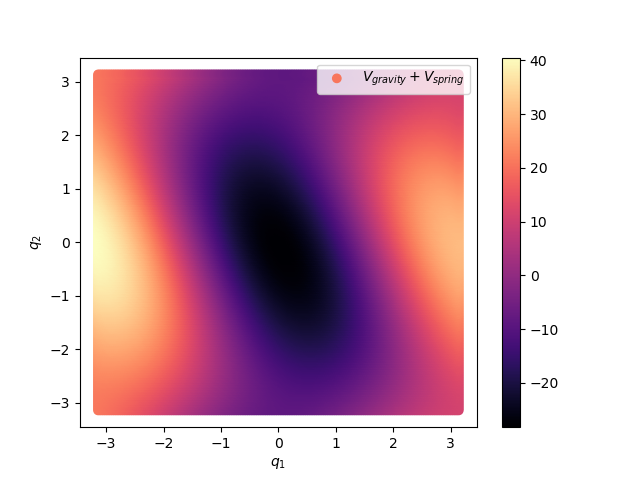}
         \caption{Open-loop potential $V_{\text{spring}}+V_{\text{gravity}}$.}
     \end{subfigure}
     \begin{subfigure}[b]{0.33\textwidth}
         \centering
          \includegraphics[width=1.0\textwidth]{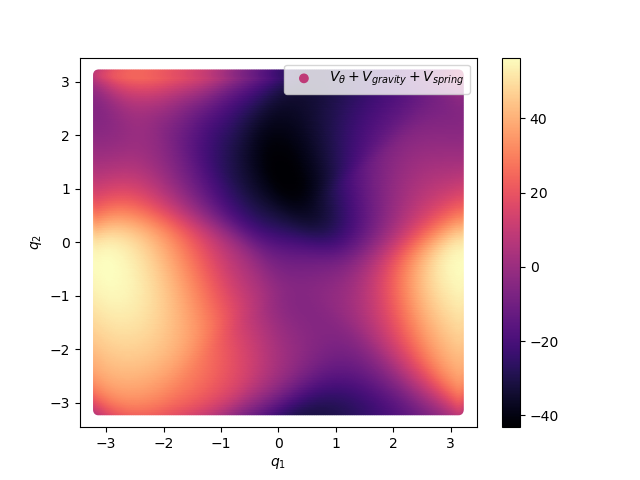}
         \caption{Overall potential $V_{\bm{\theta}}+V_{\text{spring}}+V_{\text{gravity}}$.}
     \end{subfigure}
    \caption{Potentials for T=\SI{1.75}{\second} over $\bm{q}\in[-\pi, \pi]$.}
    \label{appC:T=1.75_potentials}
\end{figure}

\begin{figure}[h!]
     \centering
     \begin{subfigure}[b]{0.33\textwidth}
         \centering
         \includegraphics[width=1.0\textwidth]{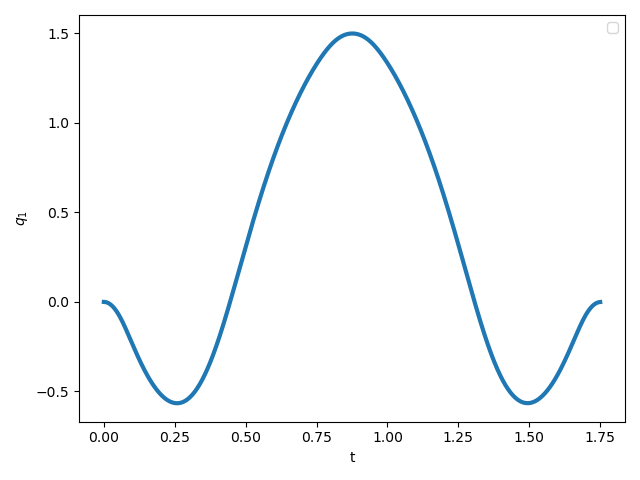}
         \caption{$q_1$ over the period.}
     \end{subfigure}
     \hfill
     \begin{subfigure}[b]{0.33\textwidth}
         \centering
         \includegraphics[width=1.0\textwidth]{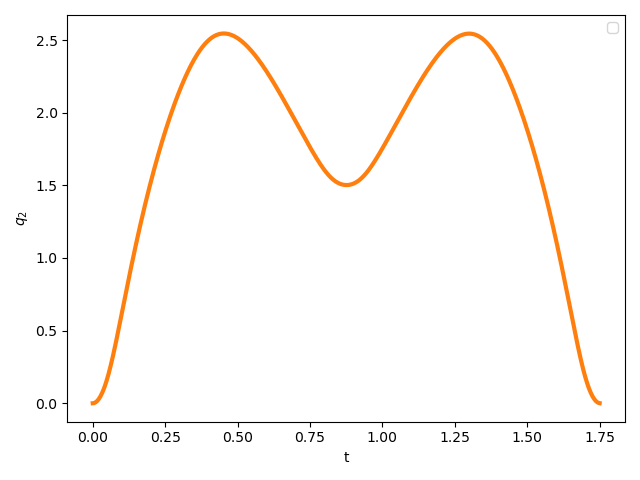}
         \caption{$q_2$ over the period.}
     \end{subfigure}
     \hfill
     \begin{subfigure}[b]{0.33\textwidth}
        \centering
        \includegraphics[width=1.0\textwidth]{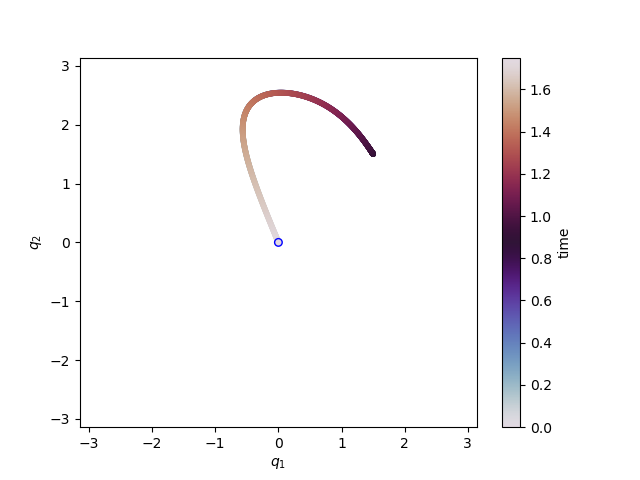}
         \caption{Trajectory in configuration space.}
     \end{subfigure}
        \caption{The time behavior of the angles $q_1$ and $q_2$ over one period.}
        \label{appC:T=1.75_q}
\end{figure}

\begin{figure}[h!]
\centering
\begin{subfigure}{.33\textwidth}
  \centering
  \includegraphics[width=1.0\linewidth]{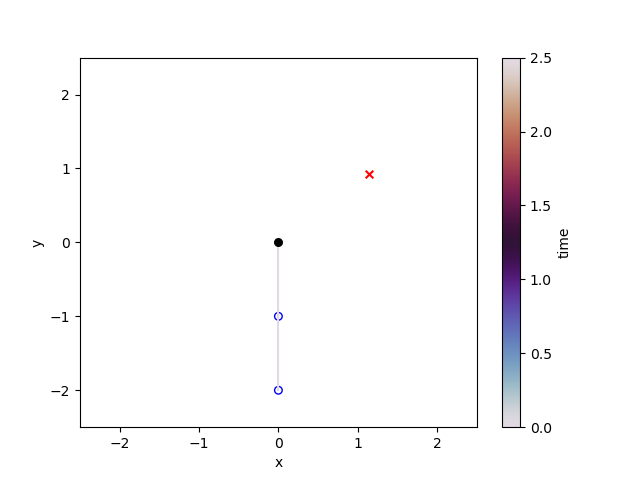}
  \caption{t=\SI{0.0}{\second}}
\end{subfigure}%
\begin{subfigure}{.33\textwidth}
  \centering
  \includegraphics[width=1.0\linewidth]{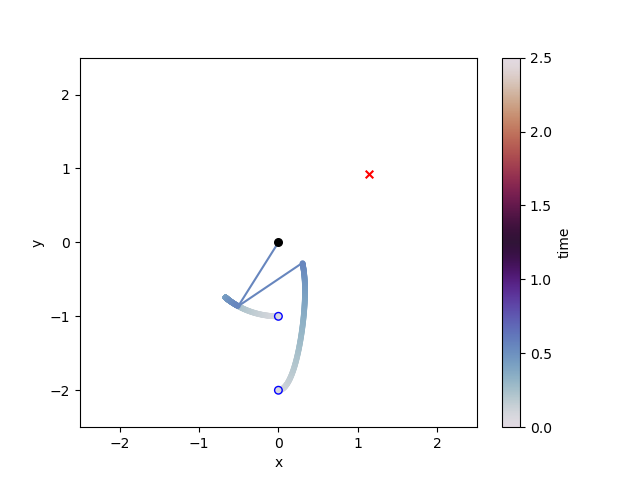}
  \caption{t=\SI{0.556}{\second}}
\end{subfigure}
\begin{subfigure}{.33\textwidth}
  \centering
  \includegraphics[width=1.0\linewidth]{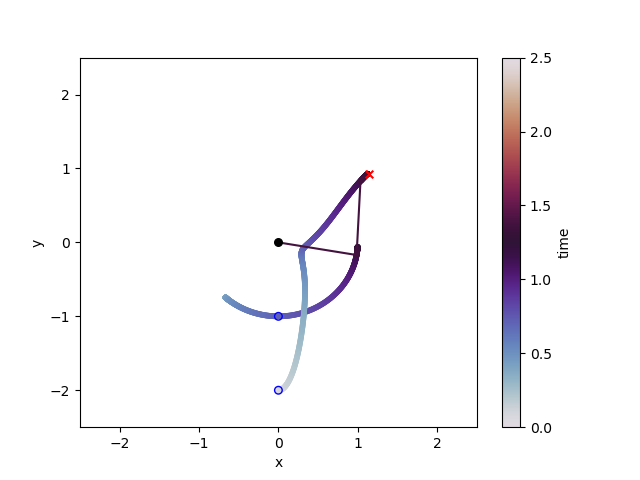}
  \caption{t=\SI{1.389}{\second}}
\end{subfigure}%
\hfill
\begin{subfigure}{.33\textwidth}
  \centering
  \includegraphics[width=1.0\linewidth]{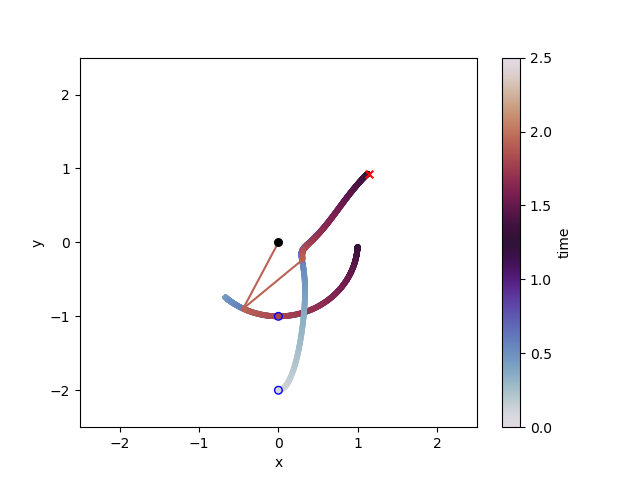}
  \caption{t=\SI{1.944}{\second}}
\end{subfigure}
\begin{subfigure}{.33\textwidth}
  \centering
  \includegraphics[width=1.0\linewidth]{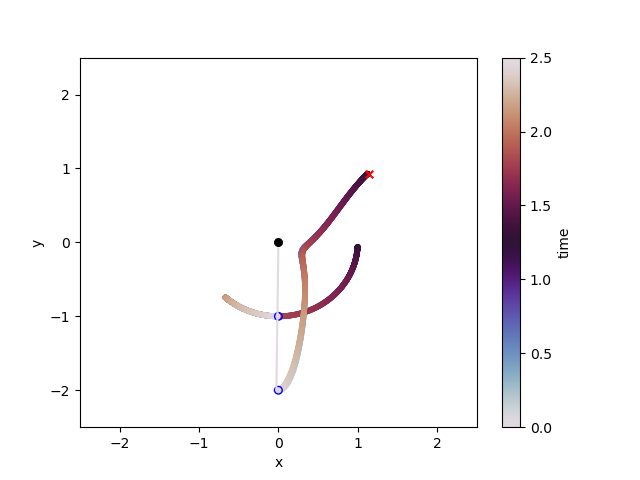}
  \caption{T=\SI{2.50}{\second}}
\end{subfigure}
\caption{Eigenmode at different time steps.}
\label{appC:T=2.5_traj}
\end{figure}

\begin{figure}[h!]
     \centering
     \begin{subfigure}[b]{0.33\textwidth}
         \centering
         \includegraphics[width=1.0\textwidth]{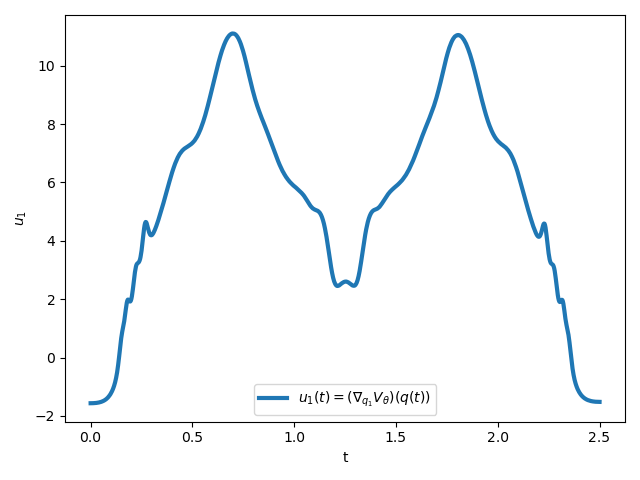}
         \caption{First component of $u(t)$.}
     \end{subfigure}
     \hfill
     \begin{subfigure}[b]{0.33\textwidth}
         \centering
         \includegraphics[width=1.0\textwidth]{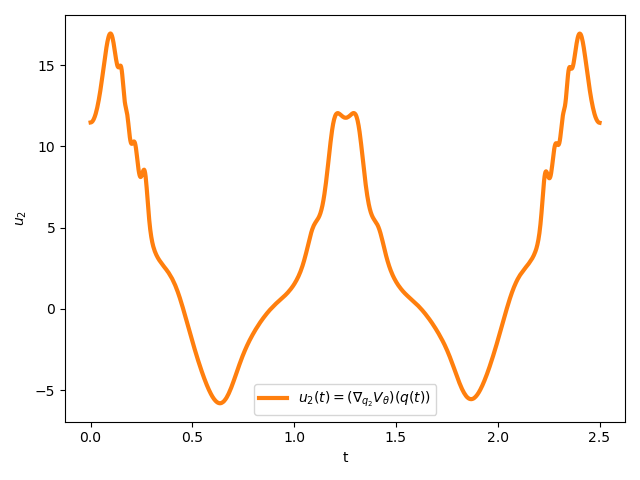}
         \caption{Second component of $u(t)$.}
     \end{subfigure}
     \hfill
     \begin{subfigure}[b]{0.33\textwidth}
        \centering
        \includegraphics[width=1.0\textwidth]{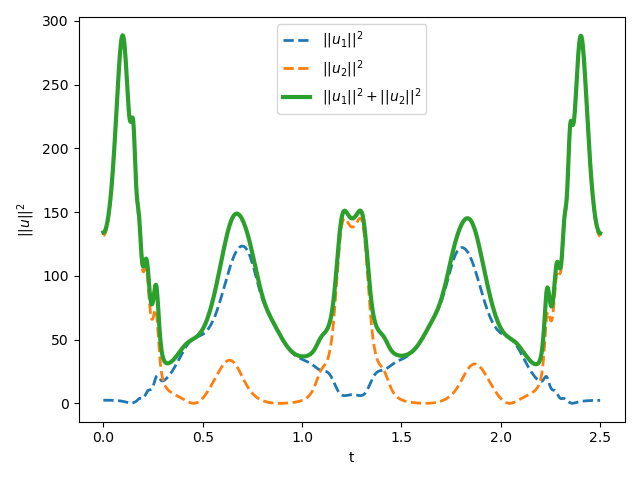}
         \caption{Squared control effort penalty.}
     \end{subfigure}
        \caption{Control inputs and control effort.}
        \label{appC:T=2.5_control_inp}
\end{figure}

\begin{figure}[h!]
     \centering
     \begin{subfigure}[b]{0.33\textwidth}
         \centering
         \includegraphics[width=1.0\textwidth]{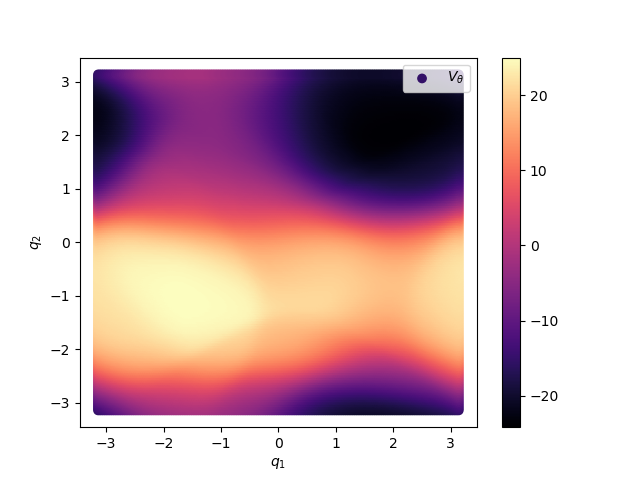}
         \caption{Learned potential $V_{\bm{\theta}}$.}
     \end{subfigure}
     \hfill
     \begin{subfigure}[b]{0.33\textwidth}
         \centering
          \includegraphics[width=1.0\textwidth]{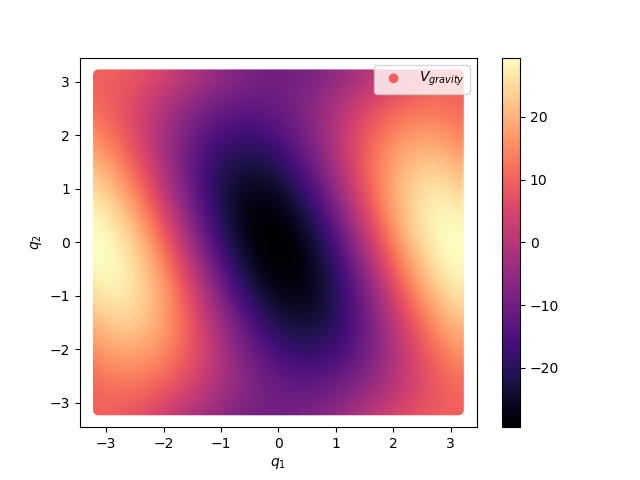}
         \caption{Gravitational potential $V_{\text{gravity}}$.}
     \end{subfigure}
     \hfill
     \begin{subfigure}[b]{0.33\textwidth}
         \centering
          \includegraphics[width=1.0\textwidth]{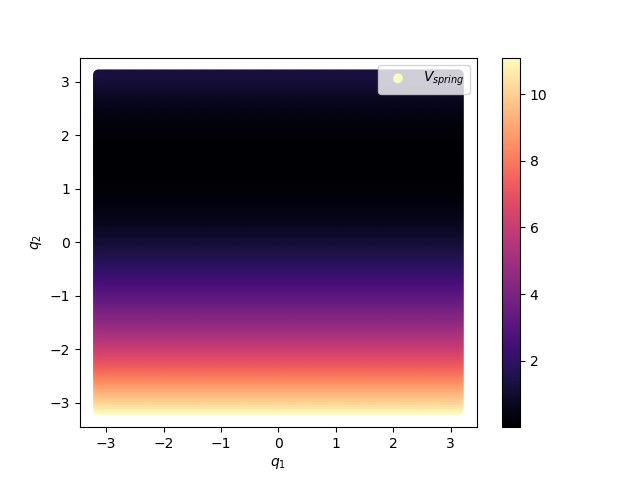}
         \caption{Spring potential $V_{spring}$}
     \end{subfigure}

     \begin{subfigure}[b]{0.33\textwidth}
         \centering
         \includegraphics[width=1.0\textwidth]{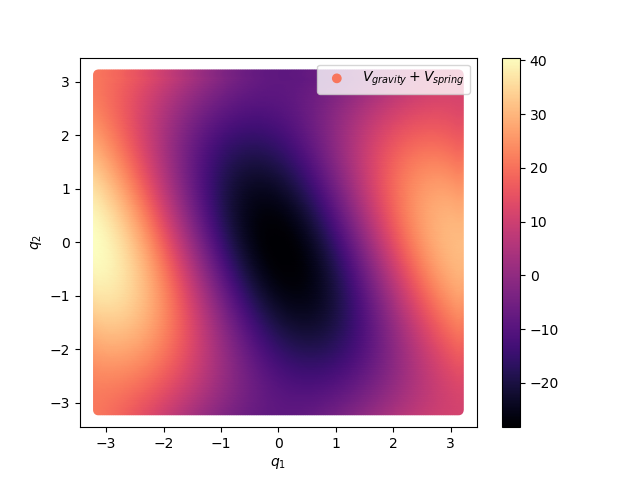}
         \caption{Open-loop potential $V_{\text{spring}}+V_{\text{gravity}}$.}
     \end{subfigure}
     \begin{subfigure}[b]{0.33\textwidth}
         \centering
          \includegraphics[width=1.0\textwidth]{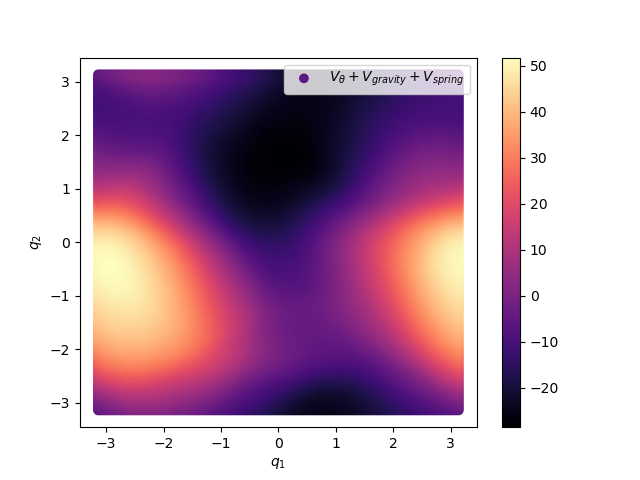}
         \caption{Overall potential $V_{\bm{\theta}}+V_{\text{spring}}+V_{\text{gravity}}$.}
     \end{subfigure}
        \caption{Potentials for T=\SI{2.50}{\second} over $\bm{q}\in[-\pi, \pi]$.}
        \label{appC:T=2.5_potentials}
\end{figure}

\begin{figure}[h!]
     \centering
     \begin{subfigure}[b]{0.33\textwidth}
         \centering
         \includegraphics[width=1.0\textwidth]{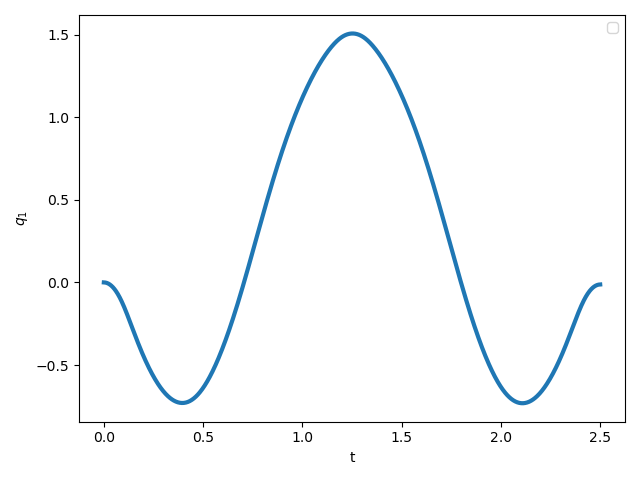}
         \caption{$q_1$ over the period.}
     \end{subfigure}
     \hfill
     \begin{subfigure}[b]{0.33\textwidth}
         \centering
         \includegraphics[width=1.0\textwidth]{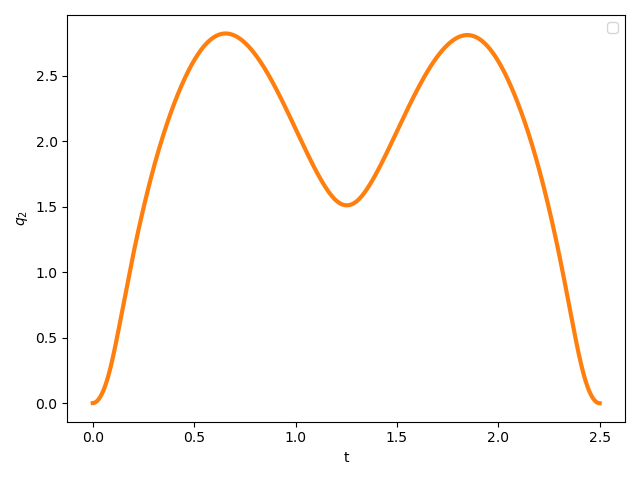}
         \caption{$q_2$ over the period.}
     \end{subfigure}
     \hfill
     \begin{subfigure}[b]{0.33\textwidth}
        \centering
        \includegraphics[width=1.0\textwidth]{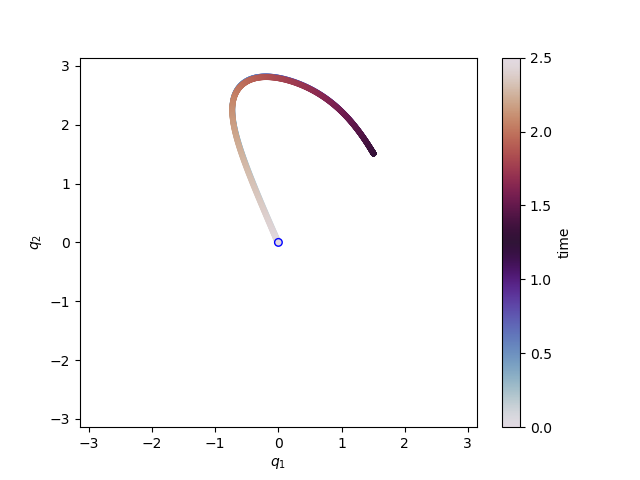}
         \caption{Trajectory in configuration space.}
     \end{subfigure}
        \caption{The time behavior of the angles $q_1$ and $q_2$ over one period.}
        \label{appC:T=2.5_q}
\end{figure}

\begin{figure}[h!]
\centering
\begin{subfigure}{.33\textwidth}
  \centering
  \includegraphics[width=1.0\linewidth]{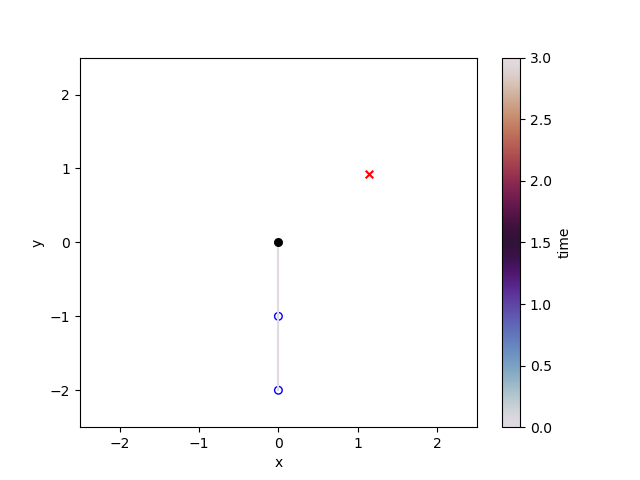}
  \caption{t=\SI{0.0}{\second}}
\end{subfigure}%
\begin{subfigure}{.33\textwidth}
  \centering
  \includegraphics[width=1.0\linewidth]{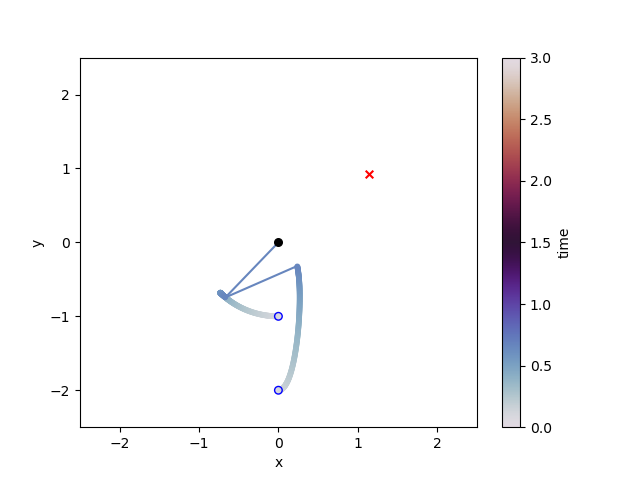}
  \caption{t=\SI{0.667}{\second}}
\end{subfigure}
\begin{subfigure}{.33\textwidth}
  \centering
  \includegraphics[width=1.0\linewidth]{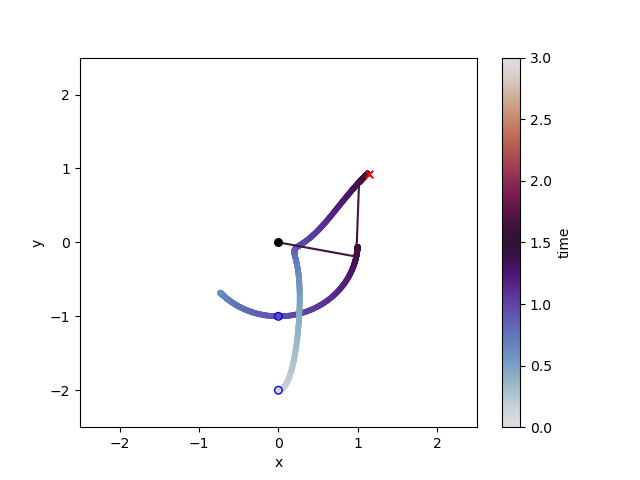}
  \caption{t=\SI{1.667}{\second}}
\end{subfigure}%
\hfill
\begin{subfigure}{.33\textwidth}
  \centering
  \includegraphics[width=1.0\linewidth]{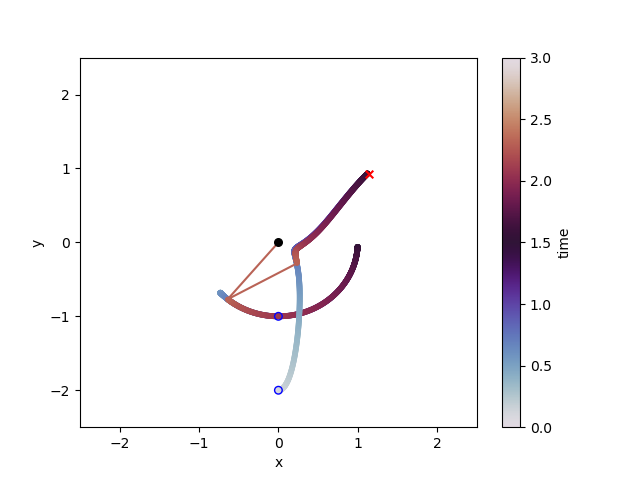}
  \caption{t=\SI{2.333}{\second}}
\end{subfigure}
\begin{subfigure}{.33\textwidth}
  \centering
  \includegraphics[width=1.0\linewidth]{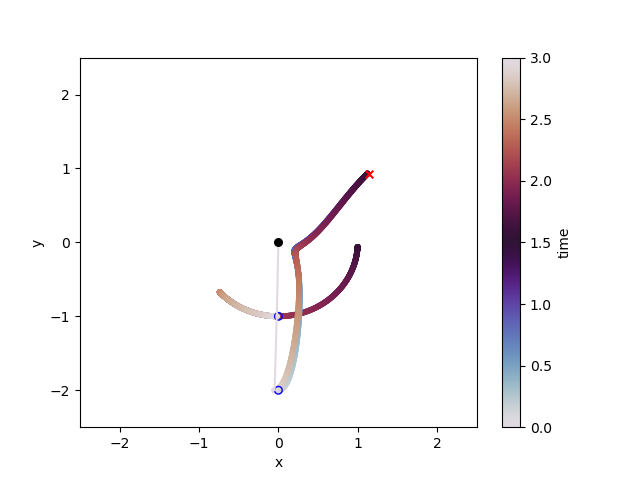}
  \caption{T=\SI{3.0}{\second}}
\end{subfigure}
\caption{Eigenmode at different time steps.}
\label{appC:T=3.0_traj}
\end{figure}

\begin{figure}[h!]
     \centering
     \begin{subfigure}[b]{0.33\textwidth}
         \centering
         \includegraphics[width=1.0\textwidth]{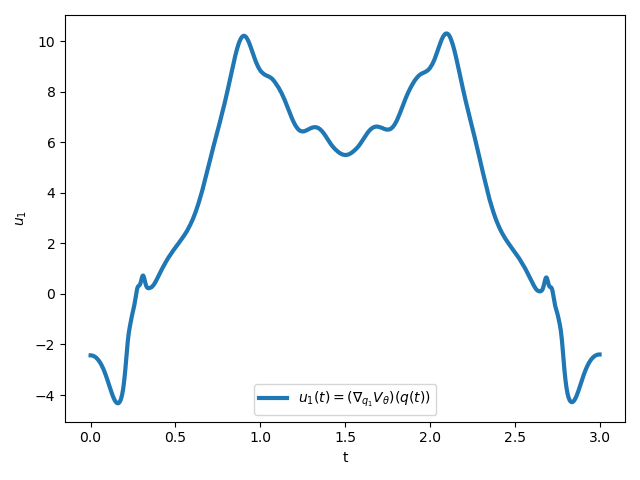}
         \caption{First component of $u(t)$.}
     \end{subfigure}
     \hfill
     \begin{subfigure}[b]{0.33\textwidth}
         \centering
         \includegraphics[width=1.0\textwidth]{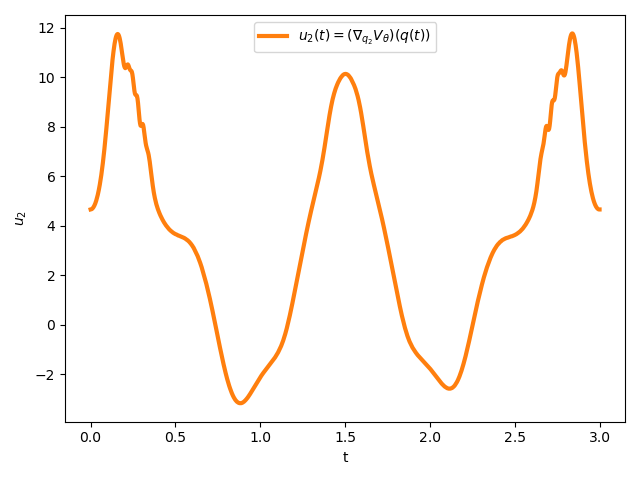}
         \caption{Second component of $u(t)$.}
     \end{subfigure}
     \hfill
     \begin{subfigure}[b]{0.33\textwidth}
        \centering
        \includegraphics[width=1.0\textwidth]{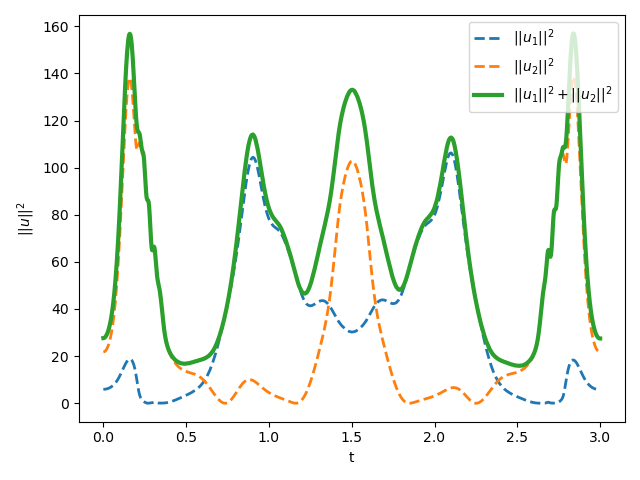}
         \caption{Squared control effort penalty.}
     \end{subfigure}
        \caption{Control inputs and control effort.}
        \label{appC:T=3.0_control_inp}
\end{figure}

\begin{figure}[h!]
     \centering
     \begin{subfigure}[b]{0.33\textwidth}
         \centering
         \includegraphics[width=1.0\textwidth]{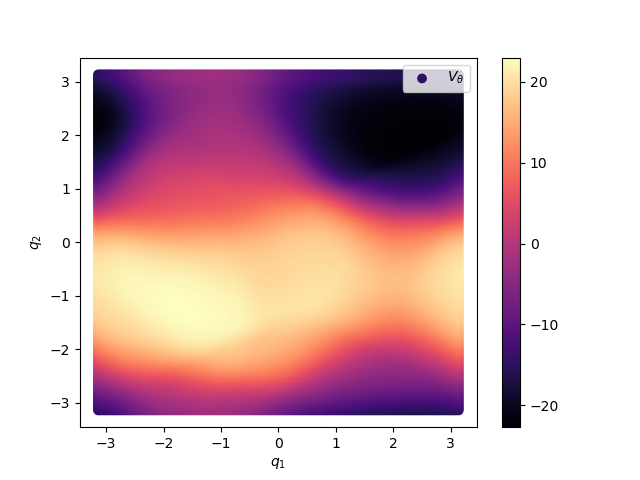}
         \caption{Learned potential $V_{\bm{\theta}}$.}
     \end{subfigure}
     \hfill
     \begin{subfigure}[b]{0.33\textwidth}
         \centering
          \includegraphics[width=1.0\textwidth]{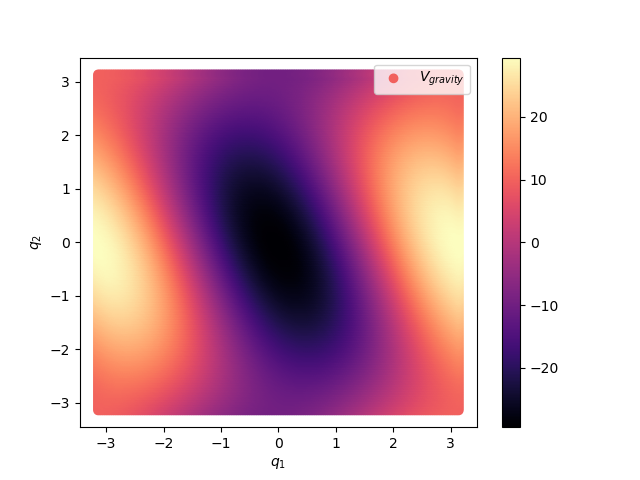}
         \caption{Gravitational potential $V_{\text{gravity}}$.}
     \end{subfigure}
     \hfill
     \begin{subfigure}[b]{0.33\textwidth}
         \centering
          \includegraphics[width=1.0\textwidth]{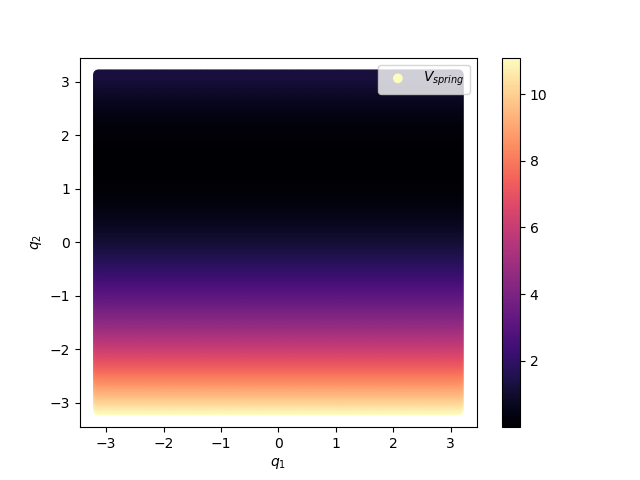}
         \caption{Spring potential $V_{\text{spring}}$}
     \end{subfigure}

     \begin{subfigure}[b]{0.33\textwidth}
         \centering
         \includegraphics[width=1.0\textwidth]{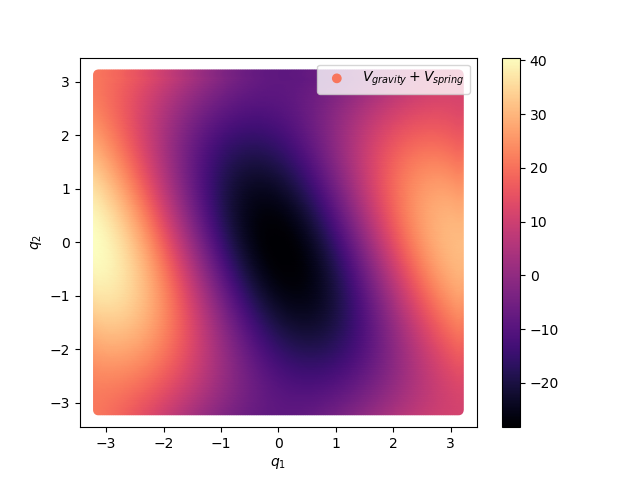}
         \caption{Spring and gravitational potential $V_{\text{spring}}+V_{\text{gravity}}$.}
     \end{subfigure}
     \begin{subfigure}[b]{0.33\textwidth}
         \centering
          \includegraphics[width=1.0\textwidth]{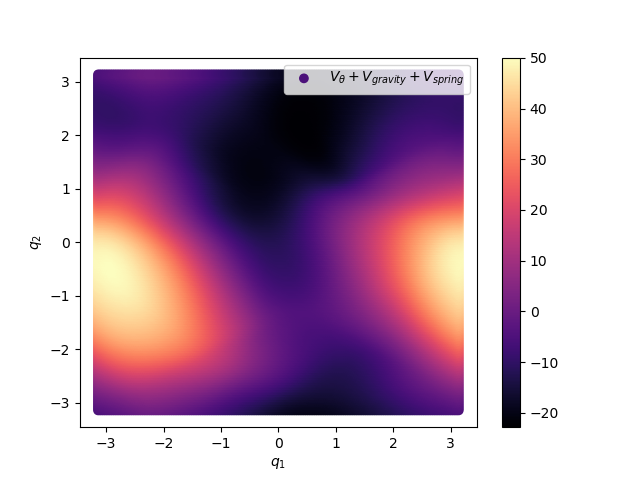}
         \caption{Overall potential $V_{\bm{\theta}}+V_{\text{spring}}+V_{\text{gravity}}$.}
     \end{subfigure}
        \caption{Potentials for T=\SI{3.0}{\second} over $\bm{q}\in[-\pi, \pi]$.}
        \label{appC:T=3.0_potential}
\end{figure}

\begin{figure}[h!]
     \centering
     \begin{subfigure}[b]{0.33\textwidth}
         \centering
         \includegraphics[width=1.0\textwidth]{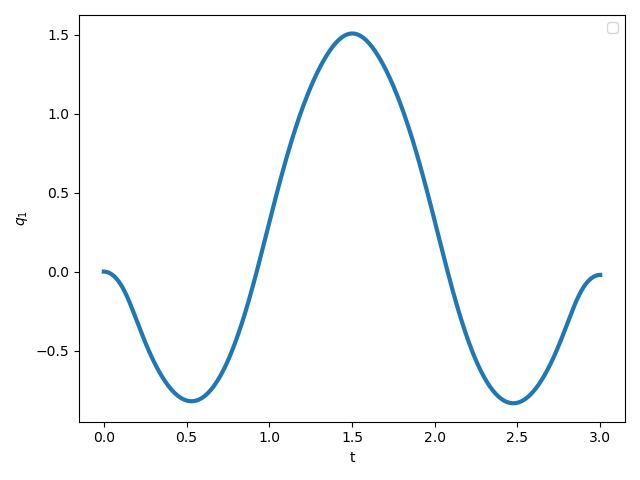}
         \caption{$q_1$ over one period.}
     \end{subfigure}
     \hfill
     \begin{subfigure}[b]{0.33\textwidth}
         \centering
         \includegraphics[width=1.0\textwidth]{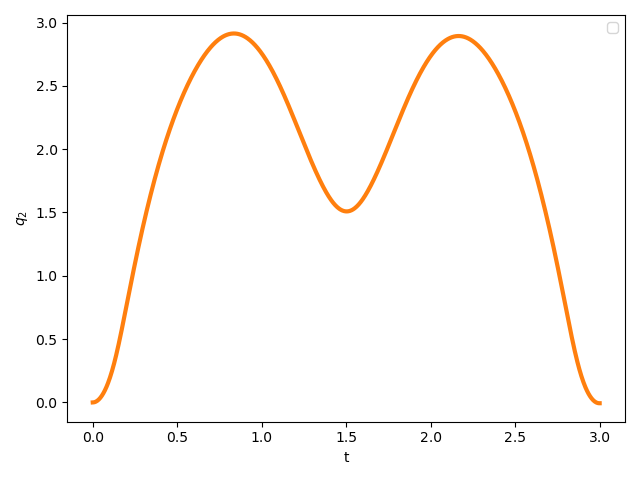}
         \caption{$q_2$ over one period.}
     \end{subfigure}
     \hfill
     \begin{subfigure}[b]{0.33\textwidth}
        \centering
        \includegraphics[width=1.0\textwidth]{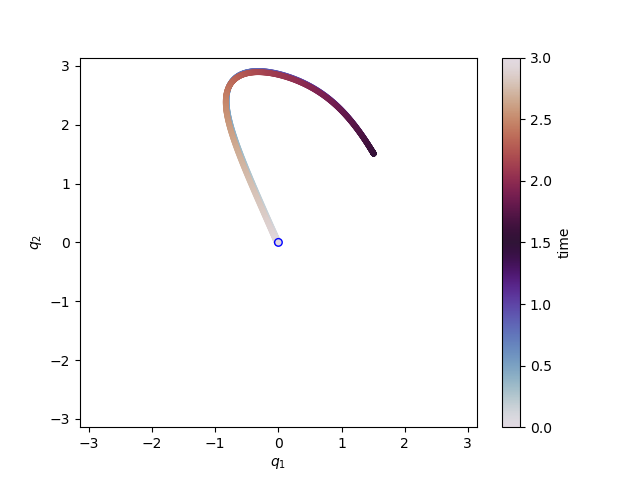}
         \caption{Trajectory in configuration space.}
     \end{subfigure}
        \caption{The time behavior of the angles $q_1$ and $q_2$ over one period.}
        \label{appC:T=3.0_q}
\end{figure}

For each value of the period $T \in[1.75, 2.5, 3.0]\si{\second}$, we show the trajectories of the double pendulum in Figure \ref{appC:T=1.75_traj}, \ref{appC:T=2.5_traj} and \ref{appC:T=3.0_traj}, the control inputs in Figure \ref{appC:T=1.75_control_inp}, \ref{appC:T=2.5_control_inp} and \ref{appC:T=3.0_control_inp}, the potentials in Figure \ref{appC:T=1.75_potentials}, \ref{appC:T=2.5_potentials} and \ref{appC:T=3.0_potential},  and the state variable over time in Figure \ref{appC:T=1.75_q}, \ref{appC:T=2.5_q} and \ref{appC:T=3.0_q}, respectively.


\clearpage




 \clearpage


\clearpage

\section{Learning the Period $T$ jointly with the Potential}\label{AppD:learnableT}
In this appendix, we briefly elaborate on a small extension of the optimization problem in \eqref{eq:optimization_problem_eigmode_formal_fixed_q0_and_fixed_T}. In the optimization problem in \eqref{eq:optimization_problem_eigmode_formal_fixed_q0_and_fixed_T}, the desired period $T$ of the eigenmode is fixed a priori\footnote{A requirement on the period length is often needed in a pick-and-place task in the context of an automatic machine where timing is crucial.}. For the experiments in Section \ref{sec:sim}, this suffices. However, in some applications it might be necessary to learn a suitable period $T$ of the eigenmode since it might be unknown. 
As a consequence, we extend the model in  \eqref{eq:optimization_problem_eigmode_formal_fixed_q0_and_fixed_T} to allow for an optimizable period $T$:
\begin{equation}
\begin{aligned}
\min_{\bm{\theta}, T} \quad & L_{\text{task}}(\bm{x}) + \beta L_{\textrm{eigen}}(\bm{x}) \\
\textrm{s.t.} \quad &
        \frac{\mathrm{d}}{\mathrm{d}t} \begin{bmatrix}
           \bm{q}(t) \\
           \bm{p}(t)
         \end{bmatrix}  = \begin{bmatrix}
            0 & \bm{I} \\ -\bm{I} & 0
         \end{bmatrix} \nabla (H+V_{\bm{\theta}})(\bm{p}, \bm{q}) \quad
        , \begin{bmatrix}
           \bm{p}(0) \\
           \bm{q}(0)
         \end{bmatrix} & = \begin{bmatrix}
           \bm{0} \\
           \bm{q}_0
         \end{bmatrix}
\end{aligned}
\label{eq:optimization_problem_eigmode_to_solve_learnableT}
\end{equation}




In the remainder of this appendix, we show numerical experiments similar to the experiments in Section \ref{sec:sim} but with a learnable period $T$. We present two such numerical experiments, each with a different pair of initial configuration $\bm{q}_0$ and target position $h^*$.

\subsection{Results}
\begin{figure}[h!]
\centering
\begin{subfigure}{.33\textwidth}
  \centering
  \includegraphics[width=1.0\linewidth]{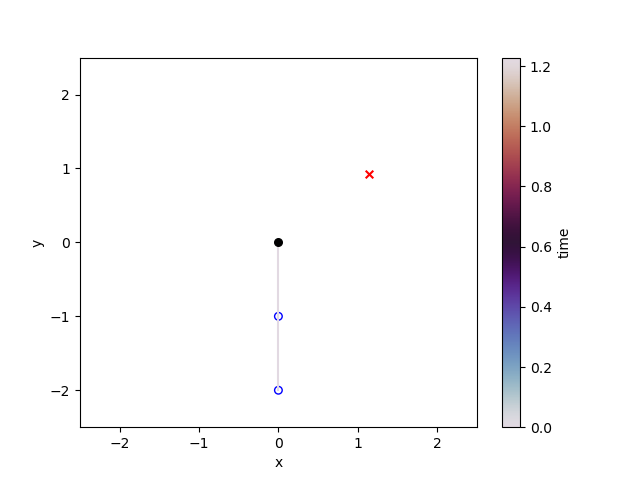}
  \caption{t=\SI{0.0}{\second}}
\end{subfigure}%
\begin{subfigure}{.33\textwidth}
  \centering
  \includegraphics[width=1.0\linewidth]{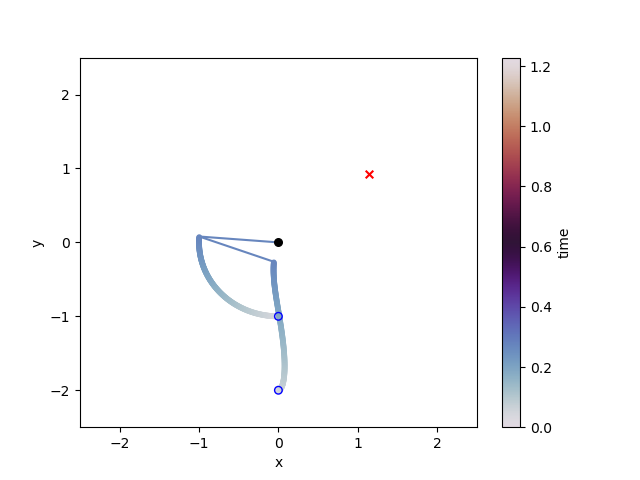}
  \caption{t=\SI{0.273}{\second}}
\end{subfigure}
\begin{subfigure}{.33\textwidth}
  \centering
  \includegraphics[width=1.0\linewidth]{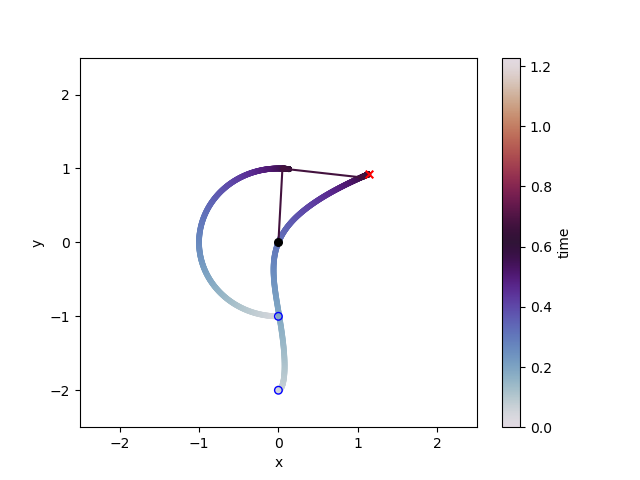}
  \caption{t=\SI{0.682}{\second}}
\end{subfigure}%
\hfill
\begin{subfigure}{.33\textwidth}
  \centering
  \includegraphics[width=1.0\linewidth]{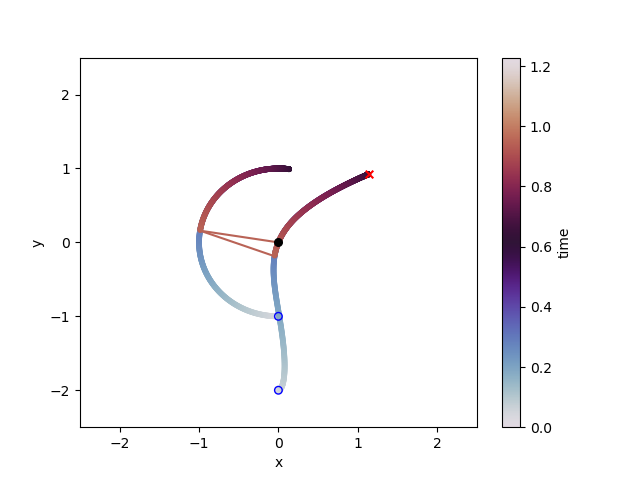}
  \caption{t=\SI{0.955}{\second}}
\end{subfigure}
\begin{subfigure}{.33\textwidth}
  \centering
  \includegraphics[width=1.0\linewidth]{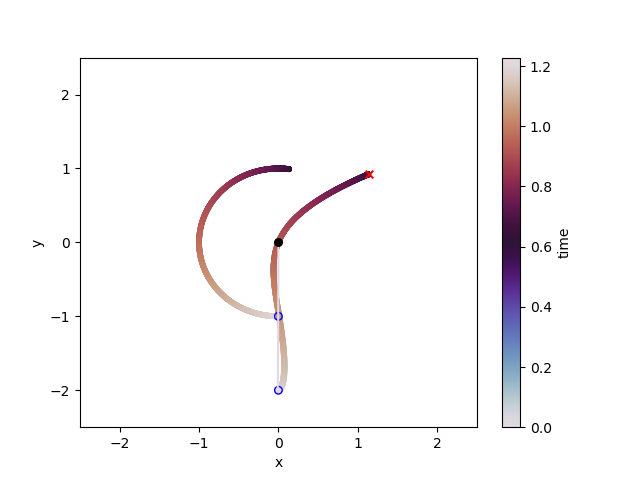}
  \caption{T=\SI{1.228}{\second}}
\end{subfigure}
\caption{Eigenmode at different time steps.}
\label{appD:conf1_traj}
\end{figure}

\begin{figure}[h!]
     \centering
     \begin{subfigure}[b]{0.33\textwidth}
         \centering
         \includegraphics[width=1.0\textwidth]{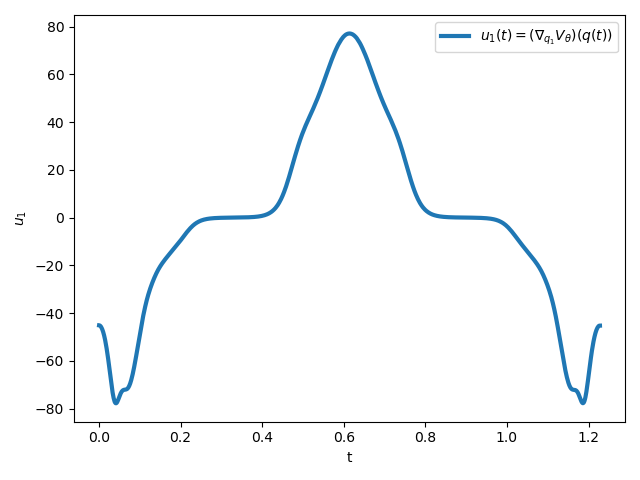}
         \caption{First component of $u(t)$.}
     \end{subfigure}
     \hfill
     \begin{subfigure}[b]{0.33\textwidth}
         \centering
         \includegraphics[width=1.0\textwidth]{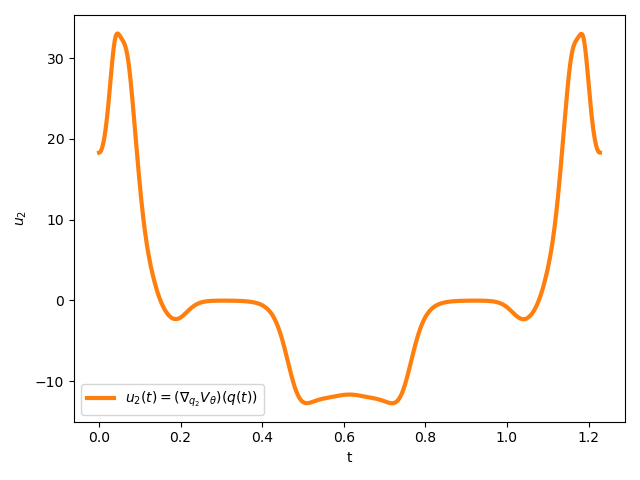}
         \caption{Second component of $u(t)$.}
     \end{subfigure}
     \hfill
     \begin{subfigure}[b]{0.33\textwidth}
        \centering
        \includegraphics[width=1.0\textwidth]{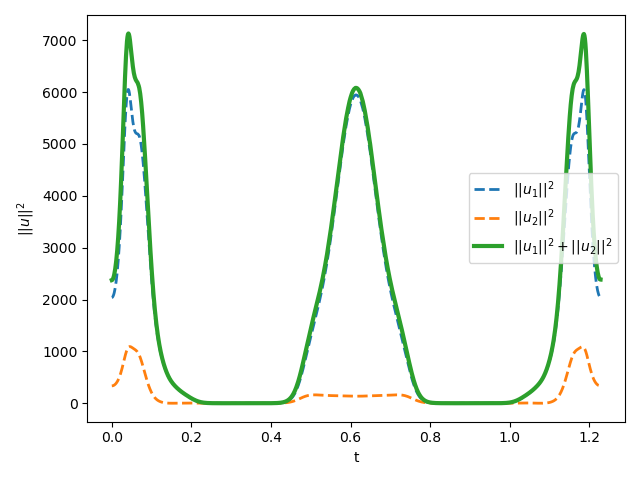}
         \caption{Squared control effort penalty.}
     \end{subfigure}
        \caption{Control inputs and control effort.}
        \label{appD:conf1_control_inp}
\end{figure}

\begin{figure}[h!]
     \centering
     \begin{subfigure}[b]{0.33\textwidth}
         \centering
         \includegraphics[width=1.0\textwidth]{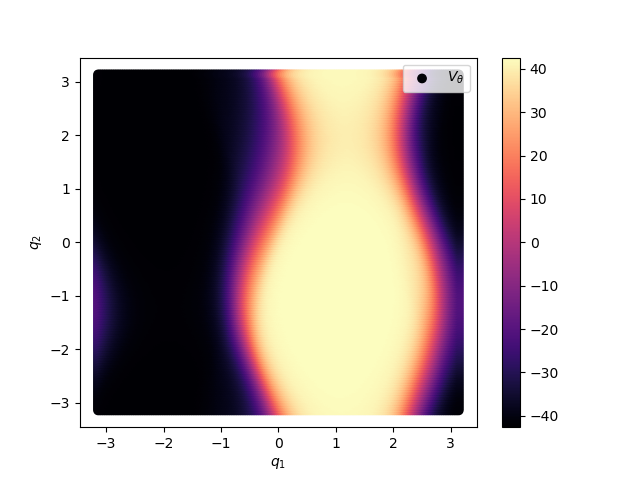}
         \caption{Learned potential $V_{\bm{\theta}}$.}
     \end{subfigure}
     \hfill
     \begin{subfigure}[b]{0.33\textwidth}
         \centering
          \includegraphics[width=1.0\textwidth]{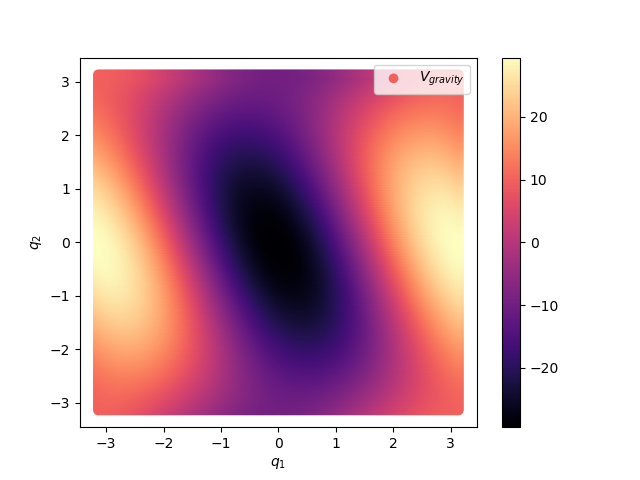}
         \caption{Gravitational potential $V_{\text{gravity}}$.}
     \end{subfigure}
     \hfill
     \begin{subfigure}[b]{0.33\textwidth}
         \centering
          \includegraphics[width=1.0\textwidth]{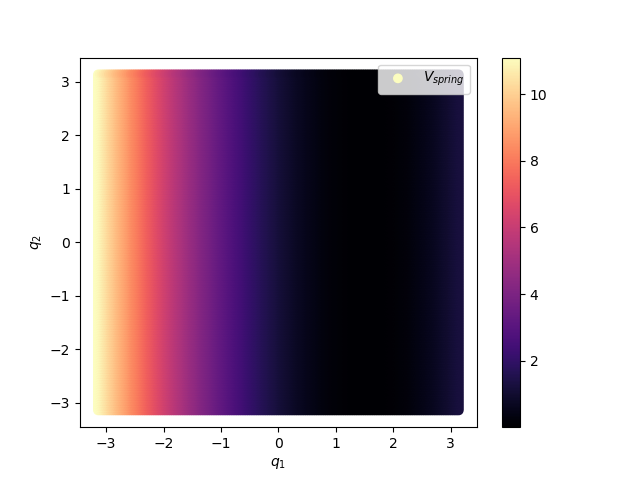}
         \caption{Spring potential $V_{spring}$}
     \end{subfigure}

     \begin{subfigure}[b]{0.33\textwidth}
         \centering
         \includegraphics[width=1.0\textwidth]{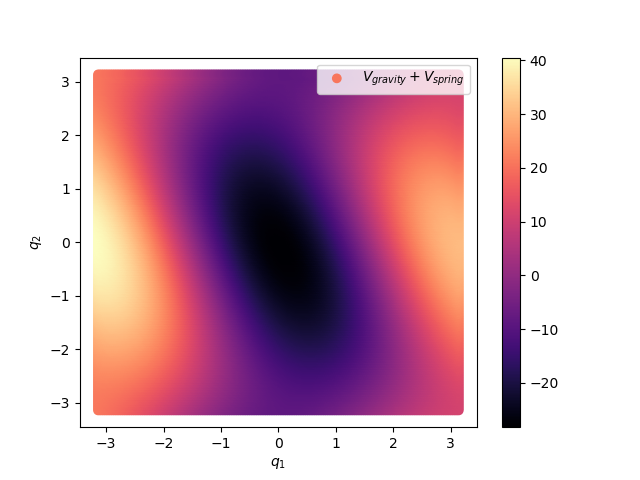}
         \caption{Open-loop potential $V_{\text{spring}}+V_{\text{gravity}}$.}
     \end{subfigure}
     \begin{subfigure}[b]{0.33\textwidth}
         \centering
          \includegraphics[width=1.0\textwidth]{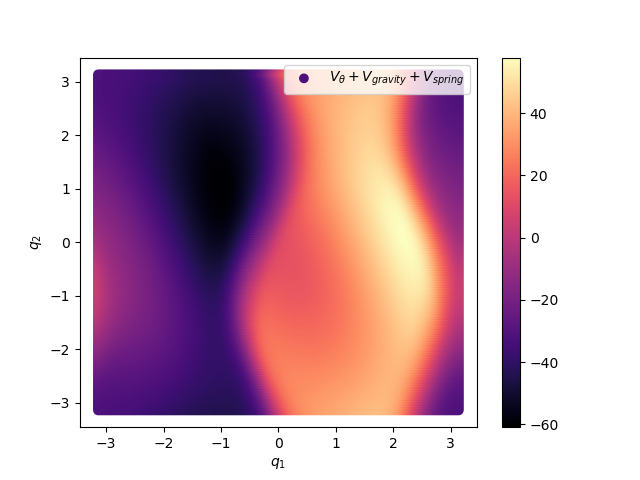}
         \caption{Overall potential $V_{\bm{\theta}}+V_{\text{spring}}+V_{\text{gravity}}$.}
     \end{subfigure}
        \caption{Potentials for T=\SI{1.228}{\second} over $\bm{q}\in[-\pi, \pi]$.}
        \label{appD:conf1_potentials}
\end{figure}

\begin{figure}[h!]
     \centering
     \begin{subfigure}[b]{0.33\textwidth}
         \centering
         \includegraphics[width=1.0\textwidth]{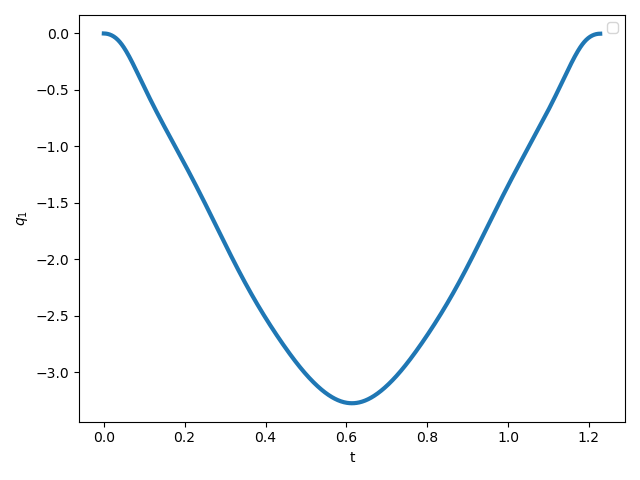}
         \caption{$q_1$ over one period.}
     \end{subfigure}
     \hfill
     \begin{subfigure}[b]{0.33\textwidth}
         \centering
         \includegraphics[width=1.0\textwidth]{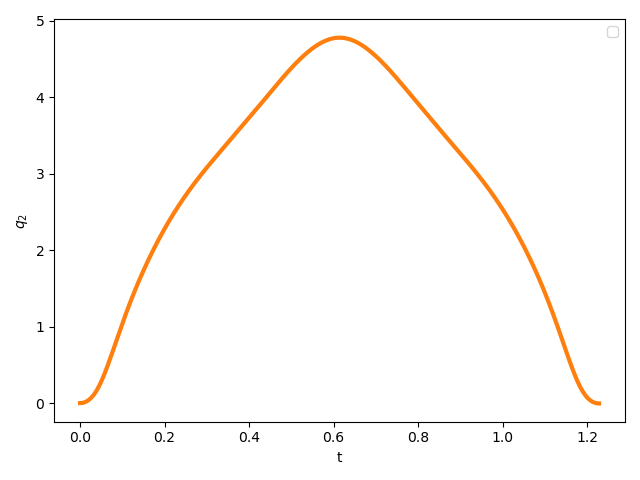}
         \caption{$q_2$ over one period.}
     \end{subfigure}
     \hfill
     \begin{subfigure}[b]{0.33\textwidth}
        \centering
        \includegraphics[width=1.0\textwidth]{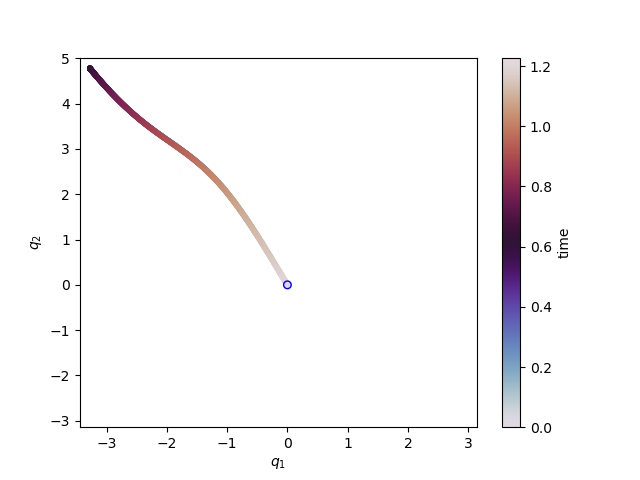}
         \caption{Trajectory in configuration space.}
     \end{subfigure}
        \caption{The time behavior of the angles $q_1$ and $q_2$ over one period.}
        \label{appD:conf1_q}
\end{figure}

\begin{figure}[h!]
\centering
\begin{subfigure}{.33\textwidth}
  \centering
  \includegraphics[width=1.0\linewidth]{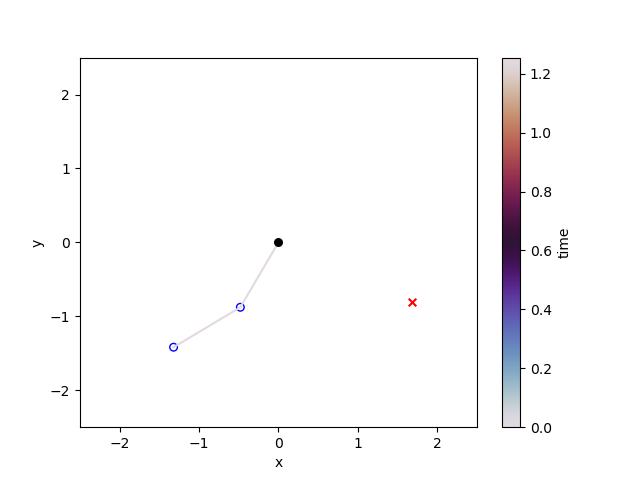}
  \caption{t=\SI{0.0}{\second}}
\end{subfigure}%
\begin{subfigure}{.33\textwidth}
  \centering
  \includegraphics[width=1.0\linewidth]{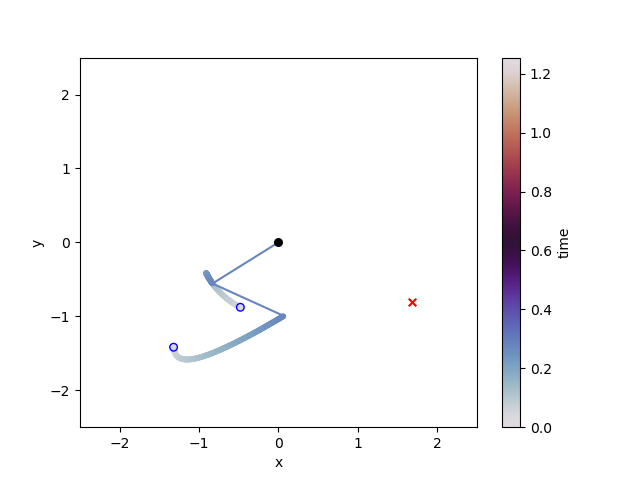}
  \caption{t=\SI{0.279}{\second}}
\end{subfigure}
\begin{subfigure}{.33\textwidth}
  \centering
  \includegraphics[width=1.0\linewidth]{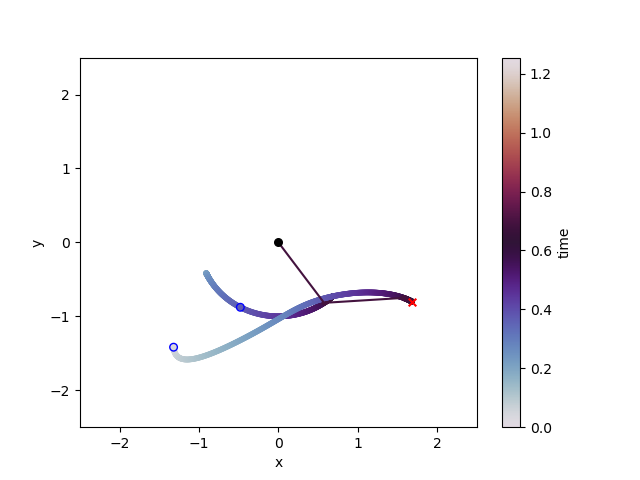}
  \caption{t=\SI{0.697}{\second}}
\end{subfigure}%
\hfill
\begin{subfigure}{.33\textwidth}
  \centering
  \includegraphics[width=1.0\linewidth]{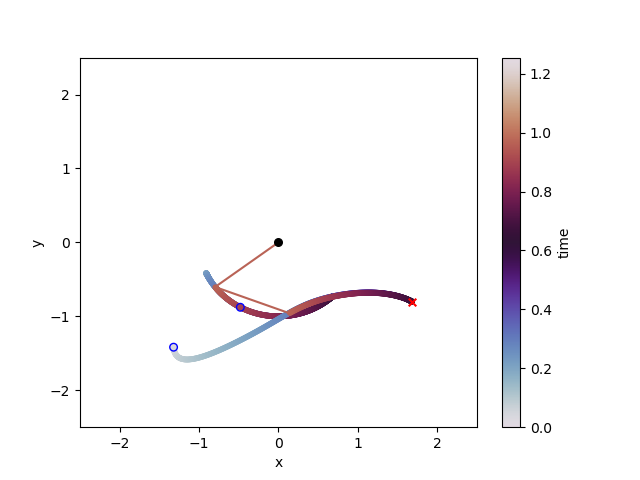}
  \caption{t=\SI{0.976}{\second}}
\end{subfigure}
\begin{subfigure}{.33\textwidth}
  \centering
  \includegraphics[width=1.0\linewidth]{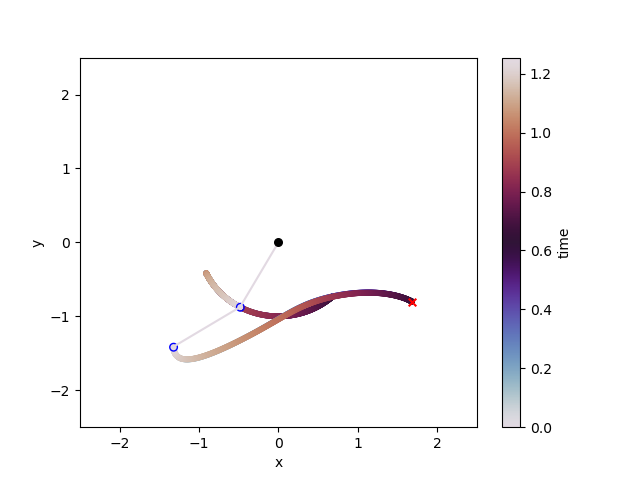}
  \caption{T=\SI{1.255}{\second}}
\end{subfigure}
\caption{Eigenmode at different time steps.}
\label{appD:conf2_traj}
\end{figure}

\begin{figure}[h!]
     \centering
     \begin{subfigure}[b]{0.33\textwidth}
         \centering
         \includegraphics[width=1.0\textwidth]{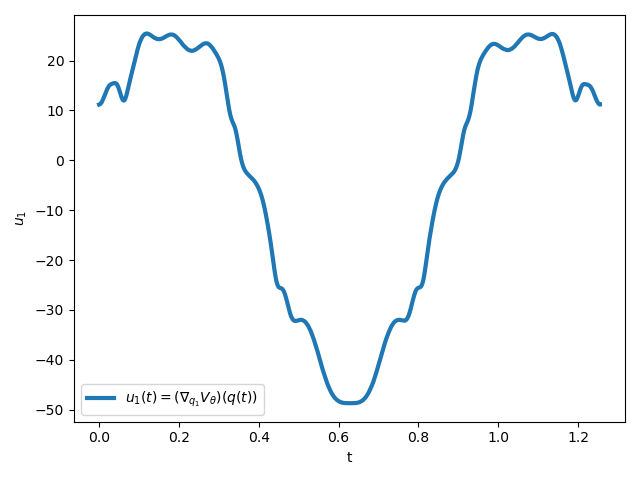}
         \caption{First component of $u(t)$.}
     \end{subfigure}
     \hfill
     \begin{subfigure}[b]{0.33\textwidth}
         \centering
         \includegraphics[width=1.0\textwidth]{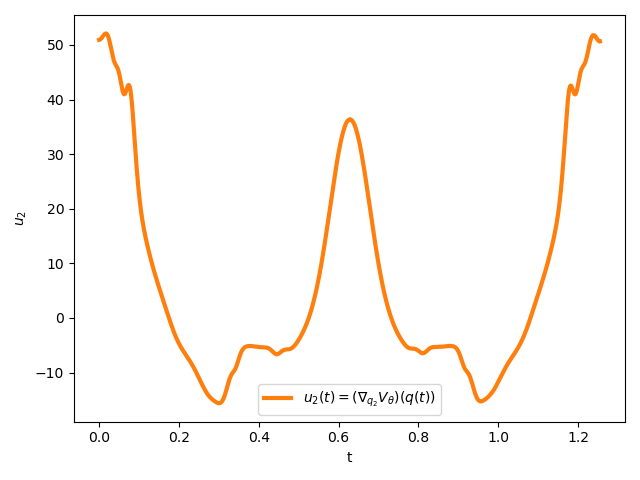}
         \caption{Second component of $u(t)$.}
     \end{subfigure}
     \hfill
     \begin{subfigure}[b]{0.33\textwidth}
        \centering
        \includegraphics[width=1.0\textwidth]{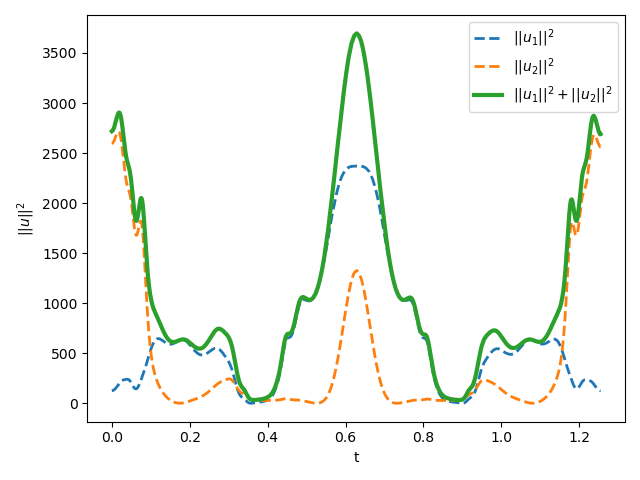}
         \caption{Squared control effort penalty.}
     \end{subfigure}
        \caption{Control inputs and control effort.}
        \label{appD:conf2_control_inp}
\end{figure}

\begin{figure}[h!]
     \centering
     \begin{subfigure}[b]{0.33\textwidth}
         \centering
         \includegraphics[width=1.0\textwidth]{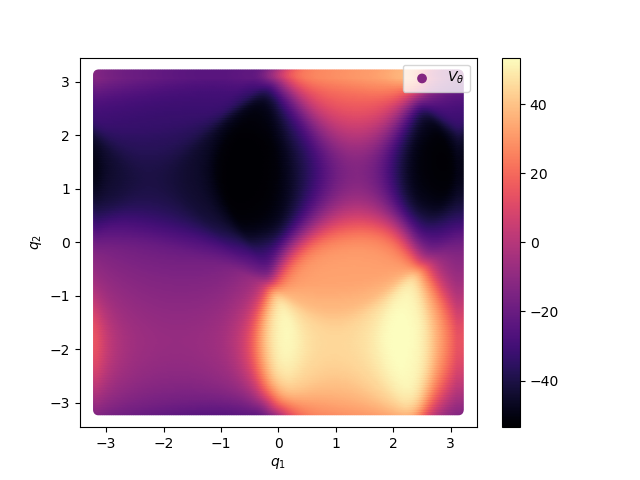}
         \caption{Learned potential $V_{\bm{\theta}}$.}
     \end{subfigure}
     \hfill
     \begin{subfigure}[b]{0.33\textwidth}
         \centering
          \includegraphics[width=1.0\textwidth]{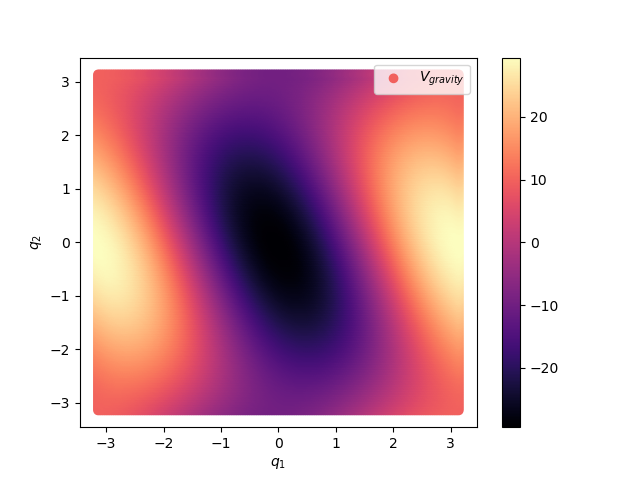}
         \caption{Gravitational potential $V_{\text{gravity}}$.}
     \end{subfigure}
     \hfill
     \begin{subfigure}[b]{0.33\textwidth}
         \centering
          \includegraphics[width=1.0\textwidth]{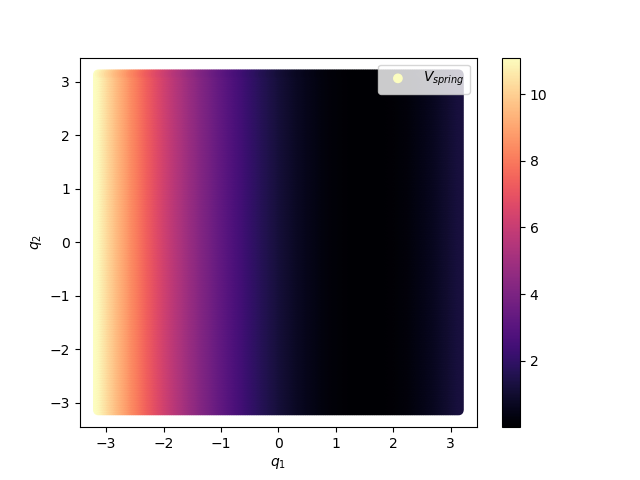}
         \caption{Spring potential $V_{spring}$}
     \end{subfigure}

     \begin{subfigure}[b]{0.33\textwidth}
         \centering
         \includegraphics[width=1.0\textwidth]{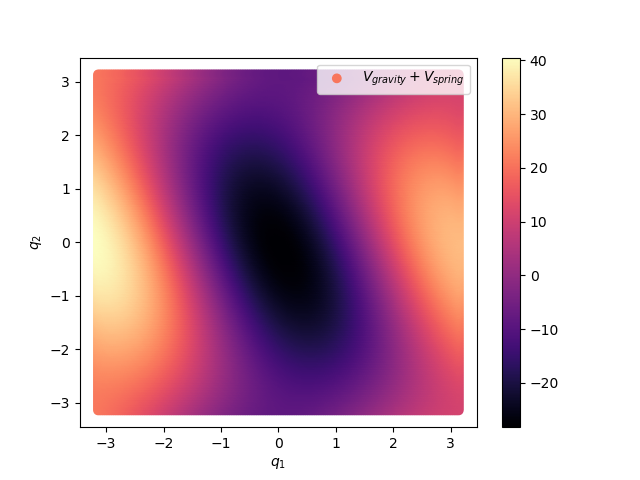}
         \caption{Open-loop potential 
$V_{\text{spring}}+V_{\text{gravity}}$.}
     \end{subfigure}
     \begin{subfigure}[b]{0.33\textwidth}
         \centering
          \includegraphics[width=1.0\textwidth]{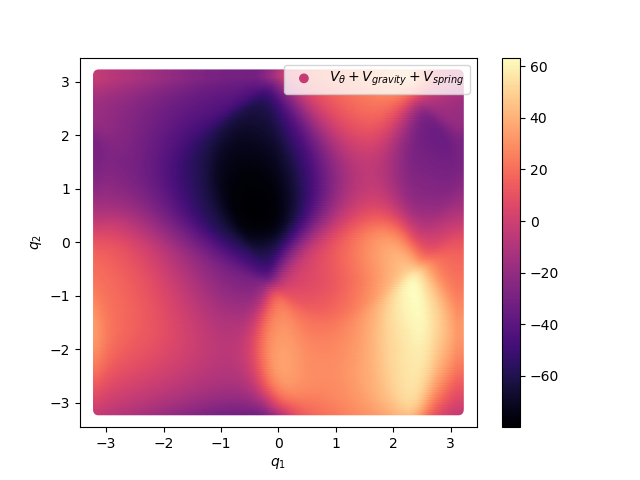}
         \caption{Overall potential $V_{\bm{\bm{\theta}}}+V_{\text{spring}}+V_{\text{gravity}}$.}
     \end{subfigure}
        \caption{Potentials for T=\SI{1.255}{\second} over $\bm{q}\in[-\pi, \pi]$.}
        \label{appD:conf2_potentials}
\end{figure}

\begin{figure}[h!]
     \centering
     \begin{subfigure}[b]{0.33\textwidth}
         \centering
         \includegraphics[width=1.0\textwidth]{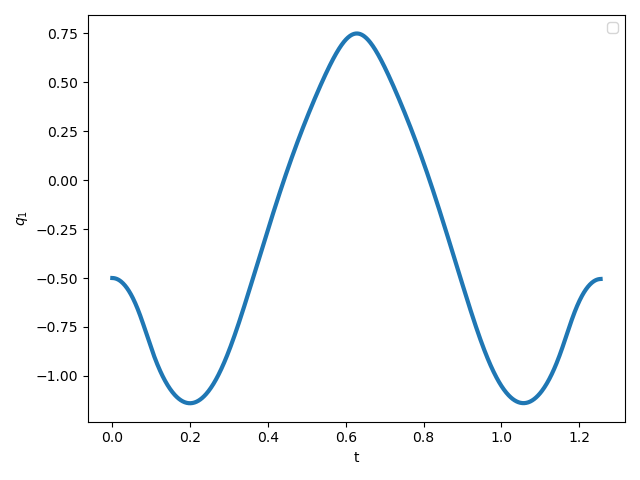}
         \caption{$q_1$ over one period.}
     \end{subfigure}
     \hfill
     \begin{subfigure}[b]{0.33\textwidth}
         \centering
         \includegraphics[width=1.0\textwidth]{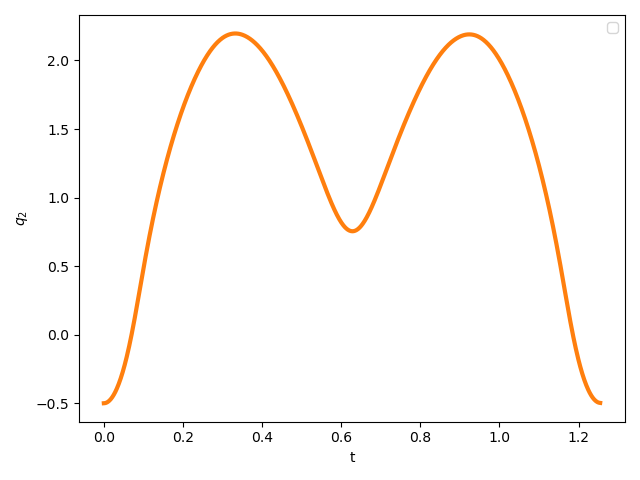}
         \caption{$q_2$ over one period.}
     \end{subfigure}
     \hfill
     \begin{subfigure}[b]{0.33\textwidth}
        \centering
        \includegraphics[width=1.0\textwidth]{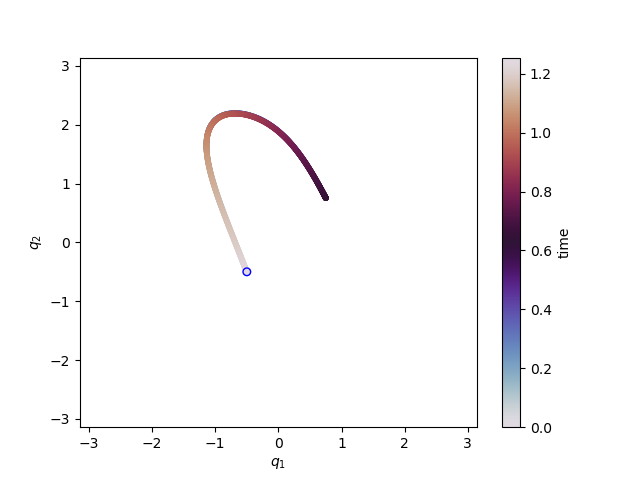}
         \caption{Trajectoryin configuration space.}
     \end{subfigure}
        \caption{The time behavior of the angles $q_1$ and $q_2$ over one period.}
        \label{appD:conf2_q}
\end{figure}

Similarly to Appendix \ref{AppC:additional_results_config1}, we show the trajectories of the double pendulum in Figure \ref{appD:conf1_traj} and \ref{appD:conf2_traj}, the control inputs in Figure \ref{appD:conf1_control_inp} and \ref{appD:conf2_control_inp}, the potentials in Figure \ref{appD:conf1_potentials} and \ref{appD:conf2_potentials}  and the state variable over time in Figure \ref{appD:conf1_q} and \ref{appD:conf2_q}, for two different initial and final positions.


\clearpage


\clearpage
\section{Implementation Details}\label{AppE:Implementation}

The neural ODE was implemented using the TorchDyn library \cite{PoliTorchDyn:PyTorch} built on top of PyTorch \cite{Paszke2019PyTorch:Library}.

The neural network $V_{\bm{\theta}}$ is composed of an input layer of dimension 2 (in the case of the double pendulum, the two input variables are $q_1$ and $q_2$, one hidden layers with 256 neurons, and an output layer of dimension 1 outputting the value of the potential $V_{\bm{\theta}}(q_1,q_2)$ for the input pair $(q_1,q_2)$. The choices of hyperparameters are shown in table \ref{tab_app:hyper_list}.

\begin{table}[h!]
\centering
\begin{tabular}{||c | c ||} 
 \hline
 Hyperparameter & Value \\ [0.5ex] 
 \hline\hline
 $V_{\bm{\theta}}$ input dimension & 2 \\ 
 \hline
 $V_{\bm{\theta}}$ output dimension & 1 \\ 
 \hline
 $V_{\bm{\theta}}$ hidden dimension & 256 \\ 
 \hline
Number of hidden layers & 1 \\ 
 \hline
Activation function hidden layer & tanh \\ 
 \hline
Activation function output layer & linear \\ 
 \hline
Optimizer & ADAM \cite{KingmaADAM:OPTIMIZATION} \\ 
 \hline
Learning rate & 1e-3 \\ 
 \hline
 Training epochs & 500 \\ 
 \hline
Computation of sensitivity & backpropagation \\ 
 \hline
  $\alpha_{\text{task}}$ & 10\\
 \hline
 $\alpha_{\text{eff}}$ & 0.0001 \\
 \hline
 $\alpha_{\text{task}}$ & 10\\
 \hline
$\lambda_1$ & 0.05\\
 \hline
  $\alpha_1$ & 0.0005\\
 \hline
 $\lambda_2$ & 0.95\\
  \hline 
 $\beta$ & 1\\
 \hline
  $\alpha_M$ & 10\\
 \hline
  $\alpha_E$ & 1\\
 \hline
Double pendulum mass link 1 & 1.0 \\ 
 \hline
Double pendulum mass link 2 & 1.0 \\ 
 \hline
Double pendulum length link 1 & 1.0 \\ 
 \hline
Double pendulum length link 2 & 1.0 \\ 
 \hline
Double pendulum spring stiffness joint 1 & 0.0 \\ 
 \hline
Double pendulum spring stiffness joint 2 & 0.5 \\  
 \hline
\end{tabular}
\caption{Hyperparameters of the experiments.}
\label{tab_app:hyper_list}
\end{table}

\end{document}